\documentclass[11pt,journal,onecolumn]{IEEEtran}
\usepackage[dvips]{graphicx}
\usepackage{amsmath,amssymb,amsthm,float,arydshln,cite,color}
\usepackage{psfrag,epsfig,setspace,wrapfig}
\usepackage[latin1]{inputenc}
\usepackage{dsfont}
\usepackage{subfigure}
\usepackage{colortbl}

\makeatletter
\def\hlinewd#1{%
  \noalign{\ifnum0=`}\fi\hrule \@height #1 \futurelet
   \reserved@a\@xhline}
\makeatother
\newsavebox{\savepar}
\newenvironment{boxit}{\begin{lrbox}{\savepar}
\begin{minipage}[b]{6in}}
{\end{minipage}\end{lrbox}\fbox{\usebox{\savepar}}}

\newcommand{\bx}{\fboxrule0.015cm\begin{boxit}\;}
\newcommand{\ex}{\;\end{boxit}\fboxrule0.1cm}
\def\wt{\widetilde}
\def\wh{\widehat}
\def\ol{\overline}
\def\col{{\mathrm{col}}}

\def\muub{\mu_{\mathrm{ub}}}

\def\bmu     {{\bm \mu}}
\def\bpi     {{\bm \pi}}
\def\bgamma  {{\bm \gamma}}
\def\bpsi    {{\bm \psi}}

\newcommand{\weq}{{\;=\;}}

\newcommand{\Ncal}{\mathcal{N}}

\newcommand{\Xcal}{\mathcal{X}}

\def\kron{\otimes}

\def\MSD{{\textrm{MSD}}}
\def\EMSE{{\textrm{EMSE}}}
\def\ER{{\textrm{ER}}}

\def\Filtration{\bm{\mathcal F}}

\def\cent{{\textrm{cent}}}
\def\dist{{\textrm{dist}}}

\def\lms{{\textrm{LMS}}}

\def\glob{\textrm{glob}}

\def\ind{{\textrm{ncop}}}
\def\av{{\textrm{av}}}
\def\sync{{\textrm{sync}}}
\def\async{{\textrm{asyn}}}

\def\lambdamax{\lambda_{\mathrm{max}}}
\def\lambdamin{\lambda_{\mathrm{min}}}

\def\tran{^{\mathsf{T}}}
\def\one{\mathds{1}}

\newcommand{\bp}{\begin{proof} \small }
\newcommand{\ep}{\end{proof} \normalsize}
\newcommand{\epx}{\end{proof} \small}
\newcommand{\bpa}{\begin{proofappx} \footnotesize }
\newcommand{\epa}{\end{proofappx} \small }
\newtheorem{example}{{\bf Example}}[section]
\newtheorem{lemma}{{\bf Lemma}}[section]
\newtheorem{definition}{{\bf Definition}}[section]
\newtheorem{assumption}{{\bf Assumption}}[section]

\newcommand{\addbox}{\addtolength{\fboxsep}{5pt}\boxed}
\newcommand{\Ex}{\mathbb{E}\hspace{0.05cm}}

\newcommand{\bm}[1]{\mbox{\boldmath $#1$}}

\newcommand{\be}{\begin{equation}}
\newcommand{\ee}{\end{equation}}
\newcommand{\bs}{\begin{subequations}}
\newcommand{\es}{\end{subequations}}
\newcommand{\bq}{\begin{eqnarray}}
\newcommand{\eq}{\end{eqnarray}}
\newcommand{\bqn}{\begin{eqnarray*}}
\newcommand{\eqn}{\end{eqnarray*}}
\newcommand{\nn}{\nonumber}
\newcommand{\ba}{\left[ \begin{array}}
\newcommand{\ea}{\\ \end{array} \right]}
\newcommand{\ben}{\begin{enumerate}}
\newcommand{\een}{\end{enumerate}}
\newcommand{\qd}{\hfill{$\blacklozenge$}}
\newcommand{\define}{\;\stackrel{\Delta}{=}\;}

\newcommand{\Tr}{\mbox{\rm {\small Tr}}}

\def\Nknk{\Ncal_k\backslash\{k\}}

\def\A{{\boldsymbol{A}}}

\def\H{{\boldsymbol{H}}}

\def\X{{\boldsymbol{X}}}

\def\a{{\boldsymbol{a}}}

\def\c{{\boldsymbol{c}}}
\def\d{{\boldsymbol{d}}}

\def\h{{\boldsymbol{h}}}

\def\s{{\boldsymbol{s}}}

\def\u{{\boldsymbol{u}}}
\def\v{{\boldsymbol{v}}}
\def\w{{\boldsymbol{w}}}
\def\x{{\boldsymbol{x}}}

\def\z{{\boldsymbol{z}}}

\def\real{{\mathchoice%
{\hbox{\rm\setbox1=\hbox{I}\copy1\kern-.45\wd1 R}}
{\hbox{\rm\setbox1=\hbox{I}\copy1\kern-.45\wd1 R}}
{\hbox{\scriptsize\rm\setbox1=\hbox{I}\copy1\kern-.45\wd1 R}}
{\hbox{\scriptsize\rm\setbox1=\hbox{I}\copy1\kern-.45\wd1 R}}}}

\def\Zint{{\mathchoice{\setbox1=\hbox{\sf Z}\copy1\kern-.75\wd1\box1}
{\setbox1=\hbox{\sf Z}\copy1\kern-.75\wd1\box1}
{\setbox1=\hbox{\scriptsize\sf Z}\copy1\kern-.75\wd1\box1}
{\setbox1=\hbox{\scriptsize\sf Z}\copy1\kern-.75\wd1\box1}}}
\newcommand{\complex}{ \hbox{\rm C\kern-0.45em\rule[.07em]{.02em}{.58em}%
\kern 0.43em}}

\begin{document}
\title{Asynchronous Adaptive Networks\thanks{{\bf To appear as book chapter in {\em Cooperative and Graph Signal Processing}, P. Djuric and C. Richard, editors, Elsevier, 2018}. The work in this book chapter was supported in part by NSF grants ECCS-1407712 and CCF-1524250. The authors are grateful to IEEE for allowing reproduction of substantial material from \cite{ProcIEEE2014,xiao2013,xiao2013AA,xiao2013BB} in this book chapter. Co-author X. Zhao was a PhD student in Electrical Engineering at UCLA. The first author is with the Ecole Polytechnique Federale de Lausanne, EPFL, CH-1015 Lausanne. Contact email: ali.sayed@epfl.ch. }
}

\author{\authorblockN{Ali H. Sayed\,\, and\,\, Xiaochuan Zhao}
%\\\vspace{0.2cm}
%Electrical Engineering Department\\
%University of California, Los Angeles\\\vspace{0.5cm}
%{\footnotesize \noindent Draft for a book chapter based on the publications \cite{ProcIEEE2014,xiao2013,xiao2013AA,xiao2013BB}}
}

\maketitle

\noindent

\begin{abstract}
The overview article \cite{ProcIEEE2014} surveyed advances related to adaptation, learning, and optimization over synchronous networks. Various distributed strategies were discussed that enable a collection of networked agents to interact locally in response to streaming data and to continually learn and adapt to track drifts in the data and models. Under reasonable technical conditions on the data, the adaptive networks were shown to be mean-square stable in the slow adaptation regime, and their mean-square-error performance and convergence rate were characterized in terms of the network topology and data statistical moments \cite{NOW2014}. Classical results for single-agent adaptation and learning were recovered as special cases. Following the works \cite{xiao2013,xiao2013AA,xiao2013BB}, this chapter complements the exposition from \cite{ProcIEEE2014} and extends the results to asynchronous networks where agents are subject to various sources of uncertainties that influence their behavior, including randomly changing topologies, random link failures, random data arrival times, and agents turning on and off randomly. In an asynchronous environment, agents may stop updating their solutions or may stop sending or receiving information in a random manner and without coordination with other agents. The presentation will reveal that the mean-square-error performance of asynchronous networks remains largely unaltered compared to synchronous networks. The results justify the remarkable resilience of cooperative networks.
\end{abstract}

\small

\section{Introduction}
Adaptive networks consist of a collection of agents with learning abilities. The agents interact with each other on a local level and diffuse information across the network to solve inference and optimization tasks in a decentralized manner. Such networks are scalable, robust to node and link failures, and are particularly suitable for learning from big data sets by tapping into the power of collaboration among distributed agents. The networks are also endowed with cognitive abilities due to the sensing abilities of their agents, their interactions with their neighbors, and the embedded feedback mechanisms for acquiring and refining information. Each agent is not only capable of sensing data and experiencing the environment directly, but it also receives information through interactions with its neighbors and processes and analyzes this information to drive its learning process.

As already indicated in \cite{ProcIEEE2014,NOW2014}, there are many good reasons for the peaked interest in networked solutions, especially in this day and age when the word ``network'' has become commonplace whether one is referring to social networks, power networks, transportation networks, biological networks, or other networks. Some of these reasons have to do with the benefits of cooperation over networks in terms of improved performance and improved robustness and resilience to failure. Other reasons deal with privacy and secrecy considerations where agents may not be comfortable sharing their data with remote fusion centers. In other situations, the data may already be available in dispersed locations, as happens with cloud computing. One may also be interested in learning and extracting information through data mining from large data sets. Decentralized learning procedures offer an attractive approach to dealing with such data sets. Decentralized mechanisms can also serve as important enablers for the design of robotic swarms, which can assist in the exploration of disaster areas.

\subsection{Asynchronous Behavior}
The survey article \cite{ProcIEEE2014} and monograph \cite{NOW2014} focused on the case of synchronous networks where data arrive at all agents in a synchronous manner and updates by the agents are also performed in a synchronous manner. The network topology was assumed to remain largely static during the adaptation process. Under these conditions, the limits of performance and stability of these networks were identified in some detail for two main classes of distributed strategies: consensus and diffusion constructions. In this chapter, we extend the overview from \cite{ProcIEEE2014} to cover {\em asynchronous} environments. In such environments, the operation of the network can suffer from the occurrence of various random events including randomly changing topologies, random link failures, random data arrival times, and agents turning on and off randomly. Agents may also stop updating their solutions or may stop sending or receiving information in a random manner and without coordination with other agents. Results in \cite{xiao2013,xiao2013AA,xiao2013BB} examined the implications of such  asynchronous events on network performance in some detail and under a fairly general model for the random events. The purpose of this chapter is to summarize the key conclusions from these works in a manner that complements the presentation from \cite{ProcIEEE2014} for the benefit of the reader. While the works \cite{xiao2013,xiao2013AA,xiao2013BB} consider a broader formulation involving complex-valued variables, we limit the discussion here to real-valued variables in order not to overload the notation and to convey the key insights more directly. Proofs and derivations are often omitted and can be found in the above references; the emphasis is on presenting the results in a motivated manner and on commenting on the insights they provide into the operation of asynchronous networks.

We indicated in \cite{xiao2013,xiao2013AA,xiao2013BB} that there already exist many useful studies in the literature on the performance of consensus strategies in the presence of asynchronous events (see, e.g., \cite{TsiBB86}--\cite{Aysal09Allerton}. There are also some studies in the context of diffusion strategies \cite{Lopes2008aax, tak2010}. However, with the exception of the latter two works, the earlier references assumed conditions that are not generally favorable for applications involving continuous {\em adaptation} and learning. For example, some of the works assumed a decaying step-size, which turns off adaptation after sufficient iterations have passed. Some other works assumed noise free data, which is a hindrance when learning from data perturbed by interferences and distortions. A third class of works focused on studying pure averaging algorithms, which are not required to respond to continuous data streaming. In the works \cite{xiao2013,xiao2013AA,xiao2013BB}, we adopted a more general asynchronous model that removes these limitations by allowing for various sources of random events and, moreover, the events are allowed to occur simultaneously. We also examined learning algorithms that respond to streaming data to enable adaptation. The main conclusion from the analysis in these works, and which will be summarized in future sections, is that asynchronous networks can still behave in a mean-square-error stable manner for sufficiently small step-sizes and, interestingly, their steady-state performance level is only slightly affected in comparison to synchronous behavior. The iterates computed by the various agents are still able to converge and hover around an agreement state with a small mean-square-error. These are reassuring results that support the intrinsic robustness and resilience of network-based cooperative solutions.

\subsection{Organization of the Chapter}
Readers will benefit more from this chapter if they review first the earlier article \cite{ProcIEEE2014}. We continue to follow a similar structure here, as well as a similar notation, since the material in both this chapter and the earlier reference \cite{ProcIEEE2014} are meant to complement each other. We organize the presentation into three main components. The first part (Sec.~\ref{sec:single}) reviews fundamental results on adaptation and learning by {\em single} stand-alone agents. The second part (Sec.~\ref{sec:batch}) covers asynchronous centralized solutions. The objective is to explain the gain in performance that results from aggregating the data from the agents and processing it centrally at a fusion center. The centralized performance is used as a frame of reference for assessing various implementations. While centralized solutions can be powerful, they nevertheless suffer from a number of limitations. First, in real-time applications where agents collect data continuously, the repeated exchange of information back and forth between the agents and the fusion center can be costly especially when these exchanges occur over wireless links or require nontrivial routing resources. Second, in some sensitive applications, agents may be reluctant to share their data with remote centers for various reasons including privacy and secrecy considerations. More importantly perhaps, centralized solutions have a critical point of failure: if the central processor fails, then this solution method collapses altogether.

For these reasons, we cover in the remaining sections of the chapter (Secs.~\ref{sec:network} and \ref{sec:asynchronous}) distributed asynchronous strategies of the consensus and diffusion types, and examine their dynamics, stability, and performance metrics. In the distributed mode of operation, agents are connected by a topology and they are permitted to share information only with their immediate neighbors.  The study of the behavior of such networked agents is more challenging than in the single-agent and centralized modes of operation due to the coupling among interacting agents and due to the fact that the networks are generally sparsely connected.

\section{Single-Agent Adaptation and Learning}
\label{sec:single}
We begin our treatment by reviewing stochastic gradient algorithms, with emphasis on their application to the problems of adaptation and learning by stand-alone agents.

\subsection{Risk and Loss Functions}
Thus, let $J(w): \real^{M\times 1} \mapsto \real$ denote a twice-differentiable real-valued (cost or utility or risk) function of a real-valued vector argument, $w\in\real^{M\times 1}$. When the variable $w$ is complex-valued, some important technical differences arise that  are beyond the scope of this chapter; they are addressed in  \cite{xiao2013,xiao2013AA,xiao2013BB,NOW2014}. Likewise, some adjustments to the arguments are needed when the risk function is non-smooth (non-differentiable); as explained in the works \cite{Ying2017a,Ying2017b}. It is sufficient for our purposes in this chapter to limit the presentation to real arguments and smooth risk functions without much loss in generality.

We denote the gradient vectors of $J(w)$ relative to $w$ and $w\tran$ by the following row and column vectors, respectively:
\bs
\begin{align}
\nabla_w \, J(w) & \define \left[
\dfrac{\partial J(w)}{\partial w_1},
\dfrac{\partial J(w)}{\partial w_2},
\dots,
\dfrac{\partial J(w)}{\partial w_M} \right] \\
\nabla_{w\tran} \, J(w) & \define \left[ \nabla_w \, J(w) \right]\tran
\end{align}
These definitions are in terms of the partial derivatives of $J(w)$ relative to the individual entries of $ w = \col\{w_1, w_2, \dots, w_M \}$, where the notation $\col\{\cdot\}$ refers to a column vector that is formed by stacking its arguments  on top of each other. Likewise, the Hessian matrix of $J(w)$ with respect to $w$ is defined as the following  $M\times M$  symmetric matrix:
\be
\nabla_w^2 \, J(w) \define \nabla_{w\tran}[\nabla_w \, J(w)] \weq \nabla_{w}[\nabla_{w\tran} \, J(w)]
\ee
\es
which is constructed from two successive gradient operations. It is common in adaptation and learning applications for the risk function $J(w)$ to be constructed as the expectation of some loss function, $Q(w;\x)$, where the {\bf boldface} variable $\x$ is used to denote some random data, say,
\be
\label{eqn:losscostdef}
J(w) \weq \Ex\,Q(w;\x)
\ee
and the expectation is evaluated over the distribution of $\x$.

{\small
\noindent \begin{example}[{\bf Mean-square-error (MSE) costs}]
\label{example:MSEcosts}{\rm
Let $\d$ denote a zero-mean scalar random variable with variance $\sigma_{d}^2 = \Ex\d^2$, and let $\u$ denote a zero-mean $1\times M$ random vector with covariance matrix $R_{u} = \Ex\u\tran\u > 0$. The combined quantities $\{\d,\u\}$ represent the random variable $\x$ referred to in \eqref{eqn:losscostdef}. The cross-covariance vector is denoted by $r_{du}=\Ex\d\u\tran$. We formulate the problem of estimating $\d$ from $\u$ in the linear least-mean-squares-error sense or, equivalently, the problem of seeking the vector $w^o$ that minimizes the quadratic cost function:
\bs
\be
\label{eqn:Jmsecostdef}
J(w) \define \Ex(\d-\u w)^2 = \sigma_{d}^2 - r_{du}\tran w - w\tran r_{du} + w\tran R_{u} w
\ee
This cost corresponds to the following choice for the loss function:
\be
Q(w;\x) \weq (\d-\u w)^2
\ee
Such quadratic costs are widely used in estimation and adaptation problems \cite{Hay02}--\cite{Kai00}. They are also widely used as quadratic risk functions in machine learning applications \cite{Bish2007,theo2008}. The gradient vector and Hessian matrix of $J(w)$ are easily seen to be:
\be
\label{eqn:JgradientHessiandef}
\nabla_w \, J(w) \weq 2 \left( R_{u} w - r_{du} \right)\tran, \;\;\;\; \nabla_w^2\,J(w) \weq 2 R_{u}
\ee
\es
\qd
}

\label{example-1}
\end{example}
}

{\small
\begin{example}[{\bf Logistic or log-loss risks}]
\label{example:loglossrisk}{\rm
Let $\bgamma$ denote a binary random variable that assumes the values $\pm 1$, and let $\h$ denote an $M\times 1$ random (feature) vector with $R_h \weq \Ex \h \h\tran$.  The combined quantities $\{ \bgamma,\h \}$ represent the random variable $\x$ referred to in \eqref{eqn:losscostdef}.  In the context of machine learning and pattern classification problems \cite{Bish2007,theo2008,hos2000}, the variable $ \bgamma $ designates the class that feature vector $\h$ belongs to. In these problems, one seeks the vector $w^o$ that minimizes the regularized logistic risk function:
\bs
\be
\label{eqn:loglossdef}
J(w) \define \dfrac{\rho}{2} \|w\|^2 + \Ex \left\{ \ln \left[ 1 + e^{-\bgamma \h\tran w} \right] \right\}
\ee
where $\rho > 0$ is some regularization parameter, $\ln(\cdot)$ is the natural logarithm function, and $\|w\|^2 \weq w\tran w$. The risk \eqref{eqn:loglossdef} corresponds to the following choice for the loss function:
\be
\label{eqn:logriskdef}
Q(w;\x) \define \dfrac{\rho}{2} \|w\|^2 + \ln \left[ 1 + e^{- \bgamma \h\tran w} \right]
\ee
Once $w^o$ is recovered, its value can be used to classify new feature vectors, say, $\{ \h_{\ell} \}$, into classes $+1$ or $-1$. This can be achieved by assigning feature vectors with $\h_{\ell}\tran w^o \geq 0$ to one class and feature vectors with $\h_{\ell}\tran w^o < 0$ to another class. It can be easily verified that:
\begin{align}
\label{eqn:loglossgradientdef}
\nabla_w \, J(w) & \weq \rho w\tran - \Ex \left\{ \bgamma \h\tran \cdot \dfrac{e^{-\bgamma \h\tran w}}
{1 + e^{-\bgamma \h\tran w}} \right\} \\
\label{eqn:loglossHessiandef}
\nabla_w^2 \, J(w) & \weq \rho I_M + \Ex \left\{ \h\h\tran \cdot \dfrac{e^{-\bgamma \h\tran w}}
{ \left( 1 + e^{- \bgamma \h\tran w} \right)^2 }\right\}
\end{align}
where $I_{M}$ denotes the identity matrix of size $M\times M$.
\es
}

\qd
\label{example-2}
\end{example}
}

\subsection{Conditions on Cost Function}
Stochastic gradient algorithms are powerful iterative procedures for solving optimization problems of the form
\be
\label{eqn:minimizeproblem}
\min_{w} \;\; J(w)
\ee
While the analysis that follows can be pursued under more relaxed conditions (see, e.g., the treatments in \cite{Pol73}--\cite{tsypkin1971}), it is sufficient for our purposes to require $J(w)$ to be strongly-convex and twice-differentiable with respect to $w$. The cost function $J(w)$ is said to be $\nu$-strongly convex if, and only if, its Hessian matrix is sufficiently bounded away from zero \cite{Poljak87,boyd2004convex,dimitri2003,nesterov2004}:
%\bs
\be
\label{eqn:nustronglyconvexdef}
J(w) \; \mbox{\rm is $\nu$-strongly convex} \Longleftrightarrow \nabla^2_w\,J(w) \geq \nu I_{M} > 0
\ee
for all $w$ and for some scalar $\nu > 0$, where the notation $A > 0$ signifies that matrix $A$ is positive-definite. Strong convexity is a useful condition in the context of adaptation and learning from streaming data because it helps guard against ill-conditioning in the algorithms; it also helps ensure that $J(w)$ has a {\em unique} global minimum, say, at location $w^o$; there will be no other minima, maxima, or saddle points. In addition, it is well-known that strong convexity helps endow stochastic-gradient algorithms with geometric convergence rates in the order of $O(\alpha^i)$, for some $0 \leq \alpha < 1$ and where $i$ is the iteration index \cite{Poljak87,Ber97}.

In many problems of interest in adaptation and learning, the cost function $J(w)$ is either already strongly-convex or can be made strongly-convex by means of regularization. For example, it is common in machine learning problems \cite{Bish2007,theo2008} and in adaptation and estimation problems \cite{Sayed08,Kai00} to incorporate regularization factors into the cost functions; these factors help ensure strong convexity. For instance, the mean-square-error cost \eqref{eqn:Jmsecostdef} is strongly convex whenever $R_{u}>0$. If $R_{u}$ happens to be singular, then the following regularized cost will be strongly convex:
\be
J(w) \define \dfrac{\rho}{2} \|w\|^2 + \Ex( \d - \u w)^2
\ee
where $\rho > 0$ is a regularization parameter similar to \eqref{eqn:loglossdef}.

Besides strong convexity, we shall also assume that the gradient vector of $J(w)$ is $\delta$-Lipschitz, namely, there exists $\delta > 0$ such that
\be
\label{eqn:deltaLipschitzdef}
\| \nabla_w\,J(w_2) - \nabla_w\,J(w_1) \| \leq \delta\,\| w_2 - w_1 \|
\ee
for all $w_1$ and $w_2$. It can be verified that for twice-differentiable costs, conditions \eqref{eqn:nustronglyconvexdef} and \eqref{eqn:deltaLipschitzdef} combined are equivalent to
\be
\label{eqn:boundedHessian}
0 < \nu I_M \leq \nabla_w^2\,J(w) \leq \delta I_M
\ee
For example, it is clear that the Hessian matrices in \eqref{eqn:JgradientHessiandef} and \eqref{eqn:loglossHessiandef} satisfy this property since
\bs
\be
\label{eqn:boundedHessianmsedef}
2 \lambdamin(R_u) I_M \leq \nabla_w^2\,J(w) \leq 2 \lambdamax(R_u) I_M
\ee
in the first case and
\be
\rho I_M \leq \nabla_w^2\,J(w) \leq (\rho + \lambdamax(R_h) ) I_M
\ee
\es
in the second case, where the notation $\lambdamin(R)$ and $\lambdamax(R)$ refers to the smallest and largest eigenvalues of the symmetric matrix argument, $R$, respectively. In summary, we will be assuming the following conditions \cite{NOW2014,Chen10b,chen2013,xiao2013}.

{%\small
\begin{assumption}[{\bf Conditions on cost function}] {\em The cost function $J(w)$ is twice-differentiable and satisfies \eqref{eqn:boundedHessian} for some positive parameters $\nu \leq \delta$. Condition \eqref{eqn:boundedHessian} is equivalent to requiring $J(w)$ to be $\nu$-strongly convex and for its gradient vector to be $\delta$-Lipschitz as in \eqref{eqn:nustronglyconvexdef} and \eqref{eqn:deltaLipschitzdef}, respectively.}

\hfill \qd

\label{assumption:costfunctions}
\end{assumption}
}

\subsection{Stochastic-Gradient Approximation}
 The traditional gradient-descent algorithm for solving \eqref{eqn:minimizeproblem} takes the form:
\be
\label{eqn:GD}
w_{i} \weq w_{i-1} - \mu \,\nabla_{w\tran} J(w_{i-1}), \quad i \geq 0
\ee
where $i \geq 0$ is an iteration index and $\mu > 0$ is a small step-size parameter. Starting from some initial condition, $w_{-1}$, the iterates $\{w_i\}$ correspond to successive estimates for the minimizer $w^o$. In order to run recursion \eqref{eqn:GD}, we need to have access to the true gradient vector. This information is generally unavailable in most instances involving learning from data. For example, when cost functions are defined as the expectations of certain loss functions as in \eqref{eqn:losscostdef}, the statistical distribution of the data $\x$ may not be known beforehand. In that case, the exact form of $J(w)$ will not be known since the expectation of $Q(w;\x)$ cannot be computed. In such situations, it is customary to replace the true gradient vector, $\nabla_{w\tran}\,J(w_{i-1})$, by an instantaneous approximation for it, and which we shall denote by $\wh{\nabla_{w\tran}\,J}(\w_{i-1})$. Doing so leads to the following {\em stochastic-gradient} recursion in lieu of \eqref{eqn:GD}:
\be
\label{eqn:SGD}
\addbox{\;\w_{i} \weq \w_{i-1} - \mu\, \wh{\nabla_{w\tran}\,J}(\w_{i-1}), \quad i\geq 0\;}
\ee
We use the {\bf boldface} notation, $\w_{i}$, for the iterates in \eqref{eqn:SGD} to highlight the fact that these iterates are now randomly perturbed versions of the values $\{w_i\}$ generated by the original recursion \eqref{eqn:GD}. The random perturbations arise from the use of the approximate gradient vector. The boldface notation is therefore meant to emphasize the random nature of the iterates in \eqref{eqn:SGD}. We refer to recursion \eqref{eqn:SGD} as a {\em synchronous} implementation since updates occur continuously over the iteration index $i$. This terminology is meant to distinguish the above recursion from its {\em asynchronous} counterpart, which is introduced later in Sec.~\ref{sec.rand.alkl1k3}.

We illustrate construction \eqref{eqn:SGD} by considering a scenario from classical adaptive filter theory \cite{Sayed08,Hay02,Widrow85}, where the gradient vector is approximated directly from data realizations. The construction will reveal why stochastic-gradient implementations of the form \eqref{eqn:SGD}, using approximate rather than exact gradient information, become naturally endowed with the ability to respond to {\em streaming} data.

{\small
\begin{example}[{\bf LMS adaptation}]
\label{example:LMS}{\rm
Let $\d(i)$ denote a streaming sequence of zero-mean random variables with variance $\sigma_{d}^2 = \Ex\d^2(i)$. Let $\u_{i}$ denote a streaming sequence of $1\times M$ independent zero-mean random vectors with covariance matrix $R_u =\Ex \u_{i}\tran\u_{i} > 0$. Both processes $\{\d(i),\u_{i}\}$ are assumed to be jointly wide-sense stationary. The cross-covariance vector between $\d(i)$ and $\u_{i}$ is denoted by $r_{du} = \Ex \d(i)\u_{i}\tran$. The data $\{\d(i),\u_i\}$ are assumed to be related via a linear regression model of the form:
\bs
\be
\label{eqn:linearmodel}
\d(i) \weq \u_i w^o + \v(i)
\ee
for some unknown parameter vector $w^o$, and where $\v(i)$ is a zero-mean white-noise process with power $\sigma_v^2 = \Ex\v^2(i)$ and assumed independent of $\u_j$ for all $i,j$. Observe that we are using parentheses to represent the time-dependency of a scalar variable, such as writing $\d(i)$, and subscripts to represent the time-dependency of a vector variable, such as writing $\u_i$. This convention will be used throughout the chapter. In a manner similar to Example~\ref{example-1}, we again pose the problem of estimating $w^o$ by minimizing the mean-square error cost
\be
J(w) \weq \Ex \left( \d(i) - \u_{i} w \right)^2 \equiv \Ex Q(w;\x_i)
\ee
where now the quantities $\{\d(i),\u_i\}$ represent the random data $\x_i$ in the definition of the loss function, $Q(w;\x_i)$. Using \eqref{eqn:GD}, the gradient-descent recursion in this case will take the form:
\be
\label{eqn:GDmse}
w_{i} \weq w_{i-1} - 2 \mu \left[ R_{u} w_{i-1} - r_{du} \right], \quad i\geq 0
\ee
The main difficulty in running this recursion is that it requires knowledge of the moments $\{r_{du},R_{u}\}$. This information is rarely available beforehand; the adaptive agent senses instead realizations $\{\d(i),\u_{i}\}$ whose statistical distributions have moments $\{r_{du},R_{u}\}$. The agent can use these realizations to approximate the moments and the true gradient vector. There are many constructions that can be used for this purpose, with different constructions leading to different adaptive algorithms \cite{Sayed03,Sayed08,Hay02,Widrow85}. It is sufficient to focus on one of the most popular adaptive algorithms, which results from using the data $\{\d(i),\u_{i}\}$ to compute {\em instantaneous} approximations for the  unavailable moments as follows:
\be
r_{du} \approx \d(i)\u_{i}\tran, \qquad
R_{u} \approx \u_{i}\tran\u_{i}
\ee
By doing so, the true gradient vector is approximated by:
\be
\label{eqn:lmsgradientapprox}
\wh{\nabla_{w\tran}\,J}(w) \weq 2 \left[\u_{i}\tran \u_{i} w - \u_{i}\tran \d(i) \right]
= \nabla_{w\tran}\,Q(w;\x_i)
\ee
Observe that this construction amounts to replacing the true gradient vector, $\nabla_{w\tran}\,J(w)$, by the gradient vector of the loss function itself (which, equivalently, amounts to dropping the expectation operator). Substituting \eqref{eqn:lmsgradientapprox} into \eqref{eqn:GDmse} leads to the well-known (synchronous) least-mean-squares (LMS, for short) algorithm \cite{Hay02,Widrow85,Sayed08,widrowhoff1960}:
\be
\label{eqn:LMS}
\addbox{\;\w_{i} \weq \w_{i-1} + 2 \mu \u_{i}\tran \left[ \d(i) - \u_{i}\w_{i-1} \right], \quad i \geq 0\;}
\ee
\es
The LMS algorithm is therefore a stochastic-gradient algorithm. By relying directly on the instantaneous data $\{\d(i),\u_i\}$, the algorithm is infused with useful tracking abilities. This is because drifts in the model $w^o$ from \eqref{eqn:linearmodel} will be reflected in the data $\{\d(i),\u_i\}$, which are used directly in the update \eqref{eqn:LMS}.

\qd}
\label{example-3}
\end{example}
}

If desired, it is also possible to employ iteration-dependent step-size sequences, $\mu(i)$, in \eqref{eqn:SGD} instead of the constant step-size $\mu$, and to require $\mu(i)$ to satisfy
\be
\label{eqn:diminishingstepsize}
\sum_{i=0}^{\infty} \mu^2(i) < \infty, \qquad \sum_{i=0}^{\infty} \mu(i) = \infty
\ee
Under some technical conditions, it is well-known that such step-size sequences ensure the convergence of $\w_i$ towards $w^o$ almost surely as $i\rightarrow \infty$ \cite{NOW2014,Poljak87,Ber97,tsypkin1971}. However, conditions \eqref{eqn:diminishingstepsize} force the step-size sequence to decay to zero, which is problematic for applications requiring continuous adaptation and learning from streaming data. This is because, in such applications, it is not unusual for the location of the minimizer, $w^o$, to drift with time. With $\mu(i)$ decaying towards zero, the stochastic-gradient algorithm \eqref{eqn:SGD} will stop updating and will not be able to track drifts in the solution. For this reason, we shall focus on constant step-sizes from this point onwards since we are interested in solutions with tracking abilities.

Now, the use of an approximate  gradient vector in  \eqref{eqn:SGD} introduces perturbations relative to the operation of the original recursion \eqref{eqn:GD}. We refer to the perturbation as gradient noise and define it as the difference:
\be
\label{eqn:gradientnoisedef}
\addbox{\;\s_{i}(\w_{i-1}) \define \wh{\nabla_{w\tran}\,J}(\w_{i-1}) - \nabla_{w\tran}\,J(\w_{i-1})\;}
\ee
The presence of this perturbation prevents the stochastic iterate, $\w_{i}$, from converging almost-surely to the minimizer $w^o$ when constant step-sizes are used. Some deterioration in performance will occur and the iterate $\w_i$ will instead fluctuate close to $w^o$. We will assess the size of these fluctuations by measuring their steady-state mean-square value (also called mean-square-deviation or MSD). It will turn out that the MSD is small and in the order of $O(\mu)$ --- see \eqref{eqn:SGDMSDOmu} further ahead. It will also turn out that stochastic-gradient algorithms converge towards their MSD levels at a geometric rate. In this way, we will be able to conclude that adaptation with small constant step-sizes can still lead to reliable performance in the presence of gradient noise, which is a reassuring result. We will also be able to conclude that adaptation with constant step-sizes is useful even for stationary environments. This is because it is generally sufficient in practice to reach an iterate $\w_{i}$ within some fidelity level from $w^o$ in a \emph{finite} number of iterations. As long as the MSD level is satisfactory, a stochastic-gradient algorithm will be able to attain satisfactory fidelity within a reasonable time frame. In comparison, although diminishing step-sizes ensure almost-sure convergence of $\w_i$ to $w^o$, they nevertheless disable tracking and can only guarantee slower than geometric rates of convergence (see, e.g., \cite{NOW2014,Poljak87,Ber97,NOW2014}). The next example from \cite{Chen10b} illustrates the nature of the gradient noise process \eqref{eqn:gradientnoisedef} in the context of mean-square-error adaptation.

{\small
\begin{example}[{\bf Gradient noise}]
\label{example:Gradientnoise}{\rm
It is clear from the expressions in Example~\ref{example-3} that the corresponding gradient noise process is
\bs
\be
\label{eqn:LMSgradientnoisedef}
\s_{i}(\w_{i-1}) \weq 2 \left( R_{u} - \u_{i}\tran \u_{i} \right) \,\wt{\w}_{i-1} - 2 \u_{i}\tran \v(i)
\ee
where we introduced the error vector:
\be
\wt{\w}_{i} \define w^o - \w_{i}
\ee
Let the symbol $\Filtration_{i-1}$  represent the collection of all possible random events generated by the past iterates $\{\w_{j}\}$ up to time $i-1$ (more formally, $\Filtration_{i-1}$ is the \emph{filtration} generated by the random process $\w_j$ for $j \leq i-1$):
\be
\Filtration_{i-1} \define \mbox{\rm filtration} \left\{\w_{-1},\,\w_{0},\,\w_{1},\dots,\w_{i-1}\right\}
\ee
It follows from the conditions on the random processes $\{\u_{i},\v(i)\}$ in Example~~\ref{example-3} that
\begin{align}
\label{eqn:martingalegradientnoise}
\Ex \left[ \s_{i}(\w_{i-1}) | \Filtration_{i-1} \right] & = 0 \\
\label{eqn:boundedvargradientnoises}
\Ex \left[ \| \s_{i}(\w_{i-1}) \|^2 | \Filtration_{i-1} \right] & \leq 4c\, \| \wt{\w}_{i-1}\|^2 + 4 \sigma_{v}^2 \, \Tr(R_{u})
\end{align}
\es
for some constant $c \geq 0$. If we take expectations of both sides of \eqref{eqn:boundedvargradientnoises}, we further conclude that the variance of the gradient noise, $\Ex \| \s_{i}(\w_{i-1})\|^2$, is bounded by the combination of two factors. The first factor depends on the quality of the iterate, $\Ex \| \wt{\w}_{i-1}\|^2$, while the second factor depends on $\sigma_{v}^2$. Therefore, even if the adaptive agent is able to approach $w^o$ with great fidelity so that $\Ex\|\wt{\w}_{i-1}\|^2$ is small, the size of the gradient noise will still depend on $\sigma_{v}^2$.

\qd}
\label{example-4}
\end{example}
}

\subsection{Conditions on Gradient Noise Process}
 In order to examine the convergence and performance properties of the stochastic-gradient recursion \eqref{eqn:SGD}, it is necessary to introduce some assumptions on the stochastic nature of the gradient noise process, $\s_i(\cdot)$.  The conditions that we introduce in the sequel are similar to conditions used earlier in the optimization literature, e.g., in \cite[pp.~95--102]{Poljak87} and \cite[p.~635]{berst00}; they are also motivated by the conditions we observed in the mean-square-error case in Example~\ref{example-4}. Following the developments in \cite{Chen10b,xiao2013,chen2013}, we let
\bs
\be
\label{eqn:Rsidef}
R_{s,i}(\w_{i-1}) \define \Ex \left[\s_i(\w_{i-1})\s_i\tran(\w_{i-1}) | \Filtration_{i-1} \right]
\ee
denote the conditional second-order moment of the gradient noise process, which generally depends on $i$. We assume that, in the limit, the covariance matrix tends to a constant value when evaluated at $w^o$ and is denoted by
\be
\label{eqn:Rsdef}
\addbox{\;R_{s} \define \lim_{i\rightarrow\infty} \Ex \left[\s_{i}(w^o) \s_{i}\tran(w^o) | \Filtration_{i-1} \right]\;}
\ee
\es
For example, comparing with expression \eqref{eqn:LMSgradientnoisedef} for mean-square-error costs, we have
\bs
\begin{align}
\label{eqn:absolutegradientnoisedef}
\s_i(w^o) & \weq - 2 \u_i\tran \v(i) \\
\label{eqn:Rsabsolutegradientnoise}
R_{s} & \weq 4 \sigma_{v}^2 R_{u}
\end{align}
\es

{%\small
\begin{assumption}[{\bf Conditions on gradient noise}]
\label{assumption:gradientnoise}  {\em
It is assumed that the first and second-order conditional moments of the gradient noise process satisfy \eqref{eqn:Rsdef} and
\bs
\begin{align}
\label{eqn:martingalegradientnoisecond}
\Ex\left[ \s_{i}(\w_{i-1}) | \Filtration_{i-1} \right] & = 0 \\
\label{eqn:boundedvargradientnoisecond}
\Ex\left[ \|{\s}_{i}(\w_{i-1}) \|^2 | \Filtration_{i-1} \right] & \leq \beta^2 \, \|\wt{\w}_{i-1}\|^2 + \sigma_{s}^2
\end{align}
\es
almost surely, for some nonnegative scalars $\beta^2$ and $\sigma_s^2$.}

\hfill
\qd
\end{assumption}
}

\noindent Condition \eqref{eqn:martingalegradientnoisecond} ensures that the approximate gradient vector is unbiased. It follows from conditions \eqref{eqn:martingalegradientnoisecond}--\eqref{eqn:boundedvargradientnoisecond} that the gradient noise process itself satisfies:
\bs
\begin{align}
\label{eqn:unbiasedgradient}
\Ex \s_{i}(\w_{i-1}) & = 0 \\
\label{eqn:boundedvargradient}
\Ex \| \s_{i}(\w_{i-1}) \|^2 & \leq \beta^2 \, \Ex \| \wt{\w}_{i-1} \|^2 + \sigma_{s}^2
\end{align}
\es
It is straightforward to verify that the gradient noise process \eqref{eqn:LMSgradientnoisedef} in the mean-square-error case satisfies conditions \eqref{eqn:martingalegradientnoisecond}--\eqref{eqn:boundedvargradientnoisecond}. Note in particular from \eqref{eqn:boundedvargradientnoises} that we can make the identifications $\sigma_{s}^2 \rightarrow 4 \sigma_{v}^2 \Tr(R_{u})$ and $\beta^2 \rightarrow 4c$.

\subsection{Random Updates}\label{sec.rand.alkl1k3}
We examined the performance of synchronous updates of the form \eqref{eqn:SGD} in some detail in \cite{ProcIEEE2014,NOW2014}. As indicated earlier, the focus of the current chapter is on extending the treatment from \cite{ProcIEEE2014} to asynchronous implementations. Accordingly, the first main digression in the exposition relative to \cite{ProcIEEE2014} occurs at this stage.

Thus, note that the stochastic-gradient recursion \eqref{eqn:SGD} employs a constant step-size parameter, $\mu > 0$. This means that this implementation expects the approximate gradient vector, $\wh{\nabla_{w\tran}\;J}(\w_{i-1})$, to be available at every iteration. Nevertheless, there are situations where data may arrive at the agent at random times, in which case the updates will also be occurring at random times. One way to capture this behavior is to model the step-size parameter as a {\em random} process, which we shall denote by the boldface notation $\bm{\mu}(i)$ (since boldface letters in our notation refer to random quantities). Doing so, allows us to replace the synchronous implementation \eqref{eqn:SGD} by the following asynchronous recursion:
\be
\label{eqn:randomstepsizeSGD}
\addbox{\;\w_i \weq \w_{i-1} - \bmu(i) \wh{\nabla_{w\tran}\,J}(\w_{i-1}), \quad i \ge 0\;}
\ee
Observe that we are attaching a time index to the step-size parameter in order to highlight that its value will now be changing randomly from one iteration to another.

For example, one particular instance discussed further ahead in Example~\ref{example-5} is when $\bm{\mu}(i)$ is a Bernoulli random variable assuming one of two possible values, say, $\bm{\mu}(i)=\mu>0$ with probability $p_{\mu}$ and $\bm{\mu}(i)=0$ with probability $1-p_{\mu}$. In this case, recursion (\ref{eqn:randomstepsizeSGD}) will be updating $p_{\mu}-$fraction of the time. We will not limit our presentation to the Bernoulli model but will allow the random step-size to assume broader probability distributions; we denote its first and second-order moments as follows.

{%\small
\begin{assumption}[{\bf Conditions on step-size process}]
\label{assumption:randomstepsize} {\em It is assumed that the stochastic process $\{ \bmu(i), i \ge 0 \}$ consists of a sequence of independent and bounded random variables, $\bmu(i) \in [0, \muub]$, where $\muub > 0$ is a constant upper bound. The mean and variance of $\bmu(i)$ are fixed over $i$, and they are denoted by:
\bs
\begin{align}
\label{eqn:meanrandomstepsize}
\bar{\mu} & \define \Ex \bmu(i) \\
\label{eqn:varrandomstepsize}
\sigma_{\mu}^2 & \define \Ex (\bmu(i) - \bar{\mu})^2
\end{align}
\es
with $\bar{\mu} > 0$ and $\sigma_{\mu}^2 \geq 0$.  Moreover, it is assumed that, for any $i$, the random variable $\bmu(i)$ is independent of any other random variable in the learning algorithm.}

\hfill \qd
\label{assum.kjla;d.;lk;1l3}
\end{assumption}
}

\noindent The following variable will play an important role in characterizing the mean-square-error performance of asynchronous updates and, hence, we introduce a unique symbol for it:
\be
\addbox{\;\mu_{x}\define \bar{\mu}\;+\;\frac{\sigma_{\mu}^2}{\bar{\mu}}\;}\label{jahd7813la.dlk}
\ee
where the subscript ``$x$'' is meant to refer to the ``asynchronous'' mode of operation. This expression captures the first and second-order moments of the random variations in the step-size parameter into a single variable, $\mu_{x}$. While the constant step-size $\mu$ determines the performance of the synchronous implementation \eqref{eqn:SGD}, it turns out that the constant variable $\mu_x$ defined above will play an equivalent role for the asynchronous implementation \eqref{eqn:randomstepsizeSGD}. Note further that the synchronous stochastic-gradient iteration \eqref{eqn:SGD} can be viewed as a special case of recursion \eqref{eqn:randomstepsizeSGD} when the variance of $\bmu(i)$ is set to zero, i.e., $\sigma_{\mu}^2=0$, and the mean value $\bar{\mu}$ is set to $\mu$. Therefore, by using these substitutions, we will able to deduce performance metrics for \eqref{eqn:SGD} from the performance metrics that we shall present for \eqref{eqn:randomstepsizeSGD}. The following two examples illustrate situations involving random updates.

{\small
\begin{example}[{\bf Random updates under a Bernoulli model}]
\label{example:Bernoulli}{\rm Assume that at every iteration $i$, the agent adopts a random ``on-off'' policy to reduce energy consumption. It tosses a coin to decide whether to enter an active learning mode or a sleeping mode. Let $0 < p_{\mu} < 1$ denote the probability of entering the active mode. During the active mode, the agent employs a step-size value $\mu$. This model is useful for  situations in which data arrives randomly at the agent: at every iteration $i$, new data is available with probability $p_{\mu}$. The random step-size process is therefore of the following type:
\bs
\be
\bmu(i) \weq \begin{cases}
\mu, & \quad {\small \mbox{{\small with probability}}} \;\; p_{\mu} \\
0, & \quad {\small \mbox{{\small with probability}}} \;\; 1 - p_{\mu} \\
\end{cases}
\ee
In this case, the mean and variance of $\bmu(i)$ are given by:
\be
\bar{\mu} \weq p_{\mu} \mu, \qquad \sigma_{\mu}^2 \weq p_{\mu}(1-p_{\mu}) \mu^2,\qquad
\mu_{x}\;=\;\mu
\label{kajdhg,1lkl1kk3}\ee
\es
\qd}
\label{example-5}
\end{example}
}

{\small
\begin{example}[{\bf Random updates under a Beta model}]
\label{example:Beta}{\rm  Since the random step-size $\bm{\mu}(i)$ is limited to a finite-length interval $[0,\muub]$, we may extend the Bernoulli model from the previous example by adopting a  more general continuous Beta distribution for $\bm{\mu}(i)$. The Beta distribution is an extension of the Bernoulli distribution. While the Bernoulli distribution assumes two discrete possibilities for the random variable, say, $\{0,\mu\}$, the Beta distribution allows for any value in the continuum $[0,\mu]$.

Thus, let $\x$ denote a generic scalar random variable that assumes values in the interval $[0,1]$ according to a Beta distribution. Then, according to this distribution, the pdf of $\x$, denoted by $f_{\x}(x;\xi,\zeta)$, is determined by two shape parameters $\{\xi,\zeta\}$ as follows  \cite{Feller1971,Hahn1994}:
\bs
\be
\label{eqn:standard_beta}
f_{\x}(x; \xi, \zeta) \weq
\begin{cases}
\displaystyle \frac{\Gamma(\xi + \zeta)}{\Gamma(\xi)\Gamma(\zeta)}\,x^{\xi - 1}(1 - x)^{\zeta - 1},  & 0 \le x \le 1 \\
0, & \mbox{otherwise} \\
\end{cases}
\ee
where $\Gamma(\cdot)$ denotes the Gamma function \cite{Abra1964,And1999}. Figure~\ref{fig:beta} plots $f_{\x}(x; \xi, \zeta)$ for two values of $\zeta$. The mean and variance of the Beta distribution \eqref{eqn:standard_beta} are given by:
\be
\label{eqn:beta_meanandvariance}
\bar{x} = \frac{\xi}{\xi + \zeta}, \qquad \sigma_x^2 = \frac{\xi \zeta}{(\xi + \zeta)^2(\xi + \zeta + 1)}
\ee
\es
We note that the classical uniform distribution over the interval $[0,1]$ is a special case of the Beta distribution for $\xi=\zeta=1$ --- see Fig.~\ref{fig:beta}. Likewise, the Bernoulli distribution with $p_{\mu}=1/2$ is recovered from the Beta distribution by letting $\xi=\zeta\rightarrow 0$.

\begin{figure}[h]
\epsfxsize 13cm \epsfclipon
\begin{center}
\leavevmode \epsffile{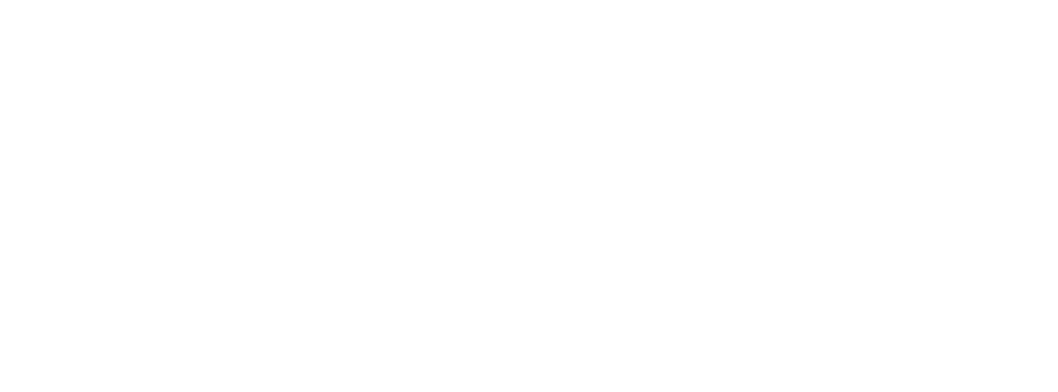}
\caption{{\small The pdf of the Beta distribution, $f_{\x}(x; \xi, \zeta)$, defined by (\ref{eqn:standard_beta}) for different values of the shape parameters $\xi$ and $\zeta$. Figure extracted with permission from \cite{xiao2013}.}}\label{fig:beta}
\end{center}
\end{figure}

In the Beta model for asynchronous adaptation, we assume that the ratio $\bmu(i) / \muub$ follows  a Beta distribution with parameters $\{\xi,\zeta\}$. Under this model, the mean and variance of the random step-size become:
\be
\bar{\mu} \weq \left(\frac{\xi}{\xi + \zeta}\right)\, \muub, \qquad
\sigma_{\mu}^2 \weq \left(\frac{\xi \zeta}{(\xi + \zeta)^2 (\xi + \zeta + 1)}\right)\, \muub^2
\ee

\qd}
\end{example}
}

\subsection{Mean-Square-Error Stability}
We now examine the convergence of the asynchronous stochastic-gradient recursion \eqref{eqn:randomstepsizeSGD}. In the statement below, the notation $a = O(\mu)$ means $a \leq b\mu$ for some constant $b$ that is independent of $\mu$.\\

{%\small
\begin{lemma}[{\bf Mean-square-error stability}]
\label{lemma:msestability}
Assume the conditions under Assumptions \ref{assumption:costfunctions}, \ref{assumption:gradientnoise}, and \ref{assumption:randomstepsize} on the cost function, the gradient noise process, and the random step-size process hold. Let $\mu_o=2\nu/(\delta^2+\beta^2)$. For any $\mu_{x}$ satisfying
\be
\label{eqn:stabilitycondstandalone}
\mu_{x} < \mu_o
\ee
it holds that $\Ex\|\wt{\w}_i\|^2$ converges exponentially (i.e., at a geometric rate) according to the recursion
\bs
\be
\label{eqn:msestabilitybound}
\Ex\|\wt{\w}_i\|^2 \;\leq\; \alpha\, \Ex\|\wt{\w}_{i-1}\|^2 + (\bar{\mu}^2 + \sigma_{\mu}^2)\sigma_s^2
\ee
where the scalar $\alpha$ satisfies $0\leq \alpha<1$ and is given by
\bq
\alpha &\define &1 - 2 \nu \bar{\mu} + (\delta^2+\beta^2)(\bar{\mu}^2 + \sigma_{\mu}^2)\nn\\
&=&1-2\nu\bar{\mu}\;+\;O(\mu_{x}^2)\label{eqn:alphadef}
\eq
It follows from \eqref{eqn:msestabilitybound}  that, for sufficiently small step-sizes:
\be
\label{eqn:SGDMSDOmu}
\limsup_{i\rightarrow\infty} \; \Ex \|\wt{\w}_i\|^2 \weq O(\mu_{x})
\ee
\es
\end{lemma}
}

\noindent \bp We subtract $w^o$ from both sides of \eqref{eqn:randomstepsizeSGD} to get
\bs
\be
\label{eqn:SGDerrorrecursion}
\wt{\w}_i \weq \wt{\w}_{i-1} + \bmu(i)\,\nabla_{w\tran}\,J(\w_{i-1}) + \bmu(i)\,\s_i(\w_{i-1})
\ee
We now appeal to the mean-value theorem \cite{NOW2014,Poljak87,Rudin76} to write:
\begin{align}
\label{eqn:bHdef}
\nabla_{w\tran}\,J(\w_{i-1}) & = -\left[ \int_{0}^{1}\nabla_w^2\,J(w_o - t \wt{\w}_{i-1}) dt \right] \, \wt{\w}_{i-1}\;\define \;- \H_{i-1} \wt{\w}_{i-1}
\end{align}
where we are introducing the {\em symmetric} and {\em random} time-variant matrix $\H_{i-1}$ to represent the integral expression. Substituting into \eqref{eqn:SGDerrorrecursion}, we get
\be
\label{eqn:SGDerrorrecursion2}
\wt{\w}_i \weq [ I_M - \bmu(i) \H_{i-1} ] \, \wt{\w}_{i-1} + \bmu(i)\,\s_i(\w_{i-1})
\ee
\es
so that from Assumption \ref{assumption:gradientnoise}:
\bq \Ex\left[ \| \wt{\w}_i \|^2 | \Filtration_{i-1} \right] &\hspace{-0.1cm}\leq\hspace{-0.1cm}& \left(\Ex\left[\| I_M - \bmu(i)\, \H_{i-1} \|^2 | \Filtration_{i-1}\right]\right) \,\|\wt{\w}_{i-1}\|^2 + \left(\Ex\bmu^2(i)\right)\,\left(\Ex \left[ \|\s_i(\w_{i-1}) \|^2 | \Filtration_{i-1} \right]\right)
\nn\\
\label{eqn:boundvarLMSstep1}
\eq
It follows from \eqref{eqn:boundedHessian} that
\bs\begin{align}
\label{eqn:bound2normstep1}
\|I_M - \bmu(i) \, \H_{i-1} \|^2 & = \left[ \rho\left( I_M - \bmu(i) \, \H_{i-1} \right) \right]^2 \nn\\
& \leq \max\left\{ \left[ 1 - \bmu(i)\, \delta \right]^2,\; \left[ 1 - \bmu(i) \, \nu \right]^2 \right\} \nn\\
& \leq 1 - 2\bmu(i)\, \nu + \bmu^2(i) \delta^2
\end{align}
since $\nu \leq \delta$. In the first line above, the notation $\rho(A)$ denotes the spectral radius of its matrix argument (i.e., $\rho(A) = \max_k |\lambda_k(A)|$ in terms of the largest magnitude eigenvalue of $A$). From \eqref{eqn:bound2normstep1} we obtain
\be
\label{eqn:boundonalpha}
\Ex\left(\| I_M - \bmu(i)\, \H_{i-1} \|^2 | \Filtration_{i-1}\right) \leq 1 - 2 \bar{\mu} \nu + (\bar{\mu}^2 + \sigma_{\mu}^2) \delta^2
\ee
\es
Taking expectations of both sides of \eqref{eqn:boundvarLMSstep1}, we arrive at \eqref{eqn:msestabilitybound}
from \eqref{eqn:boundedvargradient} and \eqref{eqn:boundonalpha} with $\alpha$ given by \eqref{eqn:alphadef}. The bound in \eqref{eqn:stabilitycondstandalone} on the moments of the random step-size ensures that $0 \leq \alpha < 1$.  For the $O(\mu_x^2)$ approximation in expression \eqref{eqn:alphadef} note from \eqref{jahd7813la.dlk} that
\be
(\bar{\mu}^2+\sigma_{\mu}^2)\;=\;\bar{\mu}\mu_x\;\leq\;\mu_x^2
\ee
Iterating recursion \eqref{eqn:msestabilitybound} gives
\bs\be
\Ex \, \|\wt{\w}_i\|^2 \leq \alpha^{i+1}\,\Ex\, \|\wt{\w}_{-1}\|^2 + \frac{(\bar{\mu}^2 + \sigma_{\mu}^2)\sigma_s^2}{1-\alpha}
\ee
Since $0\leq \alpha <1$, there exists an iteration value $I_o$ large enough such that
\be
\alpha^{i+1}\,\Ex\, \|\wt{\w}_{-1}\|^2\;\leq\;\frac{(\bar{\mu}^2 + \sigma_{\mu}^2)\sigma_s^2}{1-\alpha},\;\;\;i>I_o
\ee
It follows that the variance $\Ex\,\|\wt{\w}_i\|^2$ converges exponentially to a region that is upper bounded by $2(\bar{\mu}^2 + \sigma_{\mu}^2)\sigma_s^2/(1-\alpha)$. It can be verified that this bound does not exceed $2\mu_{x} \sigma_s^2/\nu$, which is $O(\mu_{x})$, for any $\mu_{x} < \mu_o/2$.
\es

\ep

\subsection{Mean-Square-Error Performance}
\label{subsec:MSEperformanceStandalone}
 We conclude from \eqref{eqn:SGDMSDOmu} that  the mean-square error (MSE) can be made as small as desired by using small step-sizes, $\mu_x$. In this section we derive a closed-form expression for the asymptotic mean-square error, which is more frequently called the {\em mean-square deviation} (MSD) and is defined as:
\bs\be
\addbox{\;\MSD \define \lim_{i\rightarrow\infty}\; \Ex\, \|\wt{\w}_{i}\|^2\;}\label{kahdg657813.l;}
\ee
Strictly speaking, the limit on the right-hand side of the above expression may not exist. A more accurate definition for the MSD appears in Eq.~(4.86) of \cite{NOW2014}, namely,
\be
\MSD\define \mu_{x}\cdot\left(\lim_{\mu_{x}\rightarrow 0 }\limsup_{i\rightarrow\infty}\;{1\over \mu_x}\Ex\|\widetilde{\w}_i\|^2\right)
\label{makdjl13;l;13l}\ee
\es
However, it was explained in \cite{NOW2014}[Sec. 4.5] that derivations that assume the validity of \eqref{kahdg657813.l;} still lead to the same expression for the MSD to first-order in $\mu_{x}$ as derivations that rely on the more formal definition \eqref{makdjl13;l;13l}. We therefore continue with (\ref{kahdg657813.l;}) for simplicity of presentation.
 We explain below  how an expression for the MSD can be obtained by following the energy conservation technique of \cite{NOW2014,Sayed03,Sayed08,yousef2001,naffuori}. For that purpose, we need to introduce two smoothness conditions.

\bs
{%\small
\begin{assumption}[{\bf Smoothness conditions}]
\label{assumption:smoothness} {\em
In addition to Assumptions \ref{assumption:costfunctions} and \ref{assumption:gradientnoise}, we assume that the Hessian matrix of the cost function and the noise covariance matrix defined by \eqref{eqn:Rsidef} are locally Lipschitz continuous in a small neighborhood around $w = w^o$:
\begin{align}
\label{eqn:LipschitzHessian}
\left\| \nabla_w^2\,J( w^o + \delta w) - \nabla_w^2\,J(w^o) \right\| & \leq \tau\, \|\delta w\| \\
\label{eqn:LipschitzCovariance}
\left\| R_{s,i}(w^o + \delta w) - R_{s,i}(w^o) \right\| & \leq \tau_2\, \|\delta w\|^{\kappa}
\end{align}
for small perturbations $\|\delta w\|\leq r$  and for some $\tau,\tau_2 \geq 0$ and $1 \leq \kappa \leq 2$.}

\hfill \qd
\end{assumption}
}

\smallskip

\noindent The range of values for $\kappa$ can be enlarged, e.g., to $\kappa\in(0,4]$. The only change in allowing a wider range for $\kappa$ is that the exponent of the higher-order term, $O(\mu^{{\small 3/2}}_x)$, that will appear in several performance expressions, as is the case with \eqref{eqn:SGDMSDexpression}--\eqref{eqn:SGDERexpression}, will need to be adjusted from ${3\over 2}$ to  $\min\{{3\over 2}, 1+{\kappa\over 2}\}$, without affecting the first-order term that determines the MSD \cite{xiao2013AA,NOW2014,chen2013AA}. Therefore, it is sufficient to continue with $\kappa \in [1,2]$ to illustrate the key concepts though the MSD expressions will still be valid to first-order in $\mu_x$.

\begin{figure}[t]
\epsfxsize 11cm \epsfclipon
\begin{center}
\leavevmode \epsffile{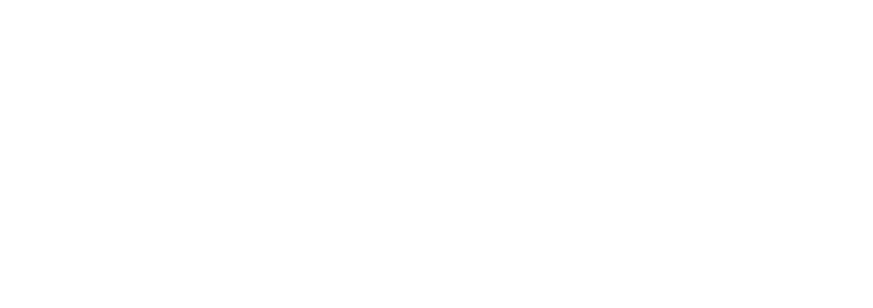}
\caption{{\small A block-diagram representation of the long-term recursion
\eqref{eqn:SGDerrorrecursionlongterm} for single-agent adaptation and learning. Figure extracted with permission  from \cite{ProcIEEE2014}.}}
\label{fig:blockdiagramlongterm}
\end{center}
\end{figure}

Using \eqref{eqn:boundedHessian}, it can be verified that condition \eqref{eqn:LipschitzHessian} translates into a global Lipschitz property relative to the minimizer $w^o$, i.e., it will also hold that \cite{NOW2014,xiao2013AA}:
\be
\label{eqn:globalLipschitzHessian}
\| \nabla_w^2\,J(w) - \nabla_w^2\,J(w^o) \| \leq \tau'\,\| w - w^o\|
\ee
\es
for all $w$ and for some $\tau'\geq 0$. For example, both conditions \eqref{eqn:LipschitzHessian}--\eqref{eqn:LipschitzCovariance} are readily satisfied by mean-square-error costs. Using property \eqref{eqn:globalLipschitzHessian}, we can now motivate a useful long-term model for the evolution of the error vector $\wt{\w}_i$ after sufficient iterations, i.e., for $i \gg 1$. Indeed, let us reconsider recursion \eqref{eqn:SGDerrorrecursion2} and introduce the deviation matrix:
\bs
\be
\label{eqn:Happroxerrordef}
\wt{\H}_{i-1} \define H - \H_{i-1}
\ee
where the constant (symmetric and positive-definite) matrix $H$ is defined as:
\be
\label{eqn:Hdef}
\addbox{\;H \define \nabla_w^2\,J(w^o)\;}
\ee
\es
Substituting \eqref{eqn:Happroxerrordef} into \eqref{eqn:SGDerrorrecursion2} gives
\bs\be
\label{eqn:SGDerrorrecursionapproxcterm}
\wt{\w}_i \weq ( I_M - \bmu(i) \, H ) \, \wt{\w}_{i-1} + \bmu(i)\,\s_i(\w_{i-1}) + \bmu(i) \, \c_{i-1}
\ee
where
\be
\c_{i-1} \define \wt{\H}_{i-1} \wt{\w}_{i-1}
\ee
\es
Using \eqref{eqn:globalLipschitzHessian} and the fact that $(\Ex\a)^2\leq \Ex\a^2$ for any real-valued random variable, we can bound the conditional expectation of the norm of the perturbation term as follows:
\bs
\begin{align}
\label{eqn:boundcprobability}
\Ex \left[\, \| \bmu(i) \, \c_{i-1} \|\, |  \Filtration_{i-1}\, \right] & =
\left(\Ex\bmu(i)\right)\,\left(\Ex \left[\, \|\c_{i-1} \|\, |  \Filtration_{i-1}\, \right]\right)\nn\\
&\leq
\sqrt{\left(\Ex\bmu^2(i)\right)}\,\Ex \left[ \|\c_{i-1} \| |  \Filtration_{i-1} \right]\nn\\
&\stackrel{(\ref{eqn:globalLipschitzHessian})}{\leq}\sqrt{\bar{\mu}^2 + \sigma_{\mu}^2} \cdot \frac{\tau'}{2} \, \|\wt{\w}_{i-1}\|^2 \nn \\
& \leq \frac{\bar{\mu}^2 + \sigma_{\mu}^2}{\bar{\mu}} \cdot \frac{\tau'}{2}\, \|\wt{\w}_{i-1}\|^2 \nn \\
& = \mu_x \cdot \frac{\tau'}{2}\, \|\wt{\w}_{i-1}\|^2
\end{align}
so that using \eqref{eqn:SGDMSDOmu}, we conclude that:
\be
\label{eqn:boundcprobability2}
\limsup_{i\rightarrow\infty} \; \Ex\,\| \bmu(i) \, \c_{i-1} \| \weq O(\mu_x^2)
\ee
We can deduce from this result that $\|\bmu(i) \, \c_{i-1}\| = O(\mu_x^2)$ asymptotically with {\em high probability} \cite{NOW2014,xiao2013AA}. To see this, let $r_c = m \mu_x^2$, for any constant integer $m\geq 1$.  Now, calling upon Markov's inequality \cite{Papoulis,Durret,Dudley03}, we conclude from \eqref{eqn:boundcprobability2} that for $i\gg 1$:
\begin{align}
\label{eqn:boundprobability3}
\Pr\left( \|\bmu(i) \, \c_{i-1} \| \,<\, r_c \right) & = 1 - \Pr\left( \|\bmu(i) \, \c_{i-1} \| \,\geq\, r_c \right)\nn\\
& \geq 1 - \frac{\Ex\|\bmu(i) \, \c_{i-1}\|}{r_c} \nn \\
& \stackrel{\rm \eqref{eqn:boundcprobability2}}{\geq} 1 - O\left( 1/m \right)
\end{align}
\es
This result shows that the probability of having $\|\bmu(i) \,\c_{i-1}\|$ bounded by $r_c$ can be made arbitrarily close to one by selecting a large enough value for $m$. Once the value for $m$ has been fixed to meet a desired confidence level, then $r_c = O(\mu_x^2)$. This analysis, along with recursion \eqref{eqn:SGDerrorrecursionapproxcterm},  motivate us to assess the mean-square performance of the error recursion \eqref{eqn:SGDerrorrecursion2} by considering instead the following long-term model, which holds with high probability after sufficient iterations $i\gg 1$:
\be
\label{eqn:SGDerrorrecursionlongterm}
\addbox{\;\wt{\w}_i \weq (I_M - \bmu(i)\, H) \wt{\w}_{i-1} + \bmu(i)\,\s_i(\w_{i-1}) + O(\mu_x^2)\;}
\ee
Working with iteration \eqref{eqn:SGDerrorrecursionlongterm} is helpful because its dynamics is driven by the constant matrix $H$ as opposed to the random matrix $\H_{i-1}$ in the original error recursion \eqref{eqn:SGDerrorrecursion2}. If desired, it can be shown that, under some technical conditions on the fourth-order moment of the gradient noise process, the MSD expression that will result from using \eqref{eqn:SGDerrorrecursionlongterm} is within $O(\mu_{x}^{3/2})$ of the actual MSD expression for the original recursion \eqref{eqn:SGDerrorrecursion2} --- see \cite{NOW2014,chen2013AA,xiao2013AA} for a formal proof of this fact. Therefore, it is sufficient to rely on the long-term model \eqref{eqn:SGDerrorrecursionlongterm} to obtain performance expressions that are accurate to first-order in $\mu_x$. Figure~\ref{fig:blockdiagramlongterm} provides a block-diagram representation for \eqref{eqn:SGDerrorrecursionlongterm}.

Before explaining how model \eqref{eqn:SGDerrorrecursionlongterm} can be used to assess the MSD, we remark that there is a second useful metric for evaluating the performance of stochastic gradient algorithms. This metric relates to the mean excess-cost; which is also called the {\em excess-risk} (ER) in the machine learning literature \cite{Bish2007,theo2008} and the excess-mean-square-error (EMSE) in the adaptive filtering literature \cite{Hay02,Widrow85,Sayed08}. We denote it by the letters ER and define it as the average fluctuation of the cost function around its minimum value:
\bs
\be
\label{eqn:ERdef}
\addbox{\;\ER \define \lim_{i\rightarrow\infty} \Ex \,\{ J(\w_{i-1}) - J(w^o) \}\;}
\ee
Using the smoothness condition \eqref{eqn:LipschitzHessian}, and the mean-value theorem \cite{Poljak87,Rudin76} again, it can be verified that \cite{xiao2013AA,NOW2014,chen2013AA}:
\be
\label{eqn:ERevaluation}
\ER \define \lim_{i\rightarrow\infty} \Ex \, \| \wt{\w}_{i-1} \|^2_{\frac{1}{2}H} + O(\mu_x^{3/2})
\ee
\es

\smallskip

{%\small
\begin{lemma}[{\bf Mean-square-error performance}]
\label{lemma:SGDMSEperformance}
Assume the conditions under Assumptions \ref{assumption:costfunctions}, \ref{assumption:gradientnoise}, \ref{assumption:randomstepsize}, and \ref{assumption:smoothness} on the cost function, the gradient noise process, and the random step-size process hold. Assume further that the asynchronous step-size parameter $\mu_x$ is sufficiently small to ensure mean-square stability as required by \eqref{eqn:stabilitycondstandalone}. Then, the {\rm MSD} and {\rm ER} metrics for the asynchronous stochastic-gradient algorithm \eqref{eqn:randomstepsizeSGD} are well-approximated to first-order in $\mu_x$ by the expressions:
\bs
\begin{align}
\label{eqn:SGDMSDexpression}
\mbox{\rm MSD}^{\rm asyn} & = \frac{\mu_x}{2}\,\Tr\left( H^{-1} R_s \right) + O(\mu_x^{3/2}) \\
\label{eqn:SGDERexpression}
\mbox{\rm ER}^{\rm asyn} & = \frac{\mu_x}{4}\,\Tr\left( R_s \right) + O(\mu_x^{3/2})
\end{align}
where $R_s$ and $H$ are defined by \eqref{eqn:Rsdef} and \eqref{eqn:Hdef}, and where we are adding the superscript ``asyn'' for clarity in order to distinguish these measures from the corresponding measures in the synchronous case (mentioned below in (\ref{synchda.,.13})--(\ref{synchda.,.13.b})). Moreover, we derived earlier in \eqref{eqn:alphadef} the following expression for the convergence rate:
\be
\label{eqn:alphadef.a}
\alpha^{\rm asyn}\;=\;1-2\nu\bar{\mu}\;+\;O(\mu_{x}^2)
\ee
\es

\end{lemma}
}

\bp
We introduce the eigen-decomposition $H = U \Lambda U\tran$, where $U$ is orthonormal and $\Lambda$ is diagonal with positive entries, and rewrite \eqref{eqn:SGDerrorrecursionlongterm} in terms of transformed quantities:
\bs
\be
\label{eqn:eigenerrorrecursion}
\ol{\w}_i \weq \left( I_M - \bmu(i) \,\Lambda \right) \ol{\w}_{i-1} + \bmu(i) \, \ol{\s}_i (\w_{i-1}) + O(\mu_x^2)
\ee
where $\ol{\w}_i=U\tran \wt{\w}_i$ and $\ol{\s}_i(\w_{i-1})=U\tran {\s}_i(\w_{i-1})$. Let $\Sigma$ denote an arbitrary $M\times M$ diagonal matrix with positive entries that we are free to choose. Then, equating the weighted squared norms of both sides of \eqref{eqn:eigenerrorrecursion} and taking expectations gives for $i\gg 1$:
\be
\label{eqn:energyeigenerrorrecursion}
\Ex \|\ol{\w}_{i}\|_{\Sigma}^2 \weq \Ex \| \ol{\w}_{i-1} \|_{\Sigma'}^2 + (\bar{\mu}^2 + \sigma_{\mu}^2) \,
\Ex \|\ol{\s}_i(\w_{i-1})\|_{\Sigma}^2 + O(\mu_x^{5/2})
\ee
where
\be
\label{eqn:Sigmaprimedef}
\Sigma' \define \Ex (I_M - \bmu(i) \, \Lambda) \Sigma (I_M - \bmu(i) \, \Lambda) = \Sigma - 2\bar{\mu} \Lambda \Sigma + O(\mu_x^2)
\ee
From \eqref{eqn:Rsdef}, \eqref{eqn:boundedvargradientnoisecond}, \eqref{eqn:SGDMSDOmu}, and \eqref{eqn:LipschitzCovariance} we obtain:
\be
\lim_{i\rightarrow\infty} \Ex \|\ol{\s}_i(\w_{i-1})\|^2_{\Sigma} \weq \Tr(U \Sigma U\tran R_s) + O(\mu_x^{\kappa/2})\ee
Therefore, substituting into \eqref{eqn:energyeigenerrorrecursion} gives for $i\rightarrow\infty$:
\be
\label{eqn:weightedMSDlimitexpression}
\lim_{i\rightarrow\infty} \Ex \|\ol{\w}_{i}\|^2_{2\Lambda\Sigma} \weq \mu_x \, \Tr( U \Sigma U\tran R_s) + O(\mu_x^{3/2})
\ee
Since we are free to choose $\Sigma$, we let $\Sigma = \frac{1}{2} \Lambda^{-1}$ and arrive at \eqref{eqn:SGDMSDexpression} since $\|\ol{\w}_i\|^2=\|\wt{\w}_i\|^2$ and $U \Sigma U\tran = \frac{1}{2} H^{-1}$. On the other hand, selecting $\Sigma = \frac{1}{4} I_M$ leads to \eqref{eqn:SGDERexpression}.
\es

\ep

\noindent We recall our earlier remark that the synchronous stochastic-gradient recursion \eqref{eqn:SGD} can be viewed as a special case of the asynchronous update \eqref{eqn:randomstepsizeSGD} by setting $\sigma_{\mu}^2=0$ and $\mu_x=\bar{\mu}\equiv \mu$. Substituting these values into \eqref{eqn:SGDMSDexpression}--\eqref{eqn:alphadef.a}, we obtain for the synchronous implementation \eqref{eqn:SGD}:
\bs
\begin{align}
\MSD^{\sync} & \weq \frac{\mu}{2}\,\Tr\left( H^{-1} R_s \right) + O(\mu^{3/2})\label{synchda.,.13} \\
\ER^{\sync} & \weq \frac{\mu}{2}\,\Tr\left(R_s \right) + O(\mu^{3/2})\label{synchda.,.13.c} \\
\alpha^{\sync} & \weq 1 - 2 \nu \mu + O(\mu^2)\label{synchda.,.13.b}
\end{align}
\es
which are the same expressions presented in \cite{ProcIEEE2014} and which agree with classical results for LMS adaptation \cite{Widrow76Proc}--\cite{Foley88TSP}. The matrices $R_s$ that appear in these expressions, and in
\eqref{eqn:SGDMSDexpression}--\eqref{eqn:alphadef.a}, were defined earlier in \eqref{eqn:Rsdef}  and they correspond to the covariance matrices of the gradient noise processes in their respective (synchronous or asynchronous) implementations. The examples that follow show how expressions \eqref{eqn:SGDMSDexpression} and \eqref{eqn:SGDERexpression} can be used to recover performance metrics for mean-square-error adaptation and learning under random updates.

{\small
\begin{example}[{\bf Performance of asynchronous LMS adaptation}]
\label{example:LMSperformance}{\rm
We reconsider the LMS recursion \eqref{eqn:LMS} albeit in an asynchronous mode of operation, namely,
\be
\label{eqn:LMS.async}
\addbox{\;\w_{i} \weq \w_{i-1} + 2 \bm{\mu}(i) \u_{i}\tran \left[ \d(i) - \u_{i}\w_{i-1} \right], \quad i \geq 0\;}
\ee
We know from Example \ref{example:LMS} and \eqref{eqn:absolutegradientnoisedef} that this situation corresponds to $H = 2R_u$ and $R_s = 4 \sigma_v^2 R_u$. Substituting into \eqref{eqn:SGDMSDexpression} and \eqref{eqn:SGDERexpression} leads to the following expressions for the MSD and EMSE of the asynchronous LMS filter:
\bs
\bq
\MSD^{\rm asyn}_{{\small \lms}} & \approx& \mu_x M \sigma_{v}^2 = O(\mu_x) \label{eqn:LMSMSDexpressioncloseform}\\
\EMSE^{\rm asyn}_{{\small \lms}} & \approx& \mu_x \sigma_{v}^2 \Tr(R_{u}) = O(\mu_x)\label{eqn:eqn:LMSEMSEexpressioncloseform}
\eq
\es
where here, and elsewhere, we will be using the notation $\approx$ to indicate that we are ignoring higher-order terms in $\mu_x$. For example, let us assume a Bernoulli update model for $\bm{\mu}(i)$ where the filter updates with probability $p_{\mu}$ using a step-size value $\mu$ or stays inactive otherwise. In this case, we conclude from
\eqref{kajdhg,1lkl1kk3} that $\mu_x=\mu$ so that the above performance expressions for the MSD and EMSE metrics will coincide with the values obtained in the synchronous case as well. In other words, the steady-state performance levels are not affected whether the algorithm learns in a synchronous or asynchronous manner. However, the convergence rate is affected since $\bar{\mu}=\mu p_{\mu}$ and, therefore,
\bs
\bq
\alpha^{\rm sync}_{\lms} &\approx& 1 - 2 \nu \mu \\
\alpha^{\rm asyn}_{\lms} &\approx& 1 - 2 \nu \mu p_{\mu} \;>\;\alpha^{\rm sync}_{\lms}\eq
\es
It follows that asynchronous LMS adaptation attains the same performance levels as synchronous LMS adaptation albeit at a slower convergence rate.

We may alternatively compare the performance of the synchronous and asynchronous implementations by fixing their convergence rates to the same value. Thus, consider now a second random update scheme with mean $\bar{\mu}$ and let us set $\mu=\bar{\mu}$. That is, the step-size used by the synchronous implementation is set equal to the mean step-size used by the asynchronous implementation. Then, in this case, we will get $\alpha^{\rm sync}_{\lms}=\alpha^{\rm asyn}_{\lms}$ so that the convergence rates coincide to first-order. However, it now holds that $\MSD^{\rm asyn}_{{\small \lms}}>\MSD^{\rm sync}_{{\small \lms}}$ since $\bar{\mu}_x>\mu$ so that some deterioration in MSD performance occurs.

\qd}
\end{example}
}

{\small
\begin{example}[{\bf Performance of asynchronous online learners}]
\label{example:performanceonlinelearner}{\rm
Consider a stand-alone learner receiving a streaming sequence of independent data vectors $\{\x_{i}, i\geq 0\}$ that arise from some fixed probability distribution $\Xcal$. The goal is to learn the vector $w^o$ that optimizes some $\nu$-strongly convex risk function $J(w)$ defined in terms of a loss function \cite{vap2000,tow2012b}:
\bs
\be
\label{eqn:wodef}
w^o \define \mbox{\rm argmin}_w \, J(w) \weq \mbox{\rm argmin}_w \, \Ex Q(w;\x_{i})	
\ee
In an asynchronous environment, the learner seeks $w^o$ by running the stochastic-gradient algorithm with random step-sizes:
\be
\label{eqn:SGDonlinelearner}
\addbox{\;\w_{i} \weq \w_{i-1} - \bmu(i) \, \nabla_{w\tran} Q(\w_{i-1};\x_{i}), \quad i\geq 0	\;}
\ee
The gradient noise vector is still given by
\be
\label{eqn:gradientnoiseonlinelearnerdef}
\s_{i}(\w_{i-1}) \weq \nabla_{w\tran} \, Q(\w_{i-1};\x_{i}) - \nabla_{w\tran}\,J(\w_{i-1})
\ee
Since $\nabla_w\, J(w^o) = 0$, and since the distribution of $\x_{i}$ is stationary, it follows that the covariance matrix of $\s_{i}(w^o)$ is constant and given by \be R_{s} = \Ex \nabla_{w\tran}\, Q(w^o;\x_{i}) \nabla_w\, Q(w^o;\x_{i})\ee The excess-risk measure that will result from this stochastic implementation is then given by \eqref{eqn:SGDERexpression} so that $\ER = O(\mu_x)$.
\es

\qd}
\end{example}
}

\section{Centralized Adaptation and Learning}
\label{sec:batch}
The discussion in the previous section establishes the mean-square stability of {\em stand-alone} adaptive agents for small step-sizes (Lemma \ref{lemma:msestability}), and provides expressions for their  MSD and ER metrics (Lemma \ref{lemma:SGDMSEperformance}). We now examine two situations involving a multitude of similar agents. In the first scenario, each agent senses data and analyzes it independently of the other agents. We refer to this mode of operation as  non-cooperative processing. In the second scenario, the agents transmit the collected data for processing at a fusion center. We refer to this mode of operation as centralized or batch processing. We motivate the discussion by considering first the case of mean-square-error costs. Subsequently, we extend the results to more general costs.

\subsection{Non-Cooperative MSE Processing}
\noindent Thus, consider {\em separate} agents, labeled $k = 1, 2, \dots, N$. Each agent, $k$, receives streaming data \bs \be \{\d_k(i),\u_{k,i}; i\geq 0\}\ee where we are using the subscript $k$ to index the data at agent $k$. We assume that the data at each agent satisfies the same statistical properties as in Example \ref{example:LMS}, and the same linear regression model \eqref{eqn:linearmodel} with a common $w^o$ albeit with noise $\v_k(i)$. We denote the statistical moments of the data at agent $k$  by
\be R_{u,k} = \Ex\u_{k,i}\tran \u_{k,i} > 0,\;\;\;\;\sigma_{v,k}^2 = \Ex\v_k^2(i)\ee
We further assume in this motivating example that the $R_{u,k}$ are uniform across the agents so that \be R_{u,k} \equiv R_u,\;\;\;k=1,2,\ldots,N\ee
 \es
 In this way, the cost $J_k(w) = \Ex (\d_k(i) -\u_{k,i} \w)^2$, which is associated with agent $k$, will satisfy a condition similar to \eqref{eqn:boundedHessian} with the corresponding parameters $\{\nu,\delta\}$ given by (cf. \eqref{eqn:boundedHessianmsedef}):
\bs
\be
\nu \weq 2\lambdamin(R_{u}), \quad
\delta \weq 2\lambda_{\max}(R_{u})
\ee
Now, if each agent runs the asynchronous LMS learning rule \eqref{eqn:LMS.async} to estimate $w^o$ on its own then, according to \eqref{eqn:LMSMSDexpressioncloseform}, each agent $k$ will attain an individual MSD level that is given by
\be
\label{eqn:MSDindexpression}
\addbox{\;\MSD_{\ind,k}^{\async} \approx \mu_x M \, \sigma_{v,k}^2, \quad k = 1, 2, \dots, N\;}
\ee
where we are further assuming that the parameter $\mu_x$ is uniform across the agents to enable a meaningful comparison.  Moreover, according to \eqref{eqn:alphadef}, agent $k$ will converge towards this level at a rate dictated by:
\be
\label{eqn:convergencerateind}
\addbox{\;\alpha_{\ind,k}^{\async} \approx 1 - 4 \mu_x \lambdamin(R_u)\;}
\ee
\es
The subscript ``ncop'' is used in \eqref{eqn:MSDindexpression} and \eqref{eqn:convergencerateind} to indicate that these expressions are for the non-cooperative mode of operation. It is seen from \eqref{eqn:MSDindexpression} that agents with noisier data (i.e., larger $\sigma_{v,k}^2$) will perform worse and have larger MSD levels than agents with cleaner data.
We are going to show in later sections that cooperation among the agents, whereby agents share information with their neighbors, can help enhance their individual performance levels.

\subsection{Centralized MSE Processing}
 Let us now contrast the above non-cooperative solution with a centralized implementation whereby, at every iteration $i$,  the $N$ agents transmit their raw data $\{\d_k(i),\u_{k,i}\}$ to a fusion center for processing.  In a \emph{synchronous} environment, once the fusion center receives the raw data, it can run a standard stochastic-gradient update of the form:
\bs
\be
\label{eqn:LMSbatchsync}
\w_i \weq \w_{i-1} + \mu \left[ \frac{1}{N} \sum_{k=1}^N 2 \u_{k,i}\tran (\d_k(i) - \u_{k,i}\w_{i-1}) \right]
\ee
where $\mu$ is the constant step-size, and the term multiplying $\mu$ can be seen to correspond to a sample average of several approximate gradient vectors. The analysis in \cite{ProcIEEE2014,NOW2014} showed that the MSD performance that results from this implementation is given by (using expression \eqref{eqn:MSDbatchrandomupdategeneral} with $H_k = 2R_u$ and $R_{s,k} = 4\sigma_{v,k}^2 R_u$):
\be
\label{eqn:MSDbatchsyncLMS}
\addbox{\;\MSD_{\cent}^{\sync} \approx  \frac{\mu M}{N} \left( \frac{1}{N} \sum_{k=1}^N \sigma_{v,k}^2\right)\;}
\ee
Moreover, using expression \eqref{eqn:alphabatchrandomupdateLMSdef} given further ahead, this centralized solution will converge towards the above MSD level at the same rate as the non-cooperative solution:
\be
\label{eqn:alphabatchsyncLMS}
\addbox{\;\alpha_{\cent}^{\sync} \approx 1 - 4 \mu \lambdamin(R_u)\;}
\ee
\es

In an asynchronous environment, there are now {\em several} random events that can interfere with the operation of the fusion center. Let us consider initially one particular random event that corresponds to the situation in which the fusion center may or may not update at any particular iteration (e.g., due to some power-saving strategy). In a manner similar to \eqref{eqn:randomstepsizeSGD}, we may represent this scenario by writing:
 \bs
\be
\label{eqn:LMSbatchrandomupdate}
\w_i \weq \w_{i-1} + \bmu(i) \left[ \frac{1}{N} \sum_{k=1}^N 2 \u_{k,i}\tran (\d_k(i) - \u_{k,i}\w_{i-1}) \right]
\ee
with a random step-size process, $\bmu(i)$; this process is again assumed to satisfy the conditions under Assumption \ref{assumption:randomstepsize}. The analysis in the sequel will show that the MSD performance that results from this implementation is given by:
\be
\label{eqn:MSDbatchrandomupdateLMS}
\MSD_{\cent}^{\async,1} \approx  \frac{\mu_x M}{N} \left( \frac{1}{N} \sum_{k=1}^N \sigma_{v,k}^2\right)
\ee
where we are adding the superscript ``1'' to indicate that this is a preliminary result pertaining to the particular asynchronous implementation \eqref{eqn:LMSbatchrandomupdate}. We will be generalizing this result soon, at which point we will drop the superscript ``1''. Likewise, using expression \eqref{eqn:alphabatchrandomupdateLMSdef} given further ahead,  this version of the asynchronous centralized solution converges towards the above MSD level at the same rate as the non-cooperative solution \eqref{eqn:convergencerateind} and the synchronous version \eqref{eqn:alphabatchsyncLMS}:
\be
\label{eqn:alphabatchrandomupdateLMS}
\alpha_{\cent}^{\async, 1} \approx 1 - 4 \mu_x \lambdamin(R_u)
\ee
\es
Observe from \eqref{eqn:MSDbatchsyncLMS} and \eqref{eqn:MSDbatchrandomupdateLMS} that the MSD level attained by the centralized solution is proportional to $1/N$ times the {\em average} noise power across all agents. This scaled average noise power can be larger than some of the individual noise variances and smaller than the remaining noise variances.
 This example shows that it does not generally hold that centralized stochastic-gradient implementations outperform all individual non-cooperative agents \cite{zhaoxiao2}.

Now, more generally, observe that in the synchronous batch processing case \eqref{eqn:LMSbatchsync}, as well as in the asynchronous implementation \eqref{eqn:LMSbatchrandomupdate}, the data collected from the various agents are equally aggregated with a weighting factor equal to $1/N$. We can incorporate a second type of random events besides the random variations in the step-size parameter. This second source of uncertainty involves the possibility of failure (or weakening) in the links connecting the individual agents to the fusion center (e.g., due to fading or outage, or perhaps due to the agents themselves deciding to enter into a sleep mode according to some power-saving policy). We can capture these possibilities by extending formulation \eqref{eqn:LMSbatchrandomupdate} in the following manner:

\bs
\be
\label{eqn:LMSbatchasync}
\addbox{\w_i= \w_{i-1} + \bmu(i) \left[ \sum_{k=1}^N 2 \bpi_k(i) \u_{k,i}\tran (\d_k(i) - \u_{k,i}\w_{i-1}) \right]}
\ee
where $\bm{\mu}(i)$ continues to be a random step-size process, but now the coefficients $\{ \bpi_k(i); k = 1, 2, \dots, N \}$ are new random fusion coefficients that satisfy:
\be
\label{eqn:fusioncoeffcond}
\addbox{\;\sum_{k=1}^{N} \bpi_k(i) \weq 1, \;\;\bpi_k(i) \ge 0\;}
\ee
\es
for every $i \ge 0$.

{%\small
\begin{assumption}[{\bf Conditions on random fusion coefficients}]
\label{assumption:randomfusion} {\em
It is assumed that $\bpi_k(i)$ is independent of $\bm{\mu}(i)$ and of other random variables in the learning algorithm. The fusion coefficients are also independent over time, namely, $\bpi_{k}(i)$ and $\bpi_{\ell}(j)$ are independent for any $i\neq j$. For the same time $i$, the coefficients $\{\bpi_k(i)\}$ are correlated over space in view of the first requirement in (\ref{eqn:fusioncoeffcond}). The mean and co-variance(s) of each $\bpi_k(i)$ are denoted by:
\bs
\begin{align}
\label{eqn:fusioncoeffmean}
\bar{\pi}_k & \define \Ex \bpi_k(i) \\
\label{eqn:fusioncoeffcov}
c_{\pi,k\ell} & \define \Ex (\bpi_k(i) - \bar{\pi}_k)(\bpi_\ell(i) - \bar{\pi}_\ell)
\end{align}
for all $k, \ell = 1, 2, \dots, N$ and all $i \ge 0$. When $k=\ell$, the scalar $c_{\pi, k k}$ corresponds to the variance of $\bpi_k(i)$ and, therefore, we shall also use the alternative notation $\sigma_{\pi,k}^2$ for this case:
\be
\sigma_{\pi,k}^2\define c_{\pi, k k}\;\geq 0
\ee
\es
}
\hfill \qd
\end{assumption}
}

\noindent It is straightforward to verify that the first and second-order moments of the coefficients $\{\bm{\pi}_k(i)\}$ satisfy:
\bs
\bq
&&\bar{\pi}_k \;\ge\; 0,\;\;\;\;\sum_{k=1}^N \bar{\pi}_k\;=\; 1\label{eqn:fusioncoeffmomentscond}\\
&&\sum_{k=1}^N c_{\pi, k\ell}\;=\;0,\;\;\;\mbox{\rm for any $\ell$}\\
&&\sum_{\ell=1}^N c_{\pi, k\ell} \;=\; 0,\;\;\;\mbox{\rm for any $k$}\label{eqn:fusioncoeffmomentscond.a}
\eq
\es

Note that the earlier asynchronous implementation \eqref{eqn:LMSbatchrandomupdate} is a special case of the more general formulation \eqref{eqn:LMSbatchasync} when the mean of the random fusion coefficients are chosen as $\bar{\pi}_k=1/N$ and their variances are set to zero, $\sigma_{\pi,k}^2=0$.  In order to enable a fair and meaningful comparison among the three centralized implementations \eqref{eqn:LMSbatchsync}, \eqref{eqn:LMSbatchrandomupdate}, and \eqref{eqn:LMSbatchasync}, we shall assume that the means of the fusion coefficients in the latest implementation are set to $\bar{\pi}_k=1/N$ (results for arbitrary mean values are listed in Example~\ref{example-1ZZYA} further ahead and appear in \cite{xiao2013BB}).
 We will show in the sequel that for this choice of $\bar{\pi}_k$, the MSD performance of implementation \eqref{eqn:LMSbatchasync} is given by:
\bs
\be
\label{eqn:MSDbatchasyncLMS}
\addbox{\;\MSD_{\cent}^{\async} \;\approx\;  \frac{\mu_x M}{N} \left[ \frac{1}{N} \sum_{k=1}^N \left(1 + N^2 \, \sigma^2_{\pi, k}\right) \sigma_{v,k}^2\right]\;}
\ee
Moreover, using expression \eqref{eqn:alphabatchasyncdef} given further ahead, this asynchronous centralized solution will converge towards the above MSD level at the same rate as the non-cooperative solution \eqref{eqn:convergencerateind} and the synchronous version \eqref{eqn:alphabatchsyncLMS}:
\be
\label{eqn:alphabatchasyncLMS}
\addbox{\;\alpha_{\cent}^{\async} \approx 1 - 4 \mu_x \lambdamin(R_u)\;}
\ee
\es

\subsection{Stochastic-Gradient Centralized Solution}\label{sec.jahdk.1llk13}
The previous two sections focused on mean-square-error costs. We now extend the conclusions to more general costs. Thus, consider a collection of $N$ agents, each with an individual convex cost function, $J_k(w)$. The objective is to determine the unique minimizer $w^o$ of the aggregate cost:
\be
\label{eqn:aggregatecostdef}
J^\glob(w) \define \sum_{k=1}^N J_k(w)
\ee
It is now the above aggregate cost, $J^\glob(w)$, that will be required to satisfy the conditions of Assumptions \ref{assumption:costfunctions} and \ref{assumption:smoothness} relative to some parameters $\{\nu_c,\delta_c,\tau_c\}$, with the subscript ``c'' used to indicate that these factors correspond to the centralized implementation. Under these conditions, the cost $J^\glob(w)$ will have a unique minimizer, which we continue to denote by $w^o$. We will not be requiring each individual cost, $J_k(w)$, to be strongly convex. It is sufficient for at least one of these costs to be strongly convex while the remaining costs can be convex; this condition ensures the strong convexity of $J^\glob(w)$.  Moreover, although some individual costs may not have a unique minimizer, we require in this exposition that $w^o$ is one of their minima so that all individual costs share a minimum at location $w^o$; the treatment in \cite{ProcIEEE2014,NOW2014,Chen2013} considers the more general case in which the minimizers of the individual costs $\{J_k(w)\}$ can be different and need not contain a common location, $w^o$.

There are many centralized solutions that can be used to determine the unique minimizer $w^o$ of \eqref{eqn:aggregatecostdef}, with some solution techniques being more powerful than other techniques. Nevertheless, we shall focus on centralized implementations of the stochastic-gradient type. The reason we consider the {\em same} class of stochastic gradient algorithms for non-cooperative, centralized, and distributed solutions in this chapter is to enable a {\em meaningful} comparison among the various implementations. Thus, we first consider a \emph{synchronous} centralized strategy of the following form:
\bs
\be
\label{eqn:batchsyncgeneral}
\w_i \weq \w_{i-1} - \mu\, \left[ \frac{1}{N} \sum_{k=1}^N \wh{\nabla_{w\tran}\,J}_k(\w_{i-1}) \right], \quad i \geq 0
\ee
with a constant step-size, $\mu$. When the fusion center employs random step-sizes, the above solution is replaced by:
\be
\label{eqn:batchrandomupdategeneral}
\w_i \weq \w_{i-1} - \bmu(i)\, \left[ \frac{1}{N} \sum_{k=1}^N \wh{\nabla_{w\tran}\,J}_k(\w_{i-1}) \right], \quad i \geq 0
\ee
where the process $\bmu(i)$ now satisfies Assumption \ref{assumption:randomstepsize}. More generally, the asynchronous implementation can employ random fusion coefficients as well such as:
\be
\label{eqn:batchasyncgeneral}
\addbox{\;\w_i \weq \w_{i-1} - \bmu(i)\, \sum_{k=1}^N \bpi_k(i) \wh{\nabla_{w\tran}\,J}_k(\w_{i-1}), \quad i \geq 0\;}
\ee
where the coefficients $\{\bpi_k(i)\}$ satisfy Assumption \ref{assumption:randomfusion} with means
\be
\label{eqn:uniformmeanpi}
\addbox{\;\bar{\pi}_k \;=\; 1/N\;}
\ee
\es

\subsection{Performance of Centralized Solution}
To examine the performance of the asynchronous implementation \eqref{eqn:batchasyncgeneral}--\eqref{eqn:uniformmeanpi}, we proceed in two steps. First, we identify the gradient noise that is present in the recursion; it is equal to the difference between the true gradient vector
 for the global cost, $J^{\rm glob}(w)$, defined by \eqref{eqn:aggregatecostdef} and its approximation. Second, we argue that \eqref{eqn:batchasyncgeneral} has a form similar to the single-agent stochastic-gradient algorithm \eqref{eqn:randomstepsizeSGD} and, therefore, invoke earlier results to write down performance metrics for the centralized solution \eqref{eqn:batchasyncgeneral}--\eqref{eqn:uniformmeanpi}.

 We start by introducing the individual gradient noise processes:
\bs
\be
\label{eqn:gradientnoisekdef}
\addbox{\;\s_{k,i}(\w_{i-1}) \define \wh{\nabla_{w\tran}\,J}_k(\w_{i-1}) - \nabla_{w\tran}\,J_k(\w_{i-1})\;}
\ee
for $k=1,2,\dots,N$. We assume that these noises satisfy conditions similar to Assumption \ref{assumption:gradientnoise} with parameters $\{\beta_{k}^2,\sigma_{s,k}^2,R_{s,k}\}$, i.e.,
\be
R_{s,k} \define \lim_{i\rightarrow\infty} \Ex \left[\s_{k,i}(w^o) \s_{k,i}\tran(w^o) | \Filtration_{i-1} \right]\ee
and
\bq
\hspace{-0.2cm} \Ex\left[ \s_{k,i}(\w_{i-1}) | \Filtration_{i-1} \right] & = &0 \\
\hspace{-0.2cm} \Ex\left[ \|{\s}_{k,i}(\w_{i-1}) \|^2 | \Filtration_{i-1} \right] & \leq &\beta_k^2 \, \|\wt{\w}_{i-1}\|^2 + \sigma_{s,k}^2
\eq
 Additionally, we assume that the gradient noise components across the agents are uncorrelated with each other:
\be
\label{eqn:uncorrelatedgradientnoise}
\Ex \left[{\s}_{k,i}(\w_{i-1}){\s}_{\ell,i}\tran(\w_{i-1}) \big| \Filtration_{i-1} \right] \weq 0, \quad \mbox{\rm all $k\neq \ell$}
\ee
\es
Using these gradient noise terms, it is straightforward to verify that recursion (\ref{eqn:batchasyncgeneral}) can be rewritten as:
\be
\w_i \;=\; \w_{i-1} - \frac{\bmu(i)}{N}\,\left[ \s_i(\w_{i-1})\;+\;\sum_{k=1}^N \nabla_{w\tran}\,J_k(\w_{i-1}) \right]\label{eqn:batchasyncgeneral2}\ee
where $\s_i(\w_{i-1})$ denotes the overall gradient noise; its expression is given by
\begin{align}
\label{eqn:gradientnoisebatchdef}
\s_{i}(\w_{i-1}) & \weq \sum_{k=1}^N \left[ N \bpi_k(i) \wh{\nabla_{w\tran}\,J}_k(\w_{i-1}) - \nabla_{w\tran}\,J_k(\w_{i-1}) \right]
\end{align}
Since iteration \eqref{eqn:batchasyncgeneral2} has the form of a stochastic gradient recursion with random update similar to \eqref{eqn:randomstepsizeSGD}, we can infer its mean-square-error behavior from Lemmas \ref{lemma:msestability} and \ref{lemma:SGDMSEperformance} if the noise process $\s_i(\w_{i-1})$ can be shown to satisfy conditions similar to Assumption \ref{assumption:gradientnoise} with some parameters $\{\beta_{c}^2,\sigma_{s}^2\}$. Indeed, starting from \eqref{eqn:gradientnoisebatchdef}, some algebra will show that
\bs
\bq
\Ex \left[ \s_{i}(\w_{i-1}) | \Filtration_{i-1} \right] &=& 0\label{eqn:martingalegradientnoisebatch}\\
\Ex \left[ \| \s_{i}(\w_{i-1}) \|^2 | \Filtration_{i-1} \right] & \leq&\beta_c^2 \| \wt{\w}_{i-1} \|^2 + \sigma_s^2
\eq
where
\bq
\beta_c^2 &\define& \sum_{k=1}^N \left[ \beta_k^2 + N^2 \sigma^2_{\pi,k} ( \beta_k^2 + \delta_c^2) \right]\label{eqn:betacandsigmasdef}\\
\sigma_s^2&\define& \sum_{k=1}^N (1 + N^2 \sigma^2_{\pi,k}) \sigma_{s,k}^2
\eq
\es
The following result now follows from Lemmas \ref{lemma:msestability} and \ref{lemma:SGDMSEperformance} \cite{xiao2013BB,NOW2014}.

{%\small
\begin{lemma}[{\bf Convergence of centralized solution}]
\label{lemma:convergencebatch}
Assume the aggregate cost \eqref{eqn:aggregatecostdef} satisfies the conditions under Assumption \ref{assumption:costfunctions} for some parameters $0 < \nu_c \leq \delta_c$. Assume further that the individual gradient noise processes defined by \eqref{eqn:gradientnoisekdef} satisfy the conditions under Assumption \ref{assumption:gradientnoise} for some parameters $\{ \beta_{k}^2, \sigma_{s,k}^2, R_{s,k} \}$, in addition to the orthogonality condition \eqref{eqn:uncorrelatedgradientnoise}. Let $\mu_o = 2 \nu_c / (\delta_c^2 + \beta_c^2)$. For any $\mu_x/N< \mu_o$, the iterates generated by the asynchronous centralized solution \eqref{eqn:batchasyncgeneral}--\eqref{eqn:uniformmeanpi} satisfy:
\bs
\be
\label{eqn:batcherrorrecurineq}
\Ex \|\wt{\w}_{i}\|^2 \;\leq\; \alpha \, \Ex \|\wt{\w}_{i-1}\|^2 + \sigma_{s}^2(\bar{\mu}^2+\sigma_{\mu}^2)/N^2\,
\ee
where the scalar $\alpha$ satisfies $0 \leq \alpha < 1$ and is given by
\be
\label{eqn:alphabatchasyncdef}
\alpha_{\rm cent}^{\rm asyn} \weq 1 - 2 \nu_c\,\left(\mu_x/N\right) + (\delta_c^2 + \beta_c^2) (\bar{\mu}^2+\sigma_{\mu}^2)/N^2
\ee
It follows from \eqref{eqn:batcherrorrecurineq} that for sufficiently small step-size parameter $\mu_x \ll 1$:
\be
\label{eqn:bounded2normbatch}
\limsup_{i\rightarrow\infty} \, \Ex \|\wt{\w}_{i}\|^2 \weq O(\mu_x)
\ee
\es
Moreover, under smoothness conditions similar to \eqref{eqn:LipschitzHessian} for $J^\glob(w)$ for some parameter $\tau_{c}\geq 0$, and similar to \eqref{eqn:LipschitzCovariance} for the individual gradient noise covariance matrices, it holds for small $\mu_x$ that:
\be
\label{eqn:MSDbatchasyncgeneral}
\mbox{\rm MSD}_{\rm cent}^{\rm asyn} \weq
\frac{\mu_x}{2N}\,\Tr\left[\left(\sum_{k=1}^N H_k\right)^{-1} \left(\sum_{k=1}^N (1 + N^2 \sigma^2_{\pi,k}) R_{s,k}\right)\right] + O\left(\mu_x^{3/2}\right)\ee
where $H_k=\nabla_w^2\,J_k(w^o)$. \hfill \qd
\end{lemma}
}

\bigskip

\noindent We can recover from the expressions in the lemma, performance results for the
particular asynchronous implementation described earlier by \eqref{eqn:batchrandomupdategeneral} by setting $\sigma_{\pi,k}^2=0$ so that
\bs\bq
\label{eqn:MSDbatchrandomupdategeneral}
\hspace{-0.4cm}\MSD_{\cent}^{\async,1} &\approx& \frac{\mu_x}{2N}\,\Tr\left[\left(\sum_{k=1}^N H_k\right)^{-1} \left(\sum_{k=1}^N R_{s,k}\right)\right]\\
\label{eqn:alphabatchrandomupdateLMSdef}
\hspace{-0.4cm}\alpha_{\cent}^{\async,1} &\approx& 1 - 2 \nu_c \mu_x/N
\eq
\es
\noindent It is seen from (\ref{eqn:MSDbatchasyncgeneral}) and (\ref{eqn:MSDbatchrandomupdategeneral}) that when the fusion center operates under the more general asynchronous policy \eqref{eqn:batchasyncgeneral}, the additional randomness in the fusion coefficients $\{\bpi_k(i)\}$ degrades the MSD performance relative to \eqref{eqn:MSDbatchrandomupdategeneral} due to the presence of the factor $1+N^2\sigma_{\pi,k}^2> 1$, i.e., to first-order in $\mu_x$ we have:
 \bs
 \be
\MSD^{\async, 1}\;<\;\MSD^{\async} \ee
In comparison, the added randomness due to the $\{\bpi_k(i)\}$ does not have significant impact on the convergence rate since it is straightforward to see that, to first order in the step-size parameter $\mu_x$:
\be
\alpha_{\cent}^{\async,1 } \approx  \alpha_{\cent}^{\async}
\ee
\es
Likewise, setting $\bar{\mu}_k=\mu$ and $\sigma_{\pi,k}^2=0$, the performance of the synchronous centralized solution \eqref{eqn:batchsyncgeneral} is obtained as a special case of the results in the lemma:
\bs\bq
\hspace{-0.2cm}\MSD_{\cent}^{\sync} &\approx& \frac{\mu}{2N}\,\Tr\left[\left(\sum_{k=1}^N H_k\right)^{-1} \left(\sum_{k=1}^N R_{s,k}\right)\right]\\
\hspace{-0.2cm}\alpha_{\cent}^{\sync}&\approx&1 - 2 \nu_c \mu/N
\eq
\es
These expressions agree with the performance results presented in \cite{ProcIEEE2014,NOW2014}.

{\small
\noindent \begin{example}[{\bf Case of general mean-values}]
{\rm Although not used in this chapter, we remark in passing that if the mean values $\{\bar{\pi}_k\}$ are not fixed at $1/N$, as was required by \eqref{eqn:uniformmeanpi}, then the MSD performance of the asynchronous
centralized solution \eqref{eqn:batchasyncgeneral} will instead be given by
\be
\label{eqn:MSDbatchasyncgeneraasxasxl}
\mbox{\rm MSD}_{\rm cent}^{\rm asyn} \weq
\frac{\mu_x N}{2}\,\Tr\left[\left(\sum_{k=1}^N H_k\right)^{-1} \left(\sum_{k=1}^N (\bar{\pi}_k^2 + \sigma_{\pi,k}^2) R_{s,k}\right)\right] + O\left(\mu_x^{3/2}\right)\ee

}

\qd
\label{example-1ZZYA}
\end{example}
}

\subsection{Comparison with Non-Cooperative Processing}
\label{subsec:comparisonbatchvsnoncoop}
We can now compare the performance of the asynchronous centralized solution \eqref{eqn:batchasyncgeneral} against the performance of  non-cooperative processing when agents act independently of each other and run the recursion:
\bs
\be
\label{eqn:SGDnoncoopk}
\addbox{\;\w_{k,i} \weq \w_{k,i-1} - \bmu(i)\,\wh{\nabla_{w\tran}\,J}_k(\w_{k,i-1}), \quad i\geq 0\;}
\ee
This comparison is {\em meaningful} when all agents share the same unique minimizer so that we can compare how well the individual agents are able to recover the same $w^o$ as the centralized solution. For this reason, we re-introduce the requirement that all individual costs $\{J_k(w)\}$ are $\nu$-strongly convex with a uniform parameter $\nu$. Since $J^\glob(w)$ is the aggregate sum of the individual costs,  then we can set the lower bound $\nu_c$ for the Hessian of $J^\glob(w)$ at $\nu_c = N \nu$. From expressions \eqref{eqn:alphadef} and \eqref{eqn:alphabatchrandomupdateLMSdef} we then conclude that, for a sufficiently small $\mu_x$, the convergence rates of the asynchronous non-cooperative solution \eqref{eqn:SGDnoncoopk} and the asynchronous centralized solution with random update \eqref{eqn:batchrandomupdategeneral} will be similar:
\bq
\alpha_{\cent}^{\async} &\approx& 1 - 2\nu_c\left(\mu_x/N\right)\;=\; 1 - 2 \nu \mu_x\; \approx \;\alpha_{\ind,k}^{\async}
\eq
Moreover, we observe from \eqref{eqn:SGDMSDexpression} that the average MSD level across $N$ non-cooperative asynchronous agents is given by
\be
\label{eqn:MSDnoncoopavg}
\MSD_{\ind, \av}^{\async} \approx \frac{\mu_x}{2N}\, \Tr\left[ \sum_{k=1}^{N} H_k^{-1} R_{s,k} \right]
\ee
so that comparing with \eqref{eqn:MSDbatchrandomupdategeneral}, some simple algebra allows us to conclude that, for small step-sizes and to first order in $\mu_x$:
\be
\label{eqn:compareMSDbatchrandomupdatevsnoncoopavg}
\MSD_{\cent}^{\async, 1} \;\leq\; \MSD_{\ind, \av}^{\async}
\ee
\es
That is, while the asynchronous centralized solution  \eqref{eqn:batchrandomupdategeneral}  with random updates need not outperform every individual non-cooperative agent in general, its performance outperforms the average performance across all non-cooperative agents.

The next example illustrates this result by considering the scenario where all agents have the same Hessian matrices at $w = w^o$, namely, $H_k \equiv H$ for $k = 1,2,\dots,N$. This situation occurs, for example, when the individual costs are identical across the agents, say, $J_k(w) \equiv J(w)$, as is common in machine learning applications. This situation also occurs for the mean-square-error costs we considered earlier in this section when the regression covariance matrices, $\{R_{u,k}\}$, are uniform across all agents, i.e., $R_{u,k} \equiv R_u$ for $k = 1,2,\dots,N$.  In these cases with uniform Hessian matrices $H_k$,  the example below establishes that the asynchronous centralized solution \eqref{eqn:batchrandomupdategeneral}  with random updates improves over the average MSD performance of the non-cooperative solution \eqref{eqn:SGDnoncoopk} by a factor of $N$.

{\small
\begin{example}[$\bm{N}$-{\bf fold improvement in MSD performance}]{\rm
Consider a collection of $N$ agents whose individual cost functions, $J_k(w)$, are $\nu$-strongly convex and are minimized at the same location $w = w^o$. The costs are also assumed to have identical Hessian matrices at $w = w^o$, i.e., $H_k \equiv H$. Then, using \eqref{eqn:MSDbatchrandomupdategeneral}, the MSD of the asynchronous centralized implementation \eqref{eqn:batchrandomupdategeneral} with random updates is given by:
\be
\label{eqn:MSDbatchrandomupdateNfoldnoncoop}
\MSD_{\cent}^{\async,1 } \approx \frac{1}{N} \left(\frac{\mu_x}{2N} \sum_{k=1}^{N} \Tr(H^{-1}  R_{s,k})\right) \approx \frac{1}{N}\, \MSD_{\ind,\av}^{\async}
\ee \hfill \qd}
\label{example-AAA}
\end{example}
}

{\small
\begin{example}[{\bf Random fusion can degrade MSD performance}]{\rm
Although the convergence rates of the asynchronous centralized solution \eqref{eqn:batchasyncgeneral} and the non-cooperative solution \eqref{eqn:SGDnoncoopk} agree to first-order in $\mu_{x}$, the relation between their MSD values is indefinite (contrary to \eqref{eqn:compareMSDbatchrandomupdatevsnoncoopavg}), as illustrated by the following example.

Consider the same setting of Example~\ref{example-AAA} and assume further that the variances of the random fusion coefficients  are uniform, i.e., $\sigma^2_{\pi,k} \equiv \sigma^2_{\pi}$. Then, using \eqref{eqn:MSDbatchasyncgeneral}, the MSD of the asynchronous centralized implementation \eqref{eqn:batchasyncgeneral} is given by
\bq
\MSD_{\cent}^{\async} &\approx &\frac{1 + N^2 \sigma^2_{\pi}}{N} \left(\frac{\mu_x}{2N} \sum_{k=1}^{N} \Tr(H^{-1} R_{s,k})\right)\;\approx \;\left(\frac{1 + N^2 \sigma^2_{\pi}}{N}\right)\, \MSD_{\ind,\av}^{\async}\label{eqn:MSDbatchasyncvsnoncoop}
\eq
Therefore, to first-order in $\mu_{\x}$, we find that
\be
\begin{cases}
\MSD_{\cent}^{\async} \leq \MSD_{\ind,\av}^{\async}, & \mbox{\rm if}\;\sigma^2_{\pi} \leq \frac{N-1}{N^2} \\
\MSD_{\cent}^{\async} > \MSD_{\ind,\av}^{\async}, & \mbox{\rm if}\;\sigma^2_{\pi} > \frac{N-1}{N^2} \\
\end{cases}
\ee
In other words, if the variance of the random fusion coefficients is large enough, then the centralized solution will generally have degraded performance relative to the non-cooperative solution (which is an expected result).

\hfill \qd}
\end{example}
}

{\small
\begin{example}[{\bf Fully-connected networks}]{\rm
In preparation for the discussion on networked agents, it is useful to describe one extreme situation where a collection of $N$ agents are fully connected to each other --- see Fig.~\ref{fig:fullyconnectedtopology}. In this case,  each agent is able to access the data from all other agents and, therefore, each individual agent can run a synchronous or asynchronous centralized implementation, say, one of the same form as \eqref{eqn:batchrandomupdategeneral}:
\be
\label{eqn:distributedfullyconnectedSGD}
\w_{k,i} \weq \w_{k,i-1} - \bm{\mu}(i) \left[ \frac{1}{N}\,\sum_{\ell=1}^N \wh{\nabla_{w\tran}\,J}_{\ell}(\w_{k,i-1}) \right], \quad i\geq 0
\ee

\bigskip
\bigskip
\bigskip
\bigskip
\bigskip
\bigskip

\begin{figure}[h]
\epsfxsize 9cm \epsfclipon
\begin{center}
\leavevmode \epsffile{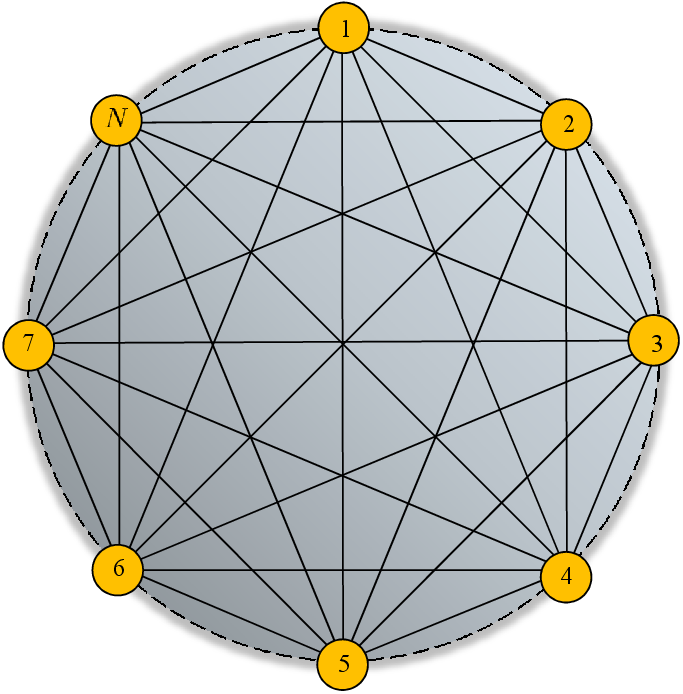}
\caption{{\small Example of a fully-connected network, where each agent can access  information from all other agents. Figure extracted  with permission from \cite{ProcIEEE2014}.}}\label{fig:fullyconnectedtopology}
\end{center}
\end{figure}

When this happens, each agent will attain the same performance level as that of the asynchronous centralized solution \eqref{eqn:batchrandomupdategeneral}. Two observations are in place. First, note from \eqref{eqn:distributedfullyconnectedSGD} that the information that agent $k$ is receiving from all other agents is their gradient vector approximations.  Obviously,  other pieces of information could be shared among the agents, such as their iterates $\{\w_{\ell,i-1}\}$. Second, note that the right-most term multiplying $\bm{\mu}(i)$ in \eqref{eqn:distributedfullyconnectedSGD}  corresponds to a convex combination of the approximate gradients from the various agents, with the combination coefficients being uniform and all equal to $1/N$. In general, there is no need for these combination weights to be identical. Even more importantly, agents do not need to have access to information from all other agents in the network. We are going to see in the sequel that interactions with a limited number of neighbors is  sufficient for the agents to attain performance that is comparable to that of the centralized solution.

Figure~\ref{fig:6topologies} shows a simple selection of connected topologies for five agents. The leftmost panel corresponds to the non-cooperative case and the rightmost panel corresponds to the fully-connected case. The panels in between illustrate some other topologies. In the coming sections, we are going to present results that allow us to answer useful questions about such networked agents such as: (a) Which topology has best performance in terms of mean-square error and convergence rate? (b) Given a connected topology, can it be made to approach  the performance of the centralized solution? (c) Which aspects of the topology influence performance? (d)
 Which aspects of the combination weights (policy) influence performance? (e) Can different topologies deliver similar performance levels? (f) Is cooperation always beneficial? and (g) If the individual agents are able to solve the inference task individually in a stable manner, does it follow that the connected network will remain stable regardless of the topology and regardless of the cooperation strategy?\\
 
 \bigskip
 \bigskip
 \bigskip

\begin{figure}[h]
\epsfxsize 11cm \epsfclipon
\begin{center}
\leavevmode \epsffile{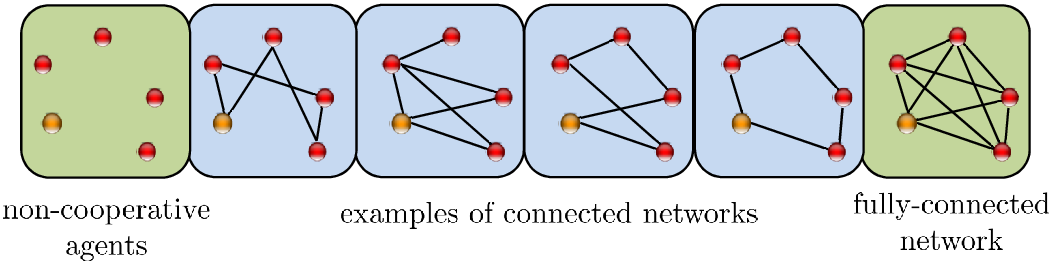}
\caption{{\small Examples of connected networks, with the left-most panel representing a collection of non-cooperative agents. Figure extracted  with permission from \cite{ProcIEEE2014}.}}\label{fig:6topologies}
\end{center}
\end{figure}

\qd}
\end{example}
}

\section{Synchronous Multi-Agent Adaptation and Learning}
\label{sec:network}
In this section, we describe distributed strategies of the consensus (e.g., \cite{bookchapter,TsiBB86,moura10} and \cite{Nedic2010}--\cite{braca08}) and diffusion (e.g., \cite{ProcIEEE2014,NOW2014,sayedSPM, bookchapter,Chen10b,  Lopes08, Cattivelli10, yusayed12bb}) types. These strategies rely solely on localized interactions among neighboring agents, and they can be used to seek the minimizer of \eqref{eqn:aggregatecostdef}. We first describe the network model.

\subsection{Strongly-Connected Networks}
\begin{figure}[h]
\epsfxsize 10cm \epsfclipon
\begin{center}
\leavevmode \epsffile{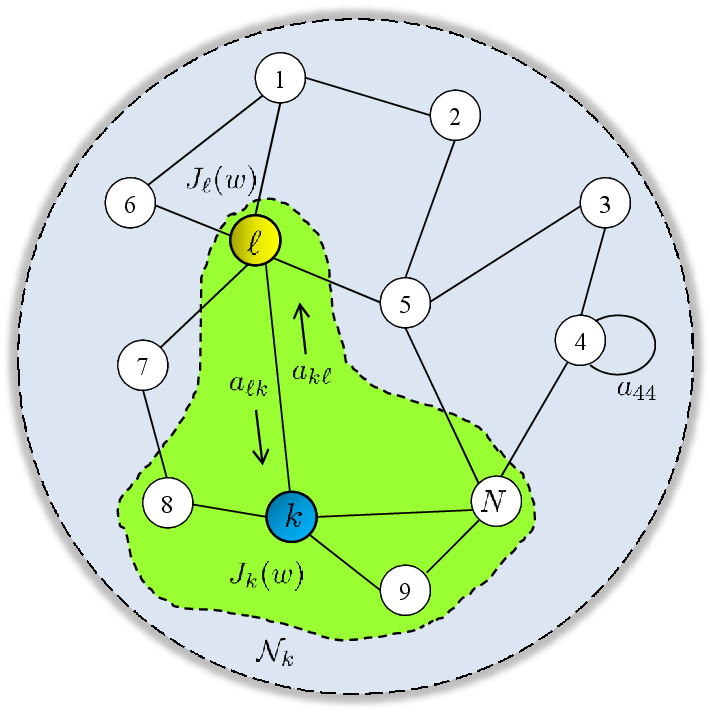}
\caption{{\small Agents that are linked by edges can share information. The neighborhood of agent $k$ is marked by the highlighted area. Figure extracted  with permission from \cite{ProcIEEE2014}.}}
\label{fig:networktopology}
\end{center}
\end{figure}

 Figure \ref{fig:networktopology} shows an example of a network consisting of $N$ connected agents, labeled $k = 1, 2, \dots, N$. Following the presentation from \cite{NOW2014,bookchapter}, the network is represented by a graph consisting of $N$ vertices (representing the agents) and a set of edges connecting the agents to each other. An edge that connects an agent to itself is called a self-loop. The neighborhood of an agent $k$ is denoted by $\Ncal_k$ and it consists of all agents that are connected to $k$ by an edge.  Any two neighboring agents $k$ and $\ell$ have the ability to share information over the edge connecting them.

\bigskip

We assume an undirected graph so that if agent $k$ is a neighbor of agent $\ell$, then agent $\ell$ is also a neighbor of agent $k$. We assign a pair of nonnegative scaling weights, $\{a_{k\ell},a_{\ell k}\}$, to the edge connecting $k$ and $\ell$. The scalar $a_{\ell k}$ is used by agent $k$ to scale the data it receives from agent $\ell$; this scaling can be interpreted as a measure of the confidence level that agent $k$ assigns to its interaction with agent $\ell$. Likewise, $a_{k \ell}$ is used by agent $\ell$ to scale   the data it receives from agent $k$.  The weights $\{a_{k\ell}, a_{\ell k}\}$ can be different so that the exchange of information between the neighboring agents $\{k,\ell\}$ need not be symmetrical. One or both weights can also be zero.

A network is said to be connected if paths with nonzero scaling weights can be found linking any two distinct agents in {\em both} directions, either directly when they are neighbors or by passing through intermediate agents when they are not neighbors. In this way, information can flow in both directions between any two agents in the network,  although the forward path from an agent $k$ to some other agent $\ell$ need not be the same as the backward path from $\ell$ to $k$.  A strongly-connected network is a connected network with at least one non-trivial self-loop, meaning that $a_{kk}>0$ for some agent $k$.

The strong connectivity of a network translates into a useful property on the combination weights. Assume we collect the coefficients $\{a_{\ell k}\}$ into an $N\times N$ matrix $A = [a_{\ell k}]$, such that the entries on the $k$-th column of $A$ contain the coefficients used by agent $k$ to scale data arriving from its neighbors $\ell \in \Ncal_k$; we set $a_{\ell k} = 0$ if $\ell \notin \Ncal_k$.  We refer to $A$ as the {\em combination} matrix or policy. It turns out that combination matrices that correspond to strongly-connected networks are {\em primitive} --- an $N\times N$ matrix $A$ with nonnegative entries is said to be primitive if  there exists some finite integer $n_o > 0$ such that all entries of $A^{n_o}$ are strictly positive \cite{NOW2014,bookchapter,horn}.

\subsection{Distributed Optimization}\label{sec.diakd.opt}
Network cooperation can be exploited to solve  adaptation, learning, and optimization problems in a decentralized manner in response to streaming data. To explain how cooperation can be achieved, we start by associating with each agent $k$ a twice-differentiable cost function $J_k(w) : \real^{M\times 1} \mapsto \real$. The objective of the network of agents is to seek the unique minimizer of the aggregate cost function, $J^\glob(w)$, defined by \eqref{eqn:aggregatecostdef}. Now, however, we seek a {\em distributed} (as opposed to a centralized) solution. In a distributed implementation, each agent $k$ can only rely on its own data and on data from its neighbors.

We continue to assume that $J^\glob(w)$ satisfies the conditions of Assumptions \ref{assumption:costfunctions} and \ref{assumption:smoothness} with parameters $\{\nu_d, \delta_d, \tau_d\}$, with the subscript ``d'' now used to indicate that these parameters are related  to the distributed implementation. Under these conditions, the cost $J^\glob(w)$ will have a unique minimizer, which we continue to denote by $w^o$. For simplicity of presentation, we will also assume in the remainder of this chapter that the individual costs $J_k(w)$ are strongly convex as well. These costs can be distinct across the agents or they can all be identical, i.e., $J_k(w) \equiv  J(w)$ for $k=1,2,\dots,N$; in the latter situation, the problem of minimizing  \eqref{eqn:aggregatecostdef}  would correspond to the case in which the agents  work together to optimize the same cost function. If we let $w_k^o$ denote the minimizer of $J_k(w)$, we continue to assume for this exposition that each $J_k(w)$ is also minimized at $w^o$:
\be
w_{k}^o\equiv w^o,\;\;\;k=1,2,\ldots,N
\label{kajdhdg6713.13}\ee
The case where the individual costs are only convex and need not be strongly convex is discussed in \cite{NOW2014,ProcIEEE2014,chen2013}; most of the results and conclusions continue to hold but the derivations become more technical. Likewise, the case in which the individual costs need not share minimizers is discussed in these references. In that case, as was already shown in \cite{Chen2013}, the iterates $\w_{k,i}$ by the individual agents will not approach the minimizer $w^o$ of (\ref{eqn:aggregatecostdef}), but rather the minimizer $w^{\star}$ of a weighted aggregate function with some positive scaling weights $\{q_k\}$. This function was defined in Eq. (57b) of \cite{ProcIEEE2014}. It was also explained in \cite{ProcIEEE2014} how this convergence property adds a useful degree of freedom to the operation of the network, and how it can be exploited advantageously to steer the network to converge to desirable limit points, including to $w^o$, through the selection of the combination policy $A$ (which determines the scaling weights $\{q_k\}$ in the weighted aggregate cost).

It is, nevertheless,  sufficient for the purposes of this chapter to continue with the case (\ref{kajdhdg6713.13}).  There are many important situations in practice where the minimizers of individual costs coincide with each other. For instance, examples abound where agents need to work cooperatively to attain a common objective such as tracking a target, locating a food source, or evading a predator (see, e.g., \cite{Tu11,catmay11,sayedSPM}). This scenario is also common in machine learning problems \cite{Bish2007,theo2008,Dekel8,Agarwal8,Predd06,Tow11} when data samples at the various agents are generated by a common distribution parameterized by some vector, $w^o$. One such situation is illustrated in the next example.

{\small
\begin{example}[{\bf Mean-square-error (MSE) networks}]
\label{example:MSEnetwork}{\rm
Consider the same setting of Example \ref{example:LMS} except that we now have $N$ agents observing streaming data $\{\d_k(i),\u_{k,i}\}$ that satisfy the regression model \eqref{eqn:linearmodel} with regression covariance matrices $R_{u,k} = \Ex\u_{k,i}\tran \u_{k,i} > 0$ and with the same unknown $w^o$, i.e.,
\bs
\be
\label{eqn:linearmodelk}
\d_k(i) \weq \u_{k,i} w^o + \v_k(i)
\ee
The individual mean-square-error costs are defined  by \be J_k(w) = \Ex(\d_k(i) - \u_{k,i} w)^2\ee
and are strongly convex in this case, with the minimizer of each $J_k(w)$ occurring at
\be
\label{eqn:wokdef}
w_k^o \define R_{u,k}^{-1} r_{du,k}, \quad k = 1,2,\dots,N
\ee
\es
If we multiply both sides of \eqref{eqn:linearmodelk} by $\u_{k,i}\tran$ from the left, and take expectations, we find that $w^o$ satisfies \be r_{du,k} = R_{u,k} w^o\ee This relation shows that the unknown $w^o$ from \eqref{eqn:linearmodelk} satisfies the same expression as $w_k^o$ in \eqref{eqn:wokdef}, for any $k = 1, 2, \dots, N$, so that we must  have $w^o = w_k^o$. Therefore, this example amounts to a situation where all costs $\{J_k(w)\}$ attain their minima at the same location, $w^o$.

We shall use the network model of this example to illustrate other results in the chapter. For ease of reference, we shall refer to strongly-connected networks with agents receiving data according to model \eqref{eqn:linearmodelk} and seeking to estimate $w^o$ by adopting the mean-square-error costs $J_k(w)$ defined above, as {\em mean-square-error {\rm (MSE)} networks}. We assume for these networks that the measurement noise process $\v_k(i)$ is temporally white and independent over space so that \be \Ex\v_k(i)\v_{\ell}(j) = \sigma_{v,k}^2 \delta_{k,\ell} \delta_{i,j}\ee in terms of the Kronecker delta $\delta_{k,\ell}$. Likewise, we assume that the regression data $\u_{k,i}$ is temporally white and independent over space so that \be \Ex\u_{k,i}\tran \u_{\ell,j} = R_{u,k} \delta_{k,\ell} \delta_{i,j}\ee Moreover, the measurement noise $\v_k(i)$ and the regression data $\u_{\ell,j}$ are independent of each other for all $k,\ell,i,j$. These statistical conditions help facilitate the analysis of such networks.

\qd}
\label{example-9A}
\end{example}
}

In the next subsections we list {\em synchronous} distributed algorithms of the consensus and diffusion types for the optimization of $J^{\rm glob}(w)$. We only list the algorithms here; for motivation and justifications, the reader may refer to the treatments in \cite{ProcIEEE2014,NOW2014}. Moreover, Sec.~V.D of \cite{ProcIEEE2014} provides commentary on several other related works in the literature, in addition to the history and evolution of the consensus and diffusion strategies.

\subsection{Synchronous Consensus Strategy}
\label{subsec:consensus}
 Let $\w_{k,i}$ denote the iterate that is available at agent $k$ at iteration $i$; this iterate serves as the estimate for $w^o$. The consensus iteration at each agent $k$ is described by the following construction
 (see, e.g., \cite{TsiBB86,moura10,Ber97,NedicOzdal,dimakis10,braca08}):
 \be
\label{eqn:consensus1step}
\addbox{\;\w_{k,i} \weq \sum_{\ell \in \Ncal_k} a_{\ell k} \, \w_{\ell,i-1} - \mu_k\, \wh{\nabla_{w\tran}\,J}_k(\w_{k,i-1})\;}
\ee
where the $\{\mu_k\}$ are individual step-size parameters,  and where the combination coefficients $\{a_{\ell k}\}$ that appear in \eqref{eqn:consensus1step} are nonnegative scalars that are required to satisfy the following conditions for each agent $k=1,2,\dots,N$:
\bs
\be
\label{eqn:conditionalk}
\addbox{\;a_{\ell k} \geq 0, \quad \sum_{\ell=1}^N a_{\ell k} = 1, \;\;\;\; a_{\ell k} = 0~\mathrm{if}~\ell \notin \Ncal_{k}\;}
\ee
Condition \eqref{eqn:conditionalk} implies that the combination matrix $A=[a_{\ell k}]$ satisfies \be A\tran \one = \one\ee where $\one$ denotes the vector with all entries equal to one. We say that $A$ is left-stochastic. One useful property of left-stochastic matrices is that their spectral radius is equal to one \cite{bookchapter,horn,golub,  BermanPF, Pillai05}: \be \rho(A) = 1\ee  \es
An equivalent representation that is useful for later analysis is to rewrite the consensus iteration \eqref{eqn:consensus1step} as shown in the following listing, where the intermediate iterate that results from the neighborhood combination is denoted by $\bpsi_{k,i-1}$. Observe that the gradient vector in the consensus implementation \eqref{eqn:consensus2steps} is evaluated at $\w_{k,i-1}$ and not $\bpsi_{k,i-1}$.

\be
\label{eqn:consensus2steps}
\begin{array}{l}
\hline
\textrm{{\bf Consensus strategy for distributed adaptation}} \\
\hline
\; \mbox{\rm {\bf for} each time instant $i\geq 0$:} \\
\; \;\; \mbox{each agent $k=1,2,\dots,N$ performs the update:} \\
\; \;\; \qquad \left\{ \begin{aligned}
\bpsi_{k,i-1} & \weq \displaystyle \sum_{\ell \in \Ncal_k} a_{\ell k}\; \w_{\ell,i-1} \\
\w_{k,i}      & \weq \displaystyle \bpsi_{k,i-1} - \mu_k\, \wh{\nabla_{w\tran}\,J}_k \left(\w_{k,i-1}\right) \\
\end{aligned} \right. \\
\; \mbox{\rm {\bf end}} \\
\hline
\end{array}
\ee

\smallskip

 We remark that one way to motivate the consensus update \eqref{eqn:consensus1step} is to start from the non-cooperative step \eqref{eqn:SGDnoncoopk} and replace the first iterate $\w_{k,i-1}$ by the convex combination used in \eqref{eqn:consensus1step}.

{\small
\begin{example}[{\bf Consensus LMS networks}]{\rm
For the MSE network of Example \ref{example:MSEnetwork}, the consensus strategy reduces to:
\be
\label{eqn:consensusLMS}
\left\{ \begin{aligned}
\bpsi_{k,i-1} & \weq \displaystyle \sum_{\ell \in \Ncal_k} a_{\ell k}\; \w_{\ell,i-1} \\
\w_{k,i} & \weq \displaystyle \bpsi_{k,i-1} + 2\mu_k\u_{k,i}\tran [\d_k(i)-\u_{k,i}\w_{k,i-1}]
\end{aligned} \right.
\ee
\qd}\label{example-consensus}
\end{example}
}

\subsection{Synchronous Diffusion Strategies}
\label{subsec:diffusion}
There is an inherent asymmetry in the consensus construction. Observe from the computation of $\w_{k,i}$ in \eqref{eqn:consensus2steps} that the update starts from $\bm{\psi}_{k,i-1}$ and corrects it by the approximate gradient vector evaluated at $\w_{k,i-1}$ (and not at $\bm{\psi}_{k,i-1}$).  This {\em asymmetry}  will be shown later, e.g., in Example~\ref{example-consensus.2}, to be problematic when the consensus strategy is used for adaptation and learning over networks. This is because the asymmetry can cause an unstable growth in the state of the network \cite{yusayed12bb} --- see also the explanations in \cite{ProcIEEE2014} and \cite{NOW2014}[Sec. 10.6].  Diffusion strategies remove the asymmetry problem.\\

\noindent \underline{\em Combine-then-Adapt (CTA) Diffusion}. In the CTA formulation of the diffusion strategy, the {\em same} iterate $\bpsi_{k,i-1}$ is used to compute $\w_{k,i}$, thus leading to description \eqref{eqn:diffusionCTA} where the gradient vector is evaluated at $\bpsi_{k,i-1}$ as well. The reason for the name ``Combine-then-Adapt'' is that the first step in \eqref{eqn:diffusionCTA} involves a combination step, while the second step involves an adaptation step. The reason for the qualification ``diffusion'' is that the use of $\bpsi_{k,i-1}$ to evaluate the gradient vector allows information to diffuse more thoroughly through the network. This is because information is not only being diffused through the aggregation of the neighborhood iterates, but also through the evaluation of the gradient vector at the aggregate state value.

\bs
\be
\label{eqn:diffusionCTA}
\begin{array}{l}
\hline
\textrm{{\bf Diffusion strategy for distributed adaptation (CTA)}} \\
\hline
\; \mbox{\rm {\bf for} each time instant $i\geq 0$:} \\
\; \;\; \mbox{each agent $k=1,2,\dots,N$ performs the update:} \\
\; \qquad \left\{ \begin{aligned}
\bpsi_{k,i-1} & \weq \displaystyle \sum_{\ell \in \Ncal_k} a_{\ell k}\; \w_{\ell,i-1} \\
\w_{k,i} & \weq \displaystyle \bpsi_{k,i-1} - \mu_k\, \wh{\nabla_{w\tran}\,J}_k\left(\bpsi_{k,i-1}\right)
\end{aligned} \right. \\
\; \mbox{\rm {\bf end}} \\
\hline
\end{array}
\ee

\bigskip

\noindent \underline{\em Adapt-then-Combine (ATC) Diffusion}. A similar implementation can be obtained by switching the order of the combination and adaptation steps in \eqref{eqn:diffusionCTA}, as shown in the listing \eqref{eqn:diffusionATC}. The structure of the CTA and ATC strategies are fundamentally identical: the difference lies in which variable we choose to correspond to the updated iterate $\w_{k,i}$. In ATC, we choose the result of the {\em combination} step to be $\w_{k,i}$, whereas in CTA we choose the result of the {\em adaptation} step to be $\w_{k,i}$.

\be
\label{eqn:diffusionATC}
\begin{array}{l}
\hline
\textrm{{\bf Diffusion strategy for distributed adaptation (ATC)}} \\
\hline
\; \mbox{\rm {\bf for} each time instant $i\geq 0$:} \\
\; \;\; \mbox{each agent $k=1,2,\dots,N$ performs the update:} \\
\; \qquad \left\{ \begin{aligned}
\bpsi_{k,i}  & \weq \displaystyle \w_{k,i-1} - \mu_k \,\wh{\nabla_{w\tran}\,J}_k(\w_{k,i-1}) \\
\w_{k,i} & \weq \displaystyle \sum_{\ell \in \Ncal_k} a_{\ell k}\; \bpsi_{\ell,i} \\
\end{aligned} \right. \\
\; \mbox{\rm {\bf end}} \\
\hline
\end{array}
\ee
\es

\bigskip
One main motivation for the introduction of the diffusion strategies (\ref{eqn:diffusionCTA}) and (\ref{eqn:diffusionATC}) is the fact that they enable {\em single} time-scale distributed learning from {\em streaming} data under {\em constant} step-size adaptation and in a stable manner \cite{xiao2013,Chen2013,Lopes08,Cattivelli10} and \cite{Lopes06}--\cite{Cattivelli08} --- see also \cite{NOW2014}[Chs.~9--11]; the diffusion strategies further allow $A$ to be left-stochastic, which permit larger modes of cooperation than doubly-stochastic policies. The CTA diffusion strategy (\ref{eqn:diffusionCTALMS}) was first introduced for mean-square-error estimation problems in \cite{Lopes08,Lopes06,Sayed07,Lopes07a}. The ATC diffusion structure (\ref{eqn:diffusionATCLMS}), with adaptation preceding combination,  appeared in the work \cite{Cattivelli07}
on adaptive distributed least-squares schemes and also in the
works \cite{Cattivelli10,Cattivelli08,Cattivelli08c,Cattivelli08a} on
distributed mean-square-error and state-space estimation methods. The CTA structure (\ref{eqn:diffusionCTA}) with an iteration dependent step-size that decays to zero,
$\mu(i)\rightarrow 0$, was employed in \cite{Sriv2011,ramdistributed,leenedic} to solve
distributed optimization problems that require all agents
to reach agreement. The ATC form (\ref{eqn:diffusionATC}), also with an iteration dependent
sequence $\mu(i)$ that decays to zero, was employed in \cite{bianchi,Stankovic11} to ensure almost-sure convergence and agreement among agents.

{\small
\begin{example}[{\bf Diffusion LMS networks}]
\label{example:diffusionLMSnetworks}{\rm
For the MSE network of Example \ref{example:MSEnetwork}, the ATC and CTA diffusion strategies reduce to:
\bs
\be
\label{eqn:diffusionCTALMS}
\left\{ \begin{aligned}
\bpsi_{k,i-1} & = \displaystyle \sum_{\ell \in \Ncal_k} a_{\ell k}\;\w_{\ell,i-1} \qquad (\mbox{\rm CTA diffusion}) \\
\w_{k,i} & = \displaystyle \bpsi_{k,i-1} + 2\mu_k \u_{k,i}\tran \left[\d_{k}(i) - \u_{k,i} \bpsi_{k,i-1}\right]\end{aligned} \right.
\ee
and
\be
\label{eqn:diffusionATCLMS}
\left\{ \begin{aligned}
\bpsi_{k,i} & = \displaystyle \w_{k,i-1} + 2\mu_k \u_{k,i}\tran \left[\d_{k}(i)-\u_{k,i}\w_{k,i-1}\right] \\
\w_{k,i} & = \displaystyle \sum_{\ell \in \Ncal_k} a_{\ell k}\; \bpsi_{\ell,i}\quad\;\;\qquad (\mbox{\rm ATC diffusion})
\end{aligned} \right.
\ee
\es
\qd}\label{example-diffusion}
\end{example}
}

{\small
\begin{example}[{\bf Diffusion logistic regression}]{\rm
We revisit the pattern classification problem from Example \ref{example:loglossrisk}, where we consider a collection of $N$ networked agents cooperating with each other to solve the logistic regression problem. Each agent receives streaming data $\{ \bgamma_k(i),\h_{k,i}\}$, where the variable $ \bgamma_k(i)$ assumes the values $\pm 1$ and designates the class that the feature vector $\h_{k,i}$ belongs to. The objective is to use the training data to determine the vector $w^o$ that minimizes the cost
\be
\label{eqn:logisticlosskdef}
J_k(w) \define \frac{\rho}{2}\|w\|^2 + \Ex\left\{\ln \left[1 + e^{-\bgamma_k(i) \h_{k,i}\tran w}\right]\right\}
\ee
under the assumption of joint wide-sense stationarity over the random data. It is straightforward to verify that the ATC diffusion strategy \eqref{eqn:diffusionATC} reduces to the following form in this case:
\be
\label{eqn:diffusionlogisticregression}
\left\{ \begin{aligned}
\bpsi_{k,i} & = \displaystyle (1 - \rho \mu_k) \w_{k,i-1} + \mu_k \, \left(\frac{ \bgamma_k(i) }{1 + e^{ \bgamma_k(i)\h_{k,i}\tran \w_{k,i-1}}} \right) \h_{k,i} \\
\w_{k,i} & = \displaystyle \sum_{\ell \in \Ncal_k} a_{\ell k}\; \bpsi_{\ell,i}
\end{aligned} \right.
\ee
\qd}
\end{example}
}

\section{Asynchronous Multi-Agent Adaptation and Learning}
\label{sec:asynchronous}
There are various ways by which asynchronous events can be introduced into the operation of a distributed strategy. Without loss in generality, we illustrate the model for asynchronous operation by describing it for the ATC diffusion strategy (\ref{eqn:diffusionATC}); similar constructions apply to CTA diffusion (\ref{eqn:diffusionCTA}) and consensus \eqref{eqn:consensus2steps}.

\subsection{Asynchronous Model}
In a first instance, we model the step-size parameters as random variables and replace \eqref{eqn:diffusionATC} by:
\be
\label{eqn:diffusionATCrandomupdate}
\left\{ \begin{aligned}
\bpsi_{k,i}  & \weq \displaystyle \w_{k,i-1} - \bmu_k(i) \,\wh{\nabla_{w\tran}\,J}_k(\w_{k,i-1}) \\
\w_{k,i} & \weq \displaystyle \sum_{\ell \in \Ncal_k} a_{\ell k}\; \bpsi_{\ell,i} \\
\end{aligned} \right.
\ee
In this model, the neighborhoods and the network topology remain fixed and only the $\bmu_k(i)$ assume random values. The step-sizes can vary across the agents and, therefore, their means and variances become agent-dependent. Moreover, the step-sizes across agents can be correlated with each other. We therefore denote the first and second-order moments of the step-size parameters by:
\bs
\bq
\bar{\mu}_k & \define& \Ex \bmu_k(i) \label{eqn:randstepsizemeanentry} \\
\sigma_{\mu,k}^2&\define&\Ex (\bmu_k(i)-\bar{\mu}_k)^2\\
c_{\mu,k\ell} &\define & \Ex (\bmu_k(i) - \bar{\mu}_k)(\bmu_\ell(i) - \bar{\mu}_\ell)\label{eqn:randstepsizecoventry}
\eq
When $\ell = k$, the scalar $c_{\mu,kk}$ coincides with the variance of $\bmu_k(i)$, i.e., $c_{\mu, kk}=\sigma_{\mu,k}^2\geq 0$. On the other hand, if the step-sizes across the agents happen to be uncorrelated, then $c_{\mu,k\ell}=0$ for $k\neq \ell$.
\es

More broadly, we can allow for random variations in the neighborhoods (and, hence, in the network structure),  and random variations in the combination coefficients as well. We capture this more general asynchronous implementation by writing:
\be
\label{eqn:diffusionATCasync}
\addbox{\;\left\{ \begin{aligned}
\bpsi_{k,i} & \weq \w_{k,i-1} - \bmu_k(i) \wh{\nabla_{w\tran}\,J}_k(\w_{k,i-1}) \\
\w_{k,i} & = \sum_{\ell\in{\small \bm{\Ncal}_{k,i}}} \a_{\ell k}(i)\;\bpsi_{\ell,i}
\end{aligned} \right.\;}
\ee
where the combination coefficients $\{\a_{\ell\,k}(i)\}$ are now \emph{random} and, moreover, the symbol $\bm{\Ncal}_{k,i}$ denotes the \emph{randomly}-changing neighborhood of agent $k$ at time $i$. These neighborhoods become random because the random variations in the combination coefficients can turn links on and off depending on the values of the $\{\a_{\ell,k}(i)\}$. We continue to require the combination coefficients $\{\a_{\ell k}(i)\}$ to satisfy the same structural constraint as given before by
 \eqref{eqn:conditionalk}, i.e.,
\be
\label{eqn:randomtopologyconstraints}
\addbox{\;\sum_{\ell\in{\small \bm{\Ncal}_{k,i}}} \a_{\ell k}(i) = 1, \;\;\mbox{and}\;\;
\begin{cases}
\a_{\ell k}(i) > 0, & \mbox{if} \; \ell \in \bm{\Ncal}_{k,i} \\
\a_{\ell k}(i) = 0, & \mbox{otherwise}
\end{cases}\;}
\ee
Since these coefficients are now random, we denote their first and second-order moments by:
\bs
\bq
\hspace{-0.3cm}\bar{a}_{\ell k} & \define& \Ex \,\a_{\ell k}(i) \label{eqn:randcombinemeanelement}
\\
\hspace{-0.3cm}\sigma^2_{a,\ell k} &\define &\Ex (\a_{\ell k}(i) - \bar{a}_{\ell k})^2\label{eqn:randcombinecoventry.2}\\
\hspace{-0.3cm}c_{a,\ell k,nm} &\define &\Ex (\a_{\ell k}(i) - \bar{a}_{\ell k})(\a_{nm}(i) - \bar{a}_{nm})\label{eqn:randcombinecoventry}
\eq
\es
When $\ell = n$ and $k=m$, the scalar $c_{a,\ell k,nm}$ coincides with the variance of $\a_{\ell k}(i)$, i.e., $c_{a,\ell k,nm}=\sigma_{a,\ell k}^2\geq 0$.

{\small
\begin{example}[{\bf Asynchronous diffusion LMS networks}]
\label{example:diffusionLMSATCnetworkasync}{\rm
For the MSE network of Example \ref{example:MSEnetwork}, the ATC diffusion strategy \eqref{eqn:diffusionATCrandomupdate} with random update reduces to
\bs
\be
\label{eqn:diffusionATCLMSrandomupdate}
\left\{ \begin{aligned}
\bpsi_{k,i} & = \displaystyle \w_{k,i-1} + 2\bmu_k(i) \u_{k,i}\tran \left[\d_{k}(i)-\u_{k,i}\w_{k,i-1}\right] \\
\w_{k,i} & = \displaystyle \sum_{\ell \in {\Ncal}_{k}} a_{\ell k}\; \bpsi_{\ell,i} \quad \mbox{(diffusion with random updates)}
\end{aligned} \right.
\ee
whereas the \emph{asynchronous} ATC diffusion strategy \eqref{eqn:diffusionATCasync} reduces to:
\be
\label{eqn:diffusionATCLMSasync}
\left\{ \begin{aligned}
\bpsi_{k,i} & = \displaystyle \w_{k,i-1} + 2\bmu_k(i) \u_{k,i}\tran \left[\d_{k}(i)-\u_{k,i}\w_{k,i-1}\right] \\
\w_{k,i} & = \displaystyle \sum_{\ell \in {\small \bm{\Ncal}_{k,i}}} \a_{\ell k}(i)\; \bpsi_{\ell,i} \quad \mbox{(asynchronous diffusion)}
\end{aligned} \right.
\ee
\es
We can view implementation \eqref{eqn:diffusionATCLMSrandomupdate} as a special case of the asynchronous update \eqref{eqn:diffusionATCLMSasync} when the variances of the random combination coefficients $\{ \a_{\ell k}(i) \}$ are set to zero.

\qd}
\end{example}
}

{\small
\begin{example}[{\bf Asynchronous consensus LMS networks}]{\rm
Similarly, for the same MSE network of Example \ref{example:MSEnetwork}, the asynchronous consensus strategy is given by
\be
\label{eqn:consensusLMS.33}
\left\{ \begin{aligned}
\bpsi_{k,i-1} & \weq \displaystyle \sum_{\ell \in {\small \bm{\Ncal}_{k,i}}} \a_{\ell k}(i)\; \w_{\ell,i-1} \\
\w_{k,i} & \weq \displaystyle \bpsi_{k,i-1} + 2\bm{\mu}_k(i)\u_{k,i}\tran [\d_k(i)-\u_{k,i}\w_{k,i-1}]
\end{aligned} \right.
\ee
\qd}\label{example-consensus.33}
\end{example}
}

\subsection{Mean Graph}
\begin{figure}[h]
\epsfxsize 11cm \epsfclipon
\begin{center}
\leavevmode \epsffile{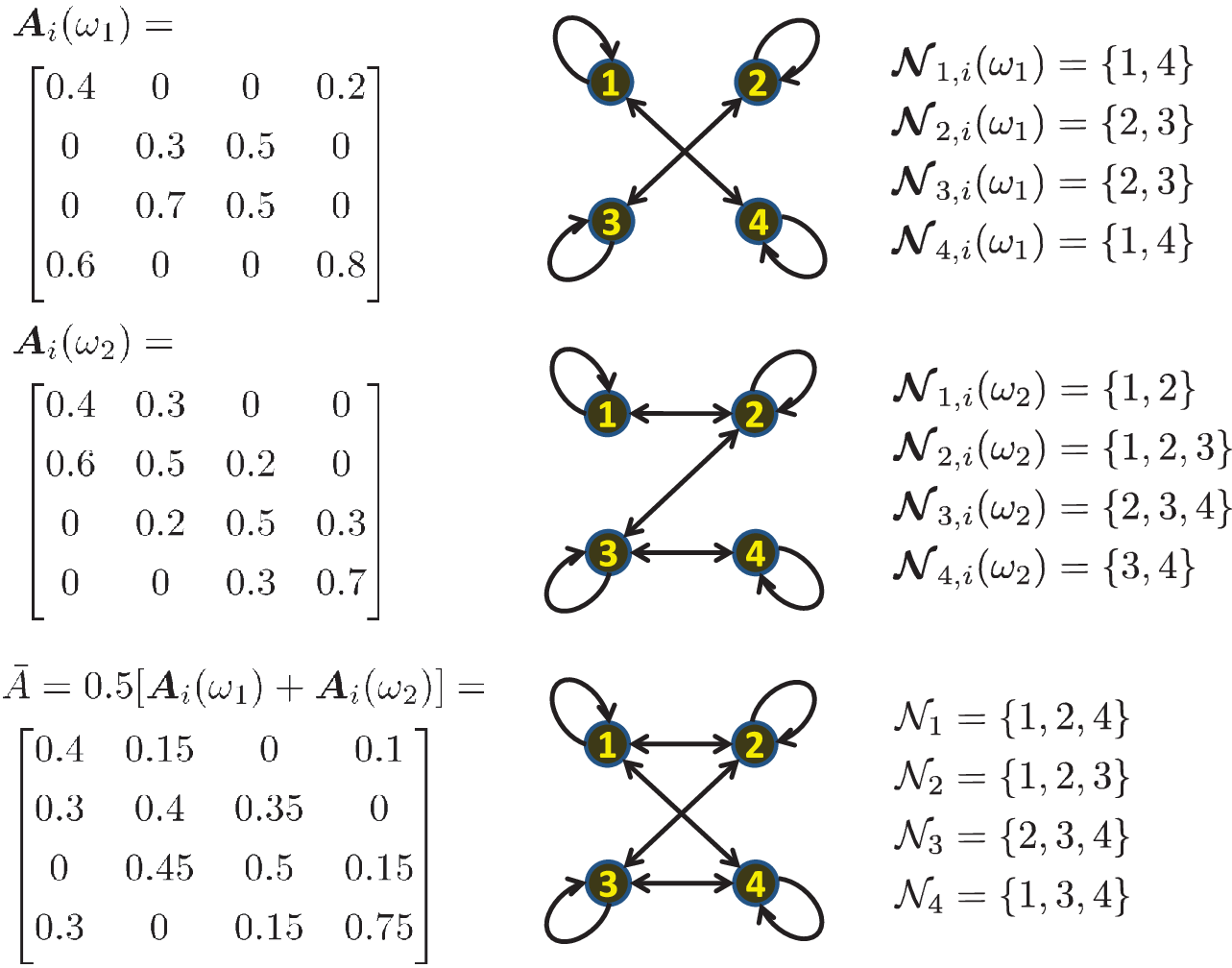}
\caption{{\small The first two rows show two equally probable realizations with the respective neighborhoods. The last row shows the resulting mean graph. This is a modified version of Fig. 6 from \cite{xiao2013}; extracted with permission.}}\label{fig:randomtopology}
\end{center}
\end{figure}

We refer to the topology that corresponds to the average combination coefficients $\{\bar{a}_{\ell k}\}$ as the {\em mean graph}, which is fixed over time. For each agent $k$, the neighborhood defined by the mean graph is denoted by $\Ncal_k$. It is straightforward to verify that the mean combination coefficients $\bar{a}_{\ell k}$ satisfy the following constraints over the mean graph (compare with \eqref{eqn:conditionalk} and \eqref{eqn:randomtopologyconstraints}):
\be
\addbox{\;\sum_{\ell \in \Ncal_k} \bar{a}_{\ell k} = 1, \quad \mbox{and} \quad
\begin{cases}
\bar{a}_{\ell k} > 0, & \mbox{if} \; \ell\in \Ncal_k \\
\bar{a}_{\ell k} = 0, & \mbox{otherwise}
\end{cases}\;}
\ee
 One example of a random network with two equally probable realizations and its mean graph is shown in Fig.~\ref{fig:randomtopology}  \cite{xiao2013}. The letter $\omega$ is used to index the sample space of the random matrix $\A_i$.

There is one useful result that relates the random neighborhoods $\{ \bm{\Ncal}_{k,i} \}$ from \eqref{eqn:diffusionATCasync} to the neighborhoods $\{ \Ncal_k \}$ from the mean graph. It is not difficult to verify that $\Ncal_k$ is equal to the \emph{union} of all possible realizations for the random neighborhoods $\bm{\Ncal}_{k,i}$; this property is already illustrated by the example of Fig.~\ref{fig:randomtopology}:
\be
\Ncal_k = \bigcup_{\omega\in\Omega} \bm{\Ncal}_{k,i}(\omega)\;
\ee
for any $k$, where $\Omega$ denotes the sample space for $\bm{\Ncal}_{k,i}$.

\subsection{Random Combination Policy}
The first and second-order moments of the combination coefficients will play an important role in characterizing the stability and mean-square-error performance of the asynchronous network (\ref{eqn:diffusionATCasync}). We collect these moments in matrix form as follows.  We first group the combination coefficients into a matrix:
\bs
\bq
\A_i &\define& [\,\a_{\ell k}(i)\,]_{\ell,k=1}^N\quad\;\;\;\;\;\;(N\times N)
\eq
\es
The sequence $\{\A_i,i\ge0\}$ represents a stochastic process consisting of left-stochastic random matrices whose entries satisfy the conditions in \eqref{eqn:randomtopologyconstraints} at every time $i$.  We subsequently introduce the mean and Kronecker-covariance matrix of $\A_i$ and  assume these quantities are constant over time; we denote them by the $N \times N$ matrix $\bar{A}$ and the $N^2 \times N^2$ matrix $C_A$, respectively:
\bs
\bq
\bar{A} & \define& \Ex \,\A_i = [\;\bar{a}_{\ell k}\;]_{\ell, k = 1}^N  \label{eqn:randcombinemean}
\\
C_A & \define &\Ex [(\A_i - \bar{A}) \kron (\A_i - \bar{A})] \label{eqn:randcombinecov}
%= [C_{a,\ell k}]_{\ell, k = 1}^N \\
%C_{a,\ell k} & \define \Ex [(\a_{\ell k}(i) - \bar{a}_{\ell k})(\A_i - \bar{A})] = [ c_{a,\ell k, n m} ]_{n, m =1}^{N}
\eq
\es
The matrix $C_A$ is \emph{not} a conventional covariance matrix and is not necessarily Hermitian. The reason for its introduction is because it captures the correlations of each entry of $\A_i$ with all other entries in $\A_i$. For example, for a network with $N=2$ agents, the entries of $\bar{A}$ and $C_A$ will be given by:
\bs
\bq
\bar{A}&=&\ba{cc}\bar{a}_{11}&\bar{a}_{12}\\\bar{a}_{21}&\bar{a}_{22}\ea\\\nn\\
C_A&=&\ba{ccccc}c_{a,11,11}&c_{a,11,12} &\vline& c_{a,12,11}& c_{a,12,12}\\
c_{a,11,21}&c_{a,11,22}&\vline&c_{a,12,21}& c_{a,12,22}\\\hline
c_{a,21,11}&c_{a,21,12}&\vline&c_{a,22,11}& c_{a,22,12}\\
c_{a,21,21}&c_{a,21,22}&\vline&c_{a,22,21}& c_{a,22,22}
\ea
\eq
\es
We thus see that the $(\ell,k)-$th block of $C_{A}$ contains the covariance coefficients of $\a_{\ell k}$ with all other entries of $\A_i$. One useful property of the matrices $\{\bar{A}, C_A\}$ so defined is that their elements are nonnegative and the following matrices are left-stochastic:
\be
\label{eqn:AACAone}
\addbox{\;\left(\bar{A}\right)\tran \one_N = \one_N,\;\;\; (\bar{A} \kron \bar{A} + C_A)\tran \one_{N^2} = \one_{N^2}\;}
\ee

{%\small
\begin{assumption}[{\bf Asynchronous network model}]
\label{assumption:asyncnetwork} {\em
It is assumed that the random processes $\{\bm{\mu}_k(i),\a_{\ell m}(j)\}$ are independent of each other for all $k, \ell, m, i,$ and $j$. They are also independent of any other random variable in the learning algorithm.}

\hfill \qd
\end{assumption}
}

The asynchronous network model described in this section covers many situations of practical interest. For example, we can choose the sample space for each step-size ${\bmu}_k(i)$ to be the binary choice $\{0, \mu\}$ to model random ``on-off'' behavior at each agent $k$ for the purpose of saving power, waiting for data, or even due to random agent failures. Similarly, we can choose the sample space for each combination coefficient ${\a}_{\ell\,k}(i)$, $\ell \in \Nknk$, to be $\{0, a_{\ell k}\}$ to model a random ``on-off'' status for the link from agent $\ell$ to agent $k$ at time $i$ for the purpose of either saving communication cost or due to random link failures.  Note that the convex constraint \eqref{eqn:randomtopologyconstraints} can always be satisfied by adjusting the value of $\a_{kk}(i)$ according to the realizations of $\{\a_{\ell k}(i); \ell \in \bm{\Ncal}_{k,i} \backslash \{k\}\}$.

{\small
\begin{example}[{\bf The spatially-uncorrelated model}]
\label{example:spatial}{\rm A useful special case of the asynchronous network model of this section is the spatially-uncorrelated model. In this case, at each iteration $i$, the random step-sizes $\{\bmu_k(i); k = 1,2,\dots, N\}$ are uncorrelated with each other across the network, and  the random combination coefficients $\{\a_{\ell k}(i); \ell \ne k, k=1,2,\dots,N\}$ are also uncorrelated with each other across the network. Then, it can be verified that the covariances $\{c_{\mu,k \ell}\}$ in \eqref{eqn:randstepsizecoventry} and $\{c_{a,\ell k, nm}\}$ in \eqref{eqn:randcombinecoventry} will be fully determined by the variances $\{\sigma_{\mu}^2,\,\sigma^2_{a,\ell k}\}$:
\bs
\be
\label{eqn:cov_mu_kl}
c_{\mu,kk}=\sigma_{\mu,k}^2,\;\;\;\;c_{\mu,k \ell} =0,\;k\neq \ell
\ee
and
\be\label{eqn:cov_a_lk_nm}
c_{a,\ell k, nm}  = \left\{\begin{array}{cl}
\;\;\;\sigma_{a,\ell k}^2, & \mbox{if} \; k = m, \ell = n, \ell \in \Nknk \\
-\sigma_{a,\ell k}^2, & \mbox{if} \; k = m = n, \ell \in \Nknk \\
-\sigma_{a,n k}^2, & \mbox{if} \; k = m = \ell, n \in \Nknk \\
\displaystyle \sum_{j \in\Nknk} \sigma_{a,j k}^2, & \mbox{if} \; k = m = \ell = n \\
0, & \mbox{otherwise} \end{array}\right.
\ee
\es

\qd}
\end{example}
}

\subsection{Perron Vectors}
Now that we introduced the network model, we can move on to examine the effect of network cooperation on performance. Some interesting patterns of behavior arise when agents cooperate to solve a global optimization problem in a distributed manner from streaming data \cite{ProcIEEE2014,NOW2014}. For example, one interesting result established in \cite{xiao2013,xiao2013AA,chen2013,Chen2013,ChenSayeAller} is that the effect of the network topology on performance is captured by the Perron vector of the combination policy. This vector turns out to summarize the influence of the topology on performance so much so that different topologies with similar Perron vectors will end up delivering similar performance. We explained the role of Perron vectors in the context of synchronous adaptation and learning in \cite{ProcIEEE2014,NOW2014}. Here we focus on asynchronous networks. In this case, two Perron vectors will be needed since the randomness in the combination policy is now represented by two moment matrices, $\bar{A}$ and $C_A$. In the synchronous case, only one Perron vector was necessary since the combination policy was fixed and described by a matrix $A$. Although we are focusing on the  asynchronous case in the sequel, we will be able to recover results for synchronous networks as special cases.

Let us first recall the definition of Perron vectors for synchronous networks, say, of the form described by \eqref{eqn:diffusionATC} with combination policy $A$. We assume the network is strongly-connected. In this case, the left-stochastic matrix $A$  will be primitive. For such primitive matrices, it follows from the Perron-Frobenius Theorem \cite{horn} that: (a) the matrix $A$ will have a {\em single} eigenvalue at one; (b)
 all other eigenvalues of $A$ will be strictly inside the unit circle so that $\rho(A)=1$; and
 (c) with proper sign scaling, all entries of the right-eigenvector of $A$ corresponding to the single eigenvalue at one will be {\em positive}. Let $p$ denote this right-eigenvector with its entries $\{p_k\}$ normalized  to add up to one, i.e.,
\be
\label{eqn:pdefsync}
\addbox{\;A p \weq p, \quad \one\tran p \weq 1, \quad p_{k} > 0,\;\; k = 1, 2,\dots,N\;}
\ee
We refer to $p$ as the \emph{Perron} eigenvector of $A$. It was explained in \cite{ProcIEEE2014,NOW2014} how the entries of this vector determine the mean-square-error performance and convergence rate of the network; these results will be revisited further ahead when we recover them as special cases of the asynchronous results.

On the other hand, for an asynchronous implementation, the individual realizations of the random combination matrix $\A_i$ in \eqref{eqn:diffusionATCasync} need not be primitive. In this context, we will require a form of primitiveness to hold on average as follows.

{
\begin{definition}[{\bf Strongly-connected asynchronous model}]
\label{assumption:primitive} {\em We say that an asynchronous model with random combination coefficients $\{\a_{\ell k}(i)\}$ is strongly-connected if the Kronecker-covariance matrix given by $\bar{A} \kron \bar{A} + C_A$ is primitive.}

 \hfill \qd
\end{definition}
}

\noindent Observe that if we set $\bar{A}=A$ and $C_A=0$, then we recover the condition for strong-connectedness in the synchronous case, namely, that $A$ should be primitive.

\begin{figure}[h]
\epsfxsize 11cm \epsfclipon
\begin{center}
\leavevmode \epsffile{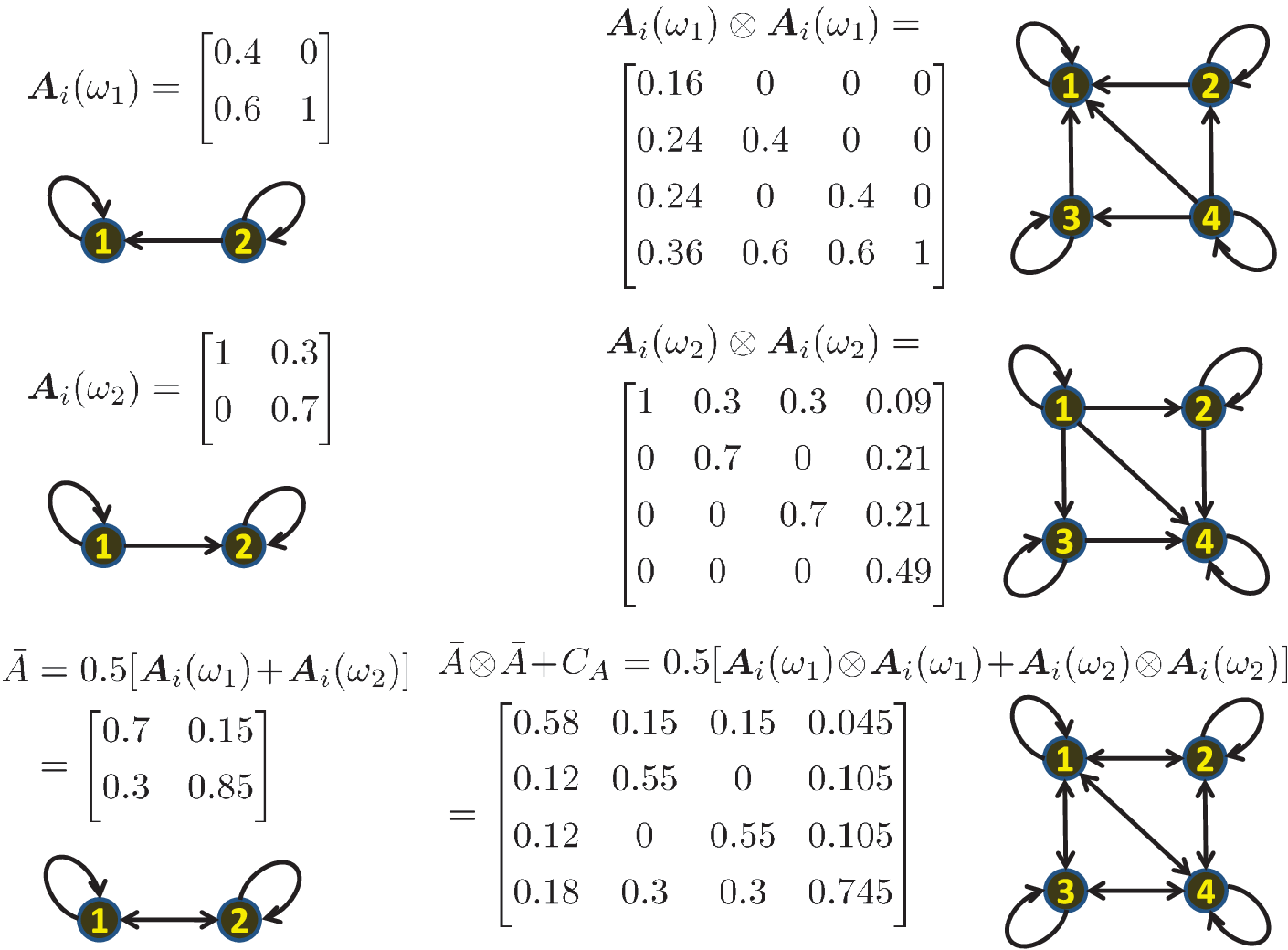}
\caption{{\small The combination policy $\A_i$ has two equally probable realizations in this example, denoted by $\{\bm{A}_i(\omega_1), \bm{A}_i(\omega_2)\}$. Observe that neither of the digraphs
$\bm{A}_i(\omega_1) \otimes \bm{A}_i(\omega_1)$ or $\bm{A}_i(\omega_2) \otimes \bm{A}_i(\omega_2)$
is strongly-connected due to the existence of the source and sink nodes. However, the digraph associated with $\Ex(\bm{A}_i \kron \bm{A}_i)$, which is the union of the first two digraphs, is strongly-connected, where information can flow in any direction through the network. This is a modified version of Fig. 7 from \cite{xiao2013AA}; extracted with permission.}}\label{fig:primitive}
\end{center}
\end{figure}

Definition \ref{assumption:primitive} means that the directed graph (digraph) associated with the matrix $\bar{A}  \kron  \bar{A}  +  C_A$ is strongly-connected (e.g., \cite[pp.~30,34]{BermanPF} and \cite{NOW2014}). It is straightforward to check from the definition of $C_A$ in \eqref{eqn:randcombinecov} that
\be
\bar{A} \kron \bar{A} + C_A\;=\;\Ex(\A_i\kron\A_i)
\ee
so that the digraph associated with $\bar{A}  \kron  \bar{A}  +  C_A$ is the union of all possible digraphs associated with the realizations of $\bm{A}_i  \kron  \bm{A}_i$ \cite[p.~29]{Bondy08}. Therefore, as explained in \cite{xiao2013,xiao2013AA,xiao2013BB}, definition \ref{assumption:primitive} amounts to an assumption that the \emph{union} of all possible digraphs associated with the realizations of $\bm{A}_i  \kron  \bm{A}_i$ is strongly-connected. As illustrated in Fig.~\ref{fig:primitive}, this condition still allows individual digraphs associated with realizations of $\bm{A}_i$ to be weakly-connected with or without self-loops or even to be disconnected \cite{xiao2013AA}.

It follows from property \eqref{eqn:AACAone} and Definition \ref{assumption:primitive}, the matrix $\bar{A} \kron \bar{A} + C_A$ is left-stochastic and primitive. It can be verified that mean matrix $\bar{A}$ is also  primitive if $\bar{A} \kron \bar{A} + C_A$ is primitive (although the converse is not true). Therefore, the matrix $\bar{A}$ is both left-stochastic and primitive. We denote the Perron eigenvector of $\bar{A} \kron \bar{A} + C_A$ by $p_c \in \real^{N^2 \times 1}$, which satisfies:
\bs
\be
\label{eqn:pdef}
\addbox{\;(\bar{A} \kron \bar{A} + C_A) p_c = p_c, \qquad  p_c\tran \one = 1\;}
\ee
Likewise, we denote the Perron eigenvector of $\bar{A}$ by $\bar{p} \in \real^{N \times 1}$, which satisfies:
\be
\label{eqn:barpdef}
\addbox{\;\bar{A} \bar{p} = \bar{p}, \qquad  \bar{p}\tran \one = 1\;}
\ee
\es
Observe that if we set $C_A=0$ and $\bar{A}=A$, then we recover the synchronous case information, namely, $\bar{p}=p$ and $p_c=p\otimes p$.

The entries of the Perron vectors $\{p_c,\bar{p}\}$ are related to each other, as illustrated by the following explanation (proofs appear in \cite{xiao2013AA}[App. VII]). Since the vector $p_c$ is of dimension $N^2\times1$, we partition it into $N$ sub-vectors of dimension $N \times 1$ each:
\bs
\be
p_c \define  \col\{ p_1, p_2, \dots, p_N \}
\ee
where $p_k$ denotes the $k$-th sub-vector. We  construct an $N \times N$ matrix $P_c$ from these sub-vectors:
\be
\label{eqn:Ppdef}
P_c \define [ p_1 \;\; p_2 \;\; \dots \;\; p_N ]\;=\;\left[p_{c,\ell k}\right]_{k,\ell=1}^N
\ee
We use $p_{\ell k}$ to denote the $(\ell,k)$-th element of matrix $P_c$, which is equal to the $\ell$-th element of $p_k$. It can be verified that the matrix $P_c$ in \eqref{eqn:Ppdef} is symmetric positive semi-definite and it satisfies \be \addbox{\;P_c \one = \bar{p},\;\;\;\;P_c\tran=P_c,\;\;\;\;P_c\geq 0\;}\ee where $\bar{p}$ is the Perron eigenvector in \eqref{eqn:barpdef}. We can further establish the following useful relations:
\be
\label{eqn:plkandpklandbarpk}
p_{c,\ell k} = p_{c,k \ell}, \quad \sum_{k=1}^N p_{c,\ell k} = \bar{p}_\ell, \quad \sum_{\ell=1}^N p_{c,\ell k} = \bar{p}_k
\ee
We can also verify that the matrix difference
\be
\label{eqn:Cpdef}
C_c \define P_c - \bar{p} \bar{p}\tran,\;\;\;C_c\tran=C_c,\;\;\;C_c\geq 0
\ee
\es
is symmetric, positive semi-definite, and satisfies $C_c \one = 0 $. Moreover, it is straightforward to verify that \be \addbox{\;c_{c,kk}\;=\;p_{c,k k} - \bar{p}_k^2\; \geq 0\;}\ee where $c_{c,kk}$ denotes the $(k,k)$-th entry in $C_c$.

\section{Asynchronous Network Performance}
We now comment on the performance of asynchronous networks and compare their metrics against both non-cooperative and centralized strategies.

\subsection{MSD Performance}
We denote the MSD performance of the individual agents and the average MSD performance across the network by: \bs\bq
\mbox{\rm MSD}_{{\rm dist}, k}&\define &\lim_{i\rightarrow\infty}\,\Ex\|\widetilde{\w}_{k,i}\|^2\\
\mbox{\rm MSD}_{\rm dist, av}&\define &\frac{1}{N}\sum_{k=1}^N\mbox{\rm MSD}_{{\rm dist},k}
\eq\es
where the error vectors are measured relative to the global optimizer, $w^{o}$. We further denote the gradient noise process at the individual agents by
\bs\be
\s_{k,i}(w)\define \widehat{\nabla_{w\tran} J}_{k}(w) \;-\;\nabla_{w\tran} J_{k}(w)\label{jakd7813.adlka}
\ee
and define
  \bq
H_k&\define&\nabla_w^2\;J_k(w^{o})\\
R_{s,k}&\define&\lim_{i\rightarrow\infty} \Ex\left[\s_{k,i}(w^{o})\s_{k,i}\tran (w^{o})\,|\,\bm{\cal F}_{i-1}\right]\label{2rpaldka}
\eq
where $\bm{\cal F}_{i-1}$ now represents the collection of all random events generated by the iterates from across all agents, $\{\w_{k,j},\,k=1,2,\ldots,N\}$, up to time $i-1$:
   \be \bm{\cal F}_{i-1}\define\mbox{\rm filtration}\{\bm{\w}_{k,-1}, \bm{\w}_{k,0},\bm{\w}_{k,1},\ldots,
\bm{\w}_{k,i-1},\;\mbox{\rm all}\; k\}\label{defined.adas.abc}\ee
\es
It is stated later in (\ref{eqn:MSDasyncATC}), under some assumptions on the gradient noise processes (\ref{jakd7813.adlka}), that for strongly-connected asynchronous diffusion networks of the form \eqref{eqn:diffusionATCrandomupdate}, and for sufficiently small step-sizes $\mu_x$ \cite{xiao2013AA,NOW2014}:
\be
\label{eqn:MSDasyncATC.1}
\MSD_{\dist,k}^{\async} \approx \MSD_{\dist,\av}^{\async} \approx\frac{1}{2} \Tr \left[ \left( \sum_{k=1}^N \bar{\mu}_k \bar{p}_k H_k \right)^{-1} \left( \sum_{k=1}^N (\bar{\mu}_k^2 + \sigma^2_{\mu,k}) p_{c,kk}  R_{s,k} \right) \right]
\ee

\noindent The same result holds for asynchronous consensus and CTA diffusion strategies.  Observe from (\ref{eqn:MSDasyncATC.1}) the interesting conclusion that the distributed strategy is able to {\em equalize} the MSD performance across all agents for sufficiently small step-sizes.  It is also instructive to compare expression \eqref{eqn:MSDasyncATC} with \eqref{eqn:MSDbatchasyncgeneral} and (\ref{eqn:MSDbatchasyncgeneraasxasxl}) in the centralized case. Observe how cooperation among agents leads to the appearance of the scaling coefficients $\{\bar{p}_k,p_{c,kk}\}$; these factors are determined by $\bar{A}$ and $C_A$.

Note further that if we set $\bar{\mu}_k=\mu$, $\sigma_{\mu,k}^2=0$, $\bar{p}_k=p_k$ and $p_{c,kk}=p_k^2$, then we recover the MSD expression for synchronous distributed strategies:
\be
\label{eqn:MSDasyncATC.1a}
\MSD_{\dist,k}^{\sync} \approx \MSD_{\dist,\av}^{\sync} \approx
 \frac{1}{2} \Tr \left[ \left( \sum_{k=1}^N \mu_k p_k H_k \right)^{-1} \left( \sum_{k=1}^N \mu_k^2 p_{k}^2 R_{s,k} \right) \right]
\ee
This result agrees with expression (62) from \cite{ProcIEEE2014}.

{\small
\begin{example}[{\bf MSE networks with random updates}]{\rm
We continue with the setting of Example \ref{example:MSEnetwork}, which deals with MSE networks. We assume the first and second-order moments of the random step-sizes are uniform, i.e., $\bar{\mu}_k \equiv \bar{\mu}$ and $c_{\mu,kk} \equiv \sigma_{\mu}^2$, and also assume uniform regression covariance matrices, i.e., $R_{u,k} \equiv R_u$ for $k = 1,2,\dots,N$. It follows that $H_k = 2 R_u \equiv H$ and $R_{s,k} = 4 \sigma_{v,k}^2 R_u$. Substituting into \eqref{eqn:MSDasyncATC}, and assuming a fixed topology with fixed combination coefficients set to $a_{\ell k}$, we conclude that the MSD performance of the diffusion strategy \eqref{eqn:diffusionATCLMSrandomupdate} with random updates  is well approximated by:
\bs
\be
\label{eqn:MSDATCLMSrandomupdate}
\MSD_{\dist,k}^{\async, 1} \; \approx \; \MSD_{\dist,\av}^{\async, 1} \;\approx\; \mu_x M \, \left(\sum_{k=1}^N {p}_{k}^2 \sigma_{v,k}^2\right)
\ee
where
\be \mu_x\;=\;\bar{\mu}\;+\;\frac{\sigma_{\mu}^2}{\bar{\mu}}\ee and ${p}_k$ is the $k-$th entry of the Perron vector $p$ defined by (\ref{eqn:pdefsync}).

If the combination matrix $A$ happens to be doubly stochastic, then its Perron eigenvector becomes $p = \one/N$. Substituting $p_k = 1/N$ into \eqref{eqn:MSDATCLMSrandomupdate} gives
\be
\label{eqn:MSDATCLMSrandomupdateuniform}
\MSD_{\dist, k}^{\async, 1} \;\approx\; \MSD_{\dist,\av}^{\async, 1} \;\approx\; \frac{\mu_x M}{N}\, \left( \frac{1}{N}\, \sum_{k=1}^N \sigma_{v,k}^2\right)
\ee
which agrees with the centralized performance \eqref{eqn:MSDbatchrandomupdateLMS}. In other words, the asynchronous diffusion strategy is able to match the performance of the centralized solution for doubly stochastic combination policies, when both implementations employ random updates. Since the centralized solution can improve the average MSD performance over non-cooperative networks, we further conclude that the diffusion strategy can also exceed the average performance of non-cooperative networks.
\es

\qd}
\label{kajd7813.13lkex}\end{example}
}

{\small
\begin{example}[{\bf Asynchronous MSE networks}]{\rm
We continue with the setting of Example~\ref{kajd7813.13lkex} except that we now employ the asynchronous LMS diffusion network \eqref{eqn:diffusionATCLMSasync}. Its MSD performance  is well approximated by:
\bs
\be
\label{eqn:MSDATCLMSasync}
\MSD_{\dist,k}^{\async} \; \approx \; \MSD_{\dist,\av}^{\async} \;\approx\; \mu_x M \, \left(\sum_{k=1}^N p_{c,kk} \sigma_{v,k}^2\right)
\ee
If the mean combination matrix $\bar{A}$  happens to be doubly stochastic, then its Perron eigenvector becomes $\bar{p} = \one/N$. Substituting $\bar{p}_k = 1/N$ into \eqref{eqn:MSDATCLMSrandomupdate}, and using $p_{c,k k} = \bar{p}_k^2 + c_{c,kk}$, where $c_{c,kk}$ is from \eqref{eqn:Cpdef}, gives
\be
\label{eqn:MSDATCLMSasyncuniform}
\MSD_{\dist, k}^{\async} \approx \MSD_{\dist,\av}^{\async} \approx \frac{\mu_x M}{N}\, \left[ \frac{1}{N}\, \sum_{k=1}^N (1 + N^2 c_{c,kk})\sigma_{v,k}^2\right]
\ee

\smallskip

\noindent It is clear that if $c_{c,kk} = \sigma^2_{\pi,k}$, then the MSD performance in \eqref{eqn:MSDATCLMSasyncuniform} will agree with the centralized performance \eqref{eqn:MSDbatchasyncLMS}. In other words, the distributed  diffusion strategy is able to match the performance of the centralized solution.
\es

\qd}
\end{example}
}

\section{Network Stability and Performance}
\label{sec.dxxx}
In this section, we examine more closely the performance and stability results that were alluded to in the earlier sections. We first examine the consensus and diffusion strategies in a unified manner, and subsequently focus on diffusion strategies due to their enhanced stability properties, as the ensuing discussion will reveal.

\subsection{MSE Networks}
We motivate the discussion by presenting first some illustrative examples with MSE networks, which involve quadratic costs. Following the examples, we extend the framework to more general costs.

{\small
\begin{example}[{\bf Error dynamics over MSE networks}] {\rm We consider the MSE network of Example~\ref{example-9A}, which involves quadratic costs with a common minimizer, $w^o$. The update equations for the non-cooperative, consensus, and diffusion strategies are given
by (\ref{eqn:SGDnoncoopk}), (\ref{eqn:consensusLMS}), and (\ref{eqn:diffusionCTALMS})--(\ref{eqn:diffusionATCLMS}).
We can group these strategies into a single unifying description by considering the
following structure in terms of three sets of combination coefficients
$\{\a_{o,\ell k}(i), \a_{1, \ell k}(i), \a_{2,\ell k}(i)\}$:
\be
\left\{\begin{array}{l}\bm{\phi}_{k,i-1}=\displaystyle \sum_{\ell\in\bm{\cal N}_{k,i}} \a_{1,\ell k}(i)\w_{\ell,i-1}\\
\bm{\psi}_{k,i} =   \displaystyle \sum_{\ell\in\bm{\cal N}_{k,i}}\hspace{-0.3cm} \a_{o,\ell k}(i)\bm{\phi}_{\ell,i-1}
+ 2\bm{\mu}_k(i) \u_{k,i}\tran \left[\d_{k}(i)-\u_{k,i}\bm{\phi}_{k,i-1}\right]\\
\w_{k,i}=\displaystyle \sum_{\ell\in\bm{\cal N}_{k,i}} \a_{2,\ell k}(i)\bm{\psi}_{\ell,i}
\end{array}\right.
\label{dynahga,da.abc}\ee

\noindent In (\ref{dynahga,da.abc}), the quantities $\{\bm{\phi}_{k,i-1},\bm{\psi}_{k,i}\}$ denote $M\times 1$
intermediate variables, while the nonnegative
entries of the $N\times N$ matrices $\A_{o,i}=[\a_{o,\ell k}(i)]$, $\A_{1,i}=[\a_{1,\ell k}(i)]$, and $\A_{2,i}=[\a_{2,\ell k}(i)]$ are assumed to
satisfy the same conditions (\ref{eqn:randomtopologyconstraints}). Any of the
combination weights $\{\a_{o,\ell k}(i), \a_{1,\ell k}(i), \a_{2, \ell k}(i)\}$ is zero
whenever $\ell\notin \bm{\cal N}_{k,i}$. Different choices for $\{\A_{o,i},\A_{1,i},\A_{2,i}\}$, including random and deterministic choices, correspond to different strategies, as the following examples reveal:
\bs\bq
\hspace{-0.4cm}\mbox{\rm non-cooperative:}\hspace{-0.2cm} && \hspace{-0.3cm} \A_{1,i}=\A_{o,i}=\A_{2,i}=I_N\label{cons.alkd131.axasva}\\
\hspace{-0.4cm}\mbox{\rm consensus:}\hspace{-0.2cm}  &&\hspace{-0.3cm} \A_{o,i}=\A_i,\;\A_{1,i}=I_N=\A_{2,i}\label{cons.alkd131.a}\\
\hspace{-0.4cm}\mbox{\rm CTA diffusion:}\hspace{-0.2cm}  &&\hspace{-0.3cm} \A_{1,i}=\A_i,\;\A_{2,i}=I_N=\A_{o,i}\\
\hspace{-0.4cm}\mbox{\rm ATC diffusion:}\hspace{-0.2cm}  &&\hspace{-0.3cm} \A_{2,i}=\A_i,\;\A_{1,i}=I_N=\A_{o,i}\label{cons.alkd131.c.a}
\eq\es
\noindent where $\A_i$ denotes some generic combination policy satisfying (\ref{eqn:randomtopologyconstraints}). We associate with each agent $k$ the following three errors:
\bs\bq
\widetilde{\w}_{k,i}&\define& w^o-\w_{k,i}\\
\widetilde{\bm{\psi}}_{k,i}&\define& w^o-\bm{\psi}_{k,i}\\
\widetilde{\bm{\phi}}_{k,i-1}&\define& w^o-\bm{\phi}_{k,i-1}
\eq
\es
\noindent which measure the deviations from the global minimizer, $w^o$. Subtracting $w^o$
from both sides of the equations in (\ref{dynahga,da.abc}) we get
\bs\be
\left\{\begin{array}{l}
\widetilde{\bm{\phi}}_{k,i-1}=\displaystyle
\sum_{\ell\in\bm{\cal N}_{k,i}} \a_{1,\ell k}(i)\;\widetilde{\w}_{\ell,i-1}\\
\widetilde{\bm{\psi}}_{k,i}=\displaystyle
\sum_{\ell\in\bm{\cal N}_{k,i}} \a_{o,\ell k}(i)\widetilde{\bm{\phi}}_{\ell,i-1}
- 2\bm{\mu}_k(i)\u_{k,i}\tran \u_{k,i}\widetilde{\bm{\phi}}_{k,i-1}\;-\; 2\bm{\mu}_k(i)\u_{k,i}\tran \v_k(i)\\
\widetilde{\w}_{k,i}=\displaystyle
\sum_{\ell\in\bm{\cal N}_{k,i}} \a_{2,\ell k}(i)\;\widetilde{\bm{\psi}}_{\ell,i}
                    \end{array}\right.
\label{dynahga,da.4.abc}\ee
\noindent In a manner similar to (\ref{eqn:LMSgradientnoisedef}), the gradient noise process at each agent $k$ is given by
\be
\s_{k,i}(\bm{\phi}_{k,i-1})\;=\;2\left(R_{u,k}-\u_{k,i}\tran\u_{k,i}\right)
\widetilde{\bm{\phi}}_{k,i-1}\;-\;2\u_{k,i}\tran\v_k(i)\label{gradinajl13.1}
\ee
\es
\noindent In order to examine the evolution of the error dynamics across the network, we collect the
error vectors from all agents into $N\times 1$ block error vectors (whose
individual entries are of size $M\times 1$ each):
\bs\be
\hspace{-0.2cm}\widetilde{\w}_{i}\define \ba{c}\widetilde{\w}_{1,i}\\
\widetilde{\w}_{2,i}\\\vdots\\
\widetilde{\w}_{N,i}
\ea,\;
\widetilde{\bm{\psi}}_{i}\define \ba{c}\widetilde{\bm{\psi}}_{1,i}\\
\widetilde{\bm{\psi}}_{2,i}\\\vdots\\
\widetilde{\bm{\psi}}_{N,i}
\ea,\;
\widetilde{\bm{\phi}}_{i-1}\define \ba{c}\widetilde{\bm{\phi}}_{1,i-1}\\
\widetilde{\bm{\phi}}_{2,i-1}\\\vdots\\
\widetilde{\bm{\phi}}_{N,i-1}
\ea
\label{extend.2.abc}\ee
\noindent Motivated by the last term in the second equation in (\ref{dynahga,da.4.abc}), and by the gradient noise terms (\ref{gradinajl13.1}), we also introduce the
following $N\times 1$ column vectors whose entries are of size $M\times 1$ each:
\be
\z_i\define \ba{c}2\u_{1,i}\tran \v_1(i)\\
2\u_{2,i}\tran \v_2(i)\\\vdots\\
2\u_{N,i}\tran \v_N(i)\ea,\;\;\;\;
\s_i\define \ba{c}\s_{1,i}(\bm{\phi}_{1,i-1})\\\s_{2,i}(\bm{\phi}_{2,i-1})\\\vdots\\
\s_{N,i}(\bm{\phi}_{N,i-1})
\ea\label{defiadk13.adlk}\ee
\es
We further introduce the Kronecker products \bs\be \bm{\cal A}_{o,i}\define \A_{o,i}\otimes
 I_{M},\quad\bm{\cal A}_{1,i}\define
\A_{1,i}\otimes I_{M},\quad\bm{\cal A}_{2,i}\define \A_{2,i}\otimes I_{M}
\label{defioas.a.abc}\ee and the following $N\times N$ {\em block}
diagonal matrices, whose individual entries are of size $M\times M$
each:
\bq \hspace{-0.4cm}\bm{\cal M}_i&\hspace{-0.3cm}\define\hspace{-0.3cm}&\mbox{\rm diag}\{\;\bm{\mu}_1(i)
I_{M},\;\bm{\mu}_2(i) I_{M},\;\ldots,\;\bm{\mu}_N(i) I_{M}\;\}\label{m.ladlk.abc}\\
\hspace{-0.4cm} \bm{\cal R}_i&\hspace{-0.3cm}\define\hspace{-0.3cm}&\mbox{\rm diag}\left\{\;2\u_{1,i}\tran \u_{1,i},\;
 2\u_{2,i}\tran \u_{2,i},\;\ldots,\;2\u_{N,i}\tran\u_{N,i}\;\right\}\eq\es
\noindent From (\ref{dynahga,da.4.abc}) we can then easily conclude that the
block network variables satisfy the recursions:
\bs\be
\left\{\begin{array}{lcl}
\widetilde{\bm{\phi}}_{i-1}&=&\bm{\cal A}_{1,i}\tran\widetilde{\w}_{i-1}\\
\widetilde{\bm{\psi}}_{i}&=&\left(\bm{\cal A}_{o,i}\tran\;-\;\bm{\cal M}_i\bm{\cal R}_i\right)
\widetilde{\bm{\phi}}_{i-1}\;-\;\bm{\cal M}_i\z_i\\
\widetilde{\w}_i&=&\bm{\cal A}_{2,i}\tran\widetilde{\bm{\psi}}_i
\end{array}\right.
\ee
so that the network weight error vector, $\widetilde{\w}_i$, evolves according to:\be
\widetilde{\w}_i\;=\;\bm{\cal A}_{2,i}\tran\left(\bm{\cal A}_{o,i}\tran-\bm{\cal M}_i\bm
{\cal R}_i\right)\bm{\cal A}_{1,i}\tran\widetilde{\w}_{i-1}\;-\;\bm{\cal
A}_{2,i}\tran\bm{\cal M}_i\z_i
\label{like.lakd.ada.abc}\ee\es
\noindent For comparison purposes, if each agent operates individually and uses the
non-cooperative strategy (\ref{eqn:SGDnoncoopk}), then the
weight error vector would instead evolve according to
the following recursion: \be
\widetilde{\w}_i\;=\;\left(I_{MN}-\bm{\cal M}_i\bm{\cal
R}_{i}\right)\widetilde{\w}_{i-1}\;-\;\bm{\cal M}_i\z_i,\;\;i\geq
0
\label{like.lakd.no.ada.abc}\ee where the matrices $\{\bm{\cal A}_{o,i},\bm{\cal A}_{1,i},\bm{\cal A}_{2,i}\}$
 do not appear any longer, and with a block diagonal coefficient matrix $\left(I_{MN}-\bm{\cal M}_i\bm{\cal
R}_{i}\right)$. It is also straightforward to verify that recursion (\ref{like.lakd.ada.abc}) can be equivalently rewritten in the following form in terms of the gradient noise vector, $\s_i$, defined by (\ref{defiadk13.adlk}):
\bs\bq
\widetilde{\w}_i&=&\bm{\cal B}_i\,\widetilde{\w}_{i-1}\;+\;\bm{\cal
A}_{2,i}\tran{\cal M}\s_i
\label{like.lakd.ada.abc.BBB}\eq
where
\bq
\hspace{-0.4cm}\bm{\cal B}_i&\define& \bm{\cal A}_{2,i}\tran\left(\bm{\cal A}_{o,i}\tran-\bm{\cal M}_i
{\cal R}\right)\bm{\cal A}_{1,i}\tran
\\
\hspace{-0.4cm}{\cal R}&\define& \Ex\bm{\cal R}_i\;=\;\mbox{\rm diag}\{2R_{u,1},\;2R_{u,2},\;\ldots,\;2R_{u,N}\}
\label{cal r.lad}
\eq\es
}
\qd\label{example-19}
\end{example}
}

{\small
\noindent \begin{example}[{\bf Mean-error behavior}]{\rm We continue with the setting of Example~\ref{example-19}. In mean-square-error analysis, we are interested in examining how the quantities $\Ex\widetilde{\w}_{i}$ and
$\Ex\|\widetilde{\w}_i\|^2$ evolve over time. If we refer back to the data model described in Example~\ref{example-9A}, where the regression data $\{\u_{k,i}\}$ were assumed to be temporally white and independent over space, then  the stochastic matrix $\bm{\cal R}_i$ appearing in (\ref{like.lakd.ada.abc})--(\ref{like.lakd.no.ada.abc}) is seen to be statistically independent of  $\widetilde{\w}_{i-1}$. We further assume that, in the unified formulation, the entries of the combination policies $\{\A_{o,i},\A_{1,i},\A_{2,i}\}$ are independent of each other (as well as over time) and of any other variable in the learning algorithm. Therefore, taking expectations of both sides of these recursions, and invoking the fact that $\u_{k,i}$ and $\v_k(i)$ are also independent of each other and have zero means
 (so that $\Ex\z_i=0$), we conclude that the mean-error vectors evolve according to the following recursions:
\bs\bq
\hspace{-0.5cm}\Ex\widetilde{\w}_i\hspace{-0.2cm}&=&\hspace{-0.2cm}
\bar{\cal B}\, \left(\Ex\widetilde{\w}_{i-1}\right)\;\;\;\quad\quad\quad\quad\quad\quad(\mbox{\rm distributed})\label{like.lakd.ada.abc.xasyacc}\\
\hspace{-0.5cm}\Ex\widetilde{\w}_i\hspace{-0.2cm}&=&\hspace{-0.2cm}\left(I_{MN}-\bar{\cal M}{\cal
R}\right)\, \left(\Ex\widetilde{\w}_{i-1}\right)\;\;\;\;\;(\mbox{\rm non-cooperative})
\label{like.lakd.no.ada.abc.2}\eq
\es
where
\bs
\bq
\bar{\cal B}&\define& \Ex\bm{\cal B}_i\;=\;\bar{\cal A}_{2}\tran\left(\bar{\cal A}_{o}\tran-\bar{\cal M}
{\cal R}\right)\bar{\cal A}_{1}\tran\\
\bar{\cal M}&=& \Ex\bm{\cal M}_i=\mbox{\rm diag}\{\;\bar{\mu}_1
I_{M},\;\ldots,\;\bar{\mu}_N I_{M}\;\}\\
\bar{\cal A}_o&=&\Ex\bm{\cal A}_{o,i}\\
\bar{\cal A}_1&=&\Ex\bm{\cal A}_{1,i}\\
\bar{\cal A}_2&=&\Ex\bm{\cal A}_{2,i}
\eq
\es
The matrix $\bar{\cal B}$ controls the dynamics of the mean weight-error vector for the distributed strategies. Observe, in particular, that $\bar{\cal B}$ reduces to the following forms for the various strategies
(non-cooperative (\ref{eqn:SGDnoncoopk}), consensus (\ref{eqn:consensusLMS}), and diffusion (\ref{eqn:diffusionCTALMS})--(\ref{eqn:diffusionATCLMS})):
\bs\bq
\bar{\cal B}_{\rm ncop}&=&I_{MN}-\bar{\cal M}
{\cal R}\label{exp.1}\\
\bar{\cal B}_{\rm cons}&=&\bar{\cal A}\tran-\bar{\cal M}
{\cal R}\label{conshdj.13lk1}\\
\bar{\cal B}_{\rm atc}&=&\bar{\cal A}\tran\left(I_{MN}-\bar{\cal M}
{\cal R}\right)\label{use.ajd61312}\\
\bar{\cal B}_{\rm cta}&=&\left(I_{MN}-\bar{\cal M}
{\cal R}\right)\bar{\cal A}\tran\label{exp.4}
\eq\es
where $\bar{\cal A}=\bar{A}\otimes I_M$ and $\bar{A}=\Ex\A_i$.
}
\qd
\label{example-20}
\end{example}
}

{\small
\noindent \begin{example}[{\bf MSE networks with uniform agents}]{\rm The results of Example~\ref{example-20}  simplify when  all agents employ step-sizes with the same mean value, $\bar{\mu}_k\equiv\bar{\mu}$, and observe regression data with the
same covariance matrix, $R_{u,k}\equiv R_u$ \cite{bookchapter,yusayed12bb}. In this case, we can  express $\bar{\cal M}$ and
${\cal R}$ from (\ref{m.ladlk.abc}) and (\ref{cal r.lad}) in Kronecker product form as follows:
\be
\bar{\cal M}=\bar{\mu} I_N\otimes I_M,\;\;\;\;\;\;{\cal R}=I_N\otimes 2R_u
\ee
so that expressions (\ref{exp.1})--(\ref{exp.4})  reduce to
\bs\bq
\bar{\cal B}_{\rm ncop}&=&I_N\otimes (I_M\;-\;2\bar{\mu} R_u)\\
\bar{\cal B}_{\rm cons}&=&\bar{A}\tran\otimes I_M\;-\;2\bar{\mu} (I_M\otimes R_u)\label{simaldh71,ad.2}\\
\bar{\cal B}_{\rm atc}&=&\bar{A}\tran\otimes (I_M-2\bar{\mu} R_u)\label{simaldh71,ad.1}\\
\bar{\cal B}_{\rm cta}&=&\bar{A}\tran\otimes (I_M-2\bar{\mu} R_u)
\label{simaldh71,ad}\eq\es
Observe  that $\bar{\cal B}_{\rm atc}=\bar{\cal B}_{\rm cta}$, so we denote these matrices
by  $\bar{\cal B}_{\rm diff}$. Using properties of the
 eigenvalues of Kronecker products of matrices, it can be easily verified that the $MN$ eigenvalues of the above $\bar{\cal B}$ matrices  are given by the following expressions in terms of the eigenvalues of the component matrices $\{\bar{A},R_u\}$ for $k=1,2,\ldots N$ and $m=1,2,\ldots,M$:
\bs\bq
\lambda(\bar{\cal B}_{\rm ncop})&=&1-2\bar{\mu}\lambda_m(R_u)\label{conajda.1}\\
\lambda(\bar{\cal B}_{\rm cons})&=&\lambda_k(\bar{A})-2\bar{\mu}\lambda_m(R_u)\label{conajda}\\
\lambda(\bar{\cal B}_{\rm diff})&=&\lambda_k(\bar{A})\, \left[1-2\bar{\mu}\lambda_m(R_u)\right]\label{j7condha.da}
\eq
\es
}
\qd
\end{example}
}

{\small
\noindent \begin{example}[{\bf Potential instability in consensus networks}]{\rm Consensus strategies can become unstable when used for adaptation purposes \cite{yusayed12bb,NOW2014}. This undesirable effect is already reflected in expressions (\ref{conajda.1})--(\ref{j7condha.da}). In particular, observe that the eigenvalues of $\bar{A}$ appear multiplying $(1-2\mu\lambda_m(R_u))$ in expression (\ref{j7condha.da}) for diffusion. As such, and since $\rho(\bar{A})=1$ for any left-stochastic matrix, we conclude for this case of uniform agents that  $\rho(\bar{\cal B}_{\rm diff})=\rho(\bar{\cal B}_{\rm ncop})$. It follows that, regardless of the choice of the mean combination policy $\bar{A}$, the diffusion strategies will be stable in the mean (i.e., $\Ex\widetilde{\w}_i$ will converge asymptotically to zero) whenever the individual non-cooperative agents are stable in the mean:
\bs \be
\mbox{\rm individual agents stable}\;\Longrightarrow\;\mbox{\rm diffusion networks stable}
\ee
\noindent The same conclusion is not true for consensus networks; the individual agents can be stable and yet the consensus network can become unstable.  This is because $\lambda_k(\bar{A})$ appears as an additive (rather than multiplicative) term in (\ref{conajda}) (see \cite{sayedSPM,yusayed12bb,NOW2014} for examples):
\be
\mbox{\rm individual agents stable}\;\nRightarrow\;\mbox{\rm consensus networks stable}
\ee\es
The fact that the combination matrix $\bar{\cal A}\tran$ appears in an additive form in (\ref{conshdj.13lk1}) is the result of the asymmetry that was mentioned earlier in the update equation for the consensus strategy. In contrast, the update equations for the diffusion strategies lead to $\bar{\cal A}\tran$ appearing in a multiplicative form in (\ref{use.ajd61312})--(\ref{exp.4}).
}
\qd
\label{example-consensus.2}
\end{example}
}

{\small
\noindent \begin{example}[{\bf Useful stability result}]{\rm It is observed from expressions (\ref{use.ajd61312})--(\ref{exp.4}) in the asynchronous case, as well as from the corresponding expressions (81c)--(81d) in the synchronous case studied in \cite{ProcIEEE2014}, that the mean stability of diffusion strategies usually involves examining the stability of a matrix product of the form:
\be
{\cal B}\define {\cal A}_2\tran {\cal D}{\cal A}_1\tran
\label{form.11as}\ee
where ${\cal D}$ is a block diagonal symmetric matrix with blocks of size $M\times M$, while ${\cal A}_1$ and ${\cal A}_2$ are Kronecker product matrices defined in terms of $N\times N$ left-stochastic matrices $A_1$ and $A_2$ as ${\cal A}_1=A_1\otimes I_M$ and ${\cal A}_2=A_2\otimes I_M$. For example, in (\ref{use.ajd61312}) we have $A_1=I_N$, $A_2=\bar{A}$, ${\cal A}_1=I_{MN}$, ${\cal A}_2=\bar{\cal A}$,  and ${\cal D}=I_{MN}-{\cal M}{\cal R}$.

Matrix products of the form (\ref{form.11as}) are induced by the cooperation mechanism that is inherent to diffusion learning. They have a useful property: it turns out that these matrix products are stable, regardless of $A_1$ and $A_2$, as long as ${\cal D}$ is stable (i.e., has all its eigenvalues strictly inside the unit disc).  This useful result is easy to establish for {\em symmetric} left-stochastic matrices $A_1$ and $A_2$, as already noted in \cite{Lopes08}. This is because for symmetric matrices, their spectral radii coincide with their $2-$induced norms and, hence,
\be
\rho(A_1)=\|A_1\|,\;\;\;\rho(A_2)=\|A_2\|
\label{lkadj.l1l3kl1k3as}\ee
Consequently, since we already know that $\rho(A_1)=\rho(A_2)=1$, it follows that
\bq
\rho({\cal B})&\leq&\|{\cal B}\|\nn\\
&\leq & \|{\cal A}_2\|\cdot \|{\cal D}\|\cdot \|{\cal A}_1\|\nn\\
&=&\rho({\cal A}_2)\cdot \rho({\cal D})\cdot \rho({\cal A}_1)\nn\\
&=&\rho(A_2)\cdot \rho({\cal D})\cdot \rho(A_1)\nn\\
&=&\rho({\cal D})\label{llkad.l1k3klakdl}
\eq
which confirms that a stable ${\cal D}$ guarantees a stable ${\cal B}$, regardless of $A_1$ and $A_2$.

The conclusion that ${\cal B}$ in (\ref{form.11as}) is stable whenever ${\cal D}$ is stable continues to hold even when the matrices $A_1$ and $A_2$ are {\em not} necessarily symmetric. However, the argument leading to (\ref{llkad.l1k3klakdl}) will need to be adjusted because property  (\ref{lkadj.l1l3kl1k3as}) need not hold anymore. This more general result was established in \cite[App. D]{bookchapter} and also in \cite[App. A, pp. 3471--3473]{zhaoxiao},  where it was shown that multiplication of a symmetric block diagonal matrix ${\cal D}$ by any (not necessarily symmetric) left-stochastic Kronecker-product transformations from left and right generally reduces the spectral radius, i.e.,
\be
\rho\left({\cal A}_2\tran {\cal D}{\cal A}_1\tran\right)\;\leq\;\rho({\cal D})\label{llkad/;l13l;1l3}
\ee
Accordingly, a stable ${\cal D}$ again ensures a stable ${\cal B}$. This conclusion was established in the above references by replacing the $2-$induced norm used to arrive at (\ref{llkad.l1k3klakdl}) by a more convenient block-maximum norm, denoted by $\|\cdot\|_{b,\infty}$ and defined as follows.

Let $x = \mathrm{col}\{x_1,x_2,\ldots,x_N\}$ denote an $N\times 1$
{\em block} column vector whose individual entries are themselves vectors of
size $M\times 1$ each. Following \cite{Ber97,Takahashi10,bookchapter}, the block
maximum norm of $x$ is denoted by $\|x\|_{b,\infty}$ and is defined
as \be
\|x\|_{b,\infty} \define  \max_{1 \le k \le N}
 \|x_k\|
\label{va.lalkd.13}\ee That is, $\|x\|_{b,\infty}$ is equal to the largest Euclidean norm of its block components (this definition extends the regular notion of the $\infty-$norm of a vector to block vectors). The vector norm (\ref{va.lalkd.13}) induces a block-maximum matrix norm. Let ${\cal A}$ denote an arbitrary $N \times N$ block matrix with individual block entries of size
$M\times M$ each. Then, the block-maximum norm of ${\cal A}$ is defined as \be \|{\cal A}\|_{b,\infty}
\define \max_{x \neq 0}
                            \frac{\|{\cal A}x\|_{b,\infty}}{\|x\|_{b,\infty}}
\label{defioa.slak}\ee The block-maximum norm has several useful properties --- see \cite{bookchapter}. In particular, when $A$ is $N\times N$ left-stochastic and ${\cal A}=A\otimes I_M$, then it can be verified that
$\|{\cal A}\tran \|_{b,\infty}=1$. Likewise, when ${\cal D}$ is block diagonal and symmetric, then $\|{\cal D}\|_{b,\infty}=\rho({\cal D})$. Consequently, repeating the argument leading to (\ref{llkad.l1k3klakdl}) and replacing the $2-$induced norm used there by the block-maximum norm we have
\be \rho({\cal B})\;\leq\; \|{\cal B}\|_{b,\infty}\;\leq\;\|{\cal A}_2\tran\|_{b,\infty}\cdot \|{\cal D}\|_{b,\infty}\cdot \|{\cal A}_1\tran\|_{b,\infty}\;=\;\rho(D)\label{lkad.l1l3kasax1}\ee
and we again conclude that a stable symmetric ${\cal D}$ guarantees a stable ${\cal B}$ for general left-stochastic matrices $A_1$ and $A_2$ since, for symmetric ${\cal D}$, it holds that $\rho({\cal D})=\|{\cal D}\|_{b,\infty}$. \\

\noindent \underline{\bf Remark 1.} The same argument (\ref{lkad.l1l3kasax1}) can be used to relax the requirement of symmetry and stability on the block diagonal matrix ${\cal D}$. Actually, as long as $\|{\cal D}\|_{b,\infty}<1$, which is guaranteed by requiring the block diagonal entries of ${\cal D}$ to have their $2-$induced norms bounded by one, we can again conclude that $\rho({\cal B})<1$ so that ${\cal B}$ is stable.\\

\noindent \underline{\bf Remark 2.} The validity of property
(\ref{llkad/;l13l;1l3}) for general left-stochastic matrices was already noted and exploited in earlier works, e.g., in  Lemma 1 of \cite{Cattivelli10} and in Lemma 2 of \cite{Cattivelli10b}. However, the statement of Lemma 1 in \cite{Cattivelli10} left out the qualification ``diagonal'' for the center matrix,  and the norm $\|\cdot\|_{\rho}$ that was used in the proof of the lemma in Appendix I of \cite{Cattivelli10} should be replaced by the $\|\cdot\|_{\infty}-$norm. Although these corrections were already noted in the references \cite{bookchapter,zhaoxiao}, we restate below the correct form of Lemma 1 from \cite{Cattivelli10} for accuracy, and provide its adjusted proof using the current notation. That lemma deals with left-stochastic matrices $A_1$ and $A_2$ prior to extension by Kronecker products. Its correct statement should read as follows (the word ``diagonal'' is missing from the statement in  \cite{Cattivelli10}).\\

\noindent {\em Restatement of Lemma 1 \cite{Cattivelli10}}: Let $A_1$, $A_2$, and $D$ denote arbitrary $N\times N$ matrices, where $A_1$ and $A_2$ have real non-negative entries, with columns adding up to one, i.e., $\one\tran A_1=\one\tran,$ $\one\tran A_2=\one\tran$. Then, the matrix $B=A_2\tran D A_1\tran$ is stable for {\em any} choice of $A_1$ and $A_2$ if, and only if, the \underline{diagonal matrix} $D$ is stable.\\

\noindent {\em Proof}: One immediate derivation is to employ the $\infty-$norm and to note that
\be
\|A_1\tran\|_{\infty}=1=\|A_2\tran\|_{\infty},\;\;\;\;\|D\|_{\infty}=\rho(D)
\ee
and, hence,
\be \rho(B)\;\leq\; \|B\|_{\infty}\;\leq\;\|A_2\tran\|_{\infty}\cdot \|D\|_{\infty}\cdot \|A_1\tran\|_{\infty}\;=\;\rho(D)\label{replace.lslad}\ee
It follows that a stable $D$ guarantees a stable $B$. A second derivation that fixes the argument from Appendix I of \cite{Cattivelli10} relies on using the $\infty-$norm instead of the $\rho-$norm used there. Thus, note that
we can alternatively argue that \bq
\|B^i\|_{\infty}&\leq & \left(\|A_2\tran\|_{\infty}\right)^i\cdot \left(\|D\|_{\infty}\right)^i\cdot \left(\|A_1\tran\|_{\infty}\right)^i\nn\\
&=&\left(\rho(D)\right)^i\;\longrightarrow\; 0,\; \mbox{\rm as $i\rightarrow \infty$}
\eq
when $D$ is stable. This completes the proof of Lemma 1 from \cite{Cattivelli10}.\\

\noindent \underline{\bf Remark 3}. These same arguments establish the validity of Lemma 2 from \cite{Cattivelli10b}, which deals with a special case involving a matrix of the form ${\cal B}={\cal A}_2\tran {\cal D}$ using the current notation (the matrix in \cite{Cattivelli10b} is denoted by ${\cal F}={\cal C}\tran{\cal M}$ with ${\cal B}\leftarrow {\cal F}$, ${\cal A}_2\leftarrow {\cal C}$, and ${\cal D}\leftarrow {\cal M}$; moreover, the matrix ${\cal M}$ is now generally non-symmetric but still block diagonal with stable entries of the form ${\cal M}_{kk}=(I-P_kS_k)F$ for some matrices $\{P_k,S_k,F\}$ defined in \cite{Cattivelli10b}. When $\|{\cal M}_{kk}\|_2<1$, these block diagonal entries will have $2-$induced norms smaller than one so that $\|{\cal D}\|_{b,\infty}<1$. Again, the $\rho-$norm used in the proof of Lemma 2 in \cite{Cattivelli10b} should be replaced by the $\|\cdot\|_{b,\infty}$ norm and the argument adjusted as in (\ref{lkad.l1l3kasax1}) or, alternatively, we note that
\be
\|{\cal B}^i\|_{b,\infty}\leq  \left(\|{\cal A}_2\tran\|_{b,\infty}\right)^i\cdot \left(\|{\cal D}\|_{b,\infty}\right)^i\;\longrightarrow\; 0,\; \mbox{\rm as $i\rightarrow \infty$}
\ee

}
\qd
\label{example-stability.condition.aa}
\end{example}
}

\vspace{-0.3cm}
\subsection{Diffusion Networks}
Given the superior stability properties of diffusion strategies for adaptation and learning over networks, we continue our presentation by focusing on this class of algorithms. The results in the previous section focus on MSE networks, which deal with mean-square-error cost functions. We now consider networks with more general costs, $\{J_k(w)\}$, and apply diffusion strategies to seek the global minimizer, $w^o$, of the aggregate cost function, $J^\glob(w)$, defined by \eqref{eqn:aggregatecostdef}. Without loss in generality, we consider ATC diffusion implementations of the form (\ref{eqn:diffusionATCasync}):
\bs
\bq
\bpsi_{k,i} & \weq& \w_{k,i-1} - \bmu_k(i) \wh{\nabla_{w\tran}\,J}_k(\w_{k,i-1}) \label{jahdk.1211a}\\
\w_{k,i} & = &\sum_{\ell\in{\small \bm{\Ncal}_{k,i}}} \a_{\ell k}(i)\;\bpsi_{\ell,i} \label{jahdk.1211a.b}
\eq
\es
Similar conclusions will apply to CTA diffusion implementations.

 {%\small
\begin{assumption}[{\bf Conditions on cost functions}]
\label{assumption:costdistributed}  {\em
The aggregate cost $J^\glob(w)$ in \eqref{eqn:aggregatecostdef} is twice differentiable and satisfies a condition similar to \eqref{eqn:boundedHessian} in Assumption \ref{assumption:costfunctions} for some positive parameters $\nu_d \le \delta_d$. Moreover, all individual costs $\{ J_k(w) \}$ are assumed to be strongly-convex with their global minimizers located at $w^o$, as indicated earlier by (\ref{kajdhdg6713.13}).}

\hfill \qd
\end{assumption}
}

As explained before following (\ref{kajdhdg6713.13}), references \cite{NOW2014,ProcIEEE2014,chen2013} present results on the case in which the individual costs are only convex and need not be strongly convex. These references also discuss the case in which the individual costs need not share minimizers. With each agent $k$ in \eqref{eqn:diffusionATCasync}, we again associate a gradient noise vector:
\be
\label{eqn:gradientnoisedistributedk}
\s_{k,i}(\w_{k,i-1}) \define \wh{\nabla_{w\tran}\,J}_{k}(\w_{k,i-1}) - \nabla_{w\tran}\,J_{k}(\w_{k,i-1})
\ee

\smallskip

{%\small
\begin{assumption}[{\bf Conditions on gradient noise}]
\label{assumption:gradientnoisedistributed} {\em
It is assumed that the first and second-order conditional moments of the gradient noise components satisfy:
\bs
\begin{align}
\Ex \left[ \s_{k,i}(\w_{k,i-1}) \,|\, \Filtration_{i-1} \right] & = 0 \\
\Ex \left[ \s_{k,i}(\w_{k,i-1}) \s_{\ell,i}\tran(\w_{\ell,i-1}) \,|\, \Filtration_{i-1} \right] & = 0, \quad \forall\; k \neq \ell \\
\Ex \left[ \| \s_{k,i}(\w_{k,i-1}) \|^2 \,|\, \Filtration_{i-1} \right] & \le \beta_k^2\, \| \wt{\w}_{k,i-1} \|^2 + \sigma_{s,k}^2
\end{align}
almost surely for some nonnegative scalars $\beta_k^2$ and $\sigma_{s,k}^2$, and where $\Filtration_{i-1}$ represents the collection of all random events generated by the iterates from across all agents, $\{\w_{\ell,j},\,\ell=1,2,\dots,N\}$, up to time $i-1$. Moreover, it is assumed that the limiting covariance matrix of $\s_{k,i}(w^o)$ exists:
\be
R_{s,k} \define \lim_{i \rightarrow \infty} \Ex \left[ \s_{k,i}(w^o) \s_{k,i}\tran(w^o)\,|\,\Filtration_{i-1} \right]
\ee
\es
}

\hfill \qd
\end{assumption}}

We collect the error vectors and gradient noises
from across all agents into $N\times 1$ {\em block} vectors,
whose individual entries are of size $M\times 1$ each:
\bq
\widetilde{\w}_{i}\define \ba{c}\widetilde{\w}_{1,i}\\
\widetilde{\w}_{2,i}\\\vdots\\
\widetilde{\w}_{N,i}
\ea,\;\;
\s_i\define \ba{c}{\s}_{1,i}\\
{\s}_{2,i}\\\vdots\\
{\s}_{N,i}
\ea\label{dkja891.lakdlk}\eq
and where we are dropping the argument $\w_{k,i-1}$ from the $\s_{k,i}(\cdot)$ for compactness of notation. Likewise, we introduce the following $N\times N$ {\em block}
diagonal matrices, whose individual entries are of size $M\times M$
each: \bs\bq \hspace{-0.8cm}\bm{\cal M}_i&\hspace{-0.2cm}=\hspace{-0.2cm}& \mbox{\rm diag}\{\;\bm{\mu}_1(i)
I_{M},\;\bm{\mu}_2(i) I_{M},\;\ldots,\;\bm{\mu}_N(i) I_{M}\;\}\label{m.ladlk.here}\\
 \hspace{-0.8cm}\bm{\cal H}_{i-1}&\hspace{-0.2cm}=\hspace{-0.2cm}&\mbox{\rm diag}\left\{\;\bm{H}_{1,i-1},\;
 \bm{H}_{2,i-1},\;\ldots,\;\bm{H}_{N,i-1}\;\right\}\label{dha7812.kad}\eq
where
\be
\H_{k,i-1}\define \int_{0}^{1}\nabla_{w}^2 J_k(w^{o}-t\widetilde{\w}_{k,i-1})dt\label{expdk.H}
\ee\es
Now, in a manner similar to (\ref{eqn:bHdef}), we can appeal to the mean-value theorem \cite{Poljak87,Rudin76,NOW2014} to note that
\bq
{\nabla_{w\tran} J}_k(\w_{k,i-1})&=&-\H_{k,i-1}\widetilde{\w}_{k,i-1}
\eq
so that the approximate gradient vector can be expressed as:
\be
\widehat{\nabla_{w\tran} J}_k(\w_{k,i-1})
=-\H_{k,i-1}\widetilde{\w}_{k,i-1}+\s_{k,i}(\w_{k,i-1})\label{jad6813laldkoque}
\ee
Subtracting $w^{o}$ from both sides of (\ref{jahdk.1211a})--(\ref{jahdk.1211a.b}), and using (\ref{jad6813laldkoque}),
we find that the network error vector evolves according to the following stochastic recursion:
\bs\be
\addbox{\;\widetilde{\w}_i\;=\;\bm{\cal B}_{i-1}\widetilde{\w}_{i-1}\;+\;\bm{\cal
A}_i\tran\bm{\cal M}_i\s_i\;}
\label{like.lakd.ada.asd.as}\ee
where
\be
\bm{\cal B}_{i-1}\define \bm{\cal A}_i\tran\left(I_{NM}-\bm{\cal M}_i\bm
{\cal H}_{i-1}\right)
\label{inad1387.adlk}\ee\es
Recursion (\ref{like.lakd.ada.asd.as}) describes the evolution of the network error vector for general convex costs, $J_k(w)$, in a manner similar to recursion (\ref{like.lakd.ada.abc.BBB}) in the mean-square-error case. However, recursion (\ref{like.lakd.ada.asd.as}) is more challenging to deal with because of the presence of the random matrix $\bm{\cal H}_{i-1}$; this matrix is replaced by the constant term ${\cal R}$ in the earlier recursion (\ref{like.lakd.ada.abc.BBB}) because that example deals with MSE networks where the individual costs, $J_k(w)$, are quadratic in $w$ and, therefore, their Hessian matrices are constant and independent of $w$. In that case, each matrix $\H_{k,i-1}$ in (\ref{expdk.H}) will evaluate to $2R_{u,k}$ and the matrix $\bm{\cal H}_{i-1}$ in (\ref{inad1387.adlk}) will coincide with the matrix ${\cal R}$ defined by (\ref{cal r.lad}).

The next statement ascertains that sufficiently small step-sizes exist that guarantee the MSE stability of the asynchronous diffusion strategy (\ref{jahdk.1211a})--(\ref{jahdk.1211a.b}) \cite{xiao2013,NOW2014}.

{%\small
\begin{lemma}[{\bf MSE network stability}]
\label{lemma:stabilitydistributed}
Consider an asynchronous network of $N$ interacting agents running the ATC diffusion strategy (\ref{jahdk.1211a})--(\ref{jahdk.1211a.b}). Assume the conditions in Assumptions \ref{assumption:asyncnetwork}, \ref{assumption:costdistributed}, and \ref{assumption:gradientnoisedistributed} hold. Let
\be
\mu_{x,\max}\;=\;\max_{1\leq k\leq N}\;\left\{\mu_{x,k}\right\}\;=\;\max_{1\leq k\leq N}\;\left(\bar{\mu}_k+\frac{\sigma_{\mu,k}^2}{\bar{\mu}_k}\right)
\ee
Then, there exists $\mu_o > 0$ such that for all $\mu_{x,\max} < \mu_o$:
\be
\label{eqn:MSDOmudistributed}
\limsup_{i\rightarrow \infty} \Ex \| \wt{\w}_{k,i} \|^2 = O(\mu_{x,\max})
\ee
%\es
%\hfill \qd
\end{lemma}
}

\bp See App. IV of \cite{xiao2013}.

\ep

\smallskip

Result \eqref{eqn:MSDOmudistributed} shows that the MSD of the network is in the order of $\mu_{x,\max}$. Therefore, sufficiently small step-sizes lead to sufficiently small MSDs. As was the case with the discussion in subsection \ref{subsec:MSEperformanceStandalone}, we can also seek a closed-form expression for the MSD performance of the asynchronous diffusion network and its agents. To do that, we first introduce the analog of Assumption \ref{assumption:smoothness} for the network case.

{%\small
\begin{assumption}[{\bf Smoothness Conditions}]
\label{assumption:smoothnessnetwork}
 {\em The Hessian matrix of the individual cost functions $\{ J_k(w) \}$, and the noise covariance matrices defined for each agent in a manner similar to \eqref{eqn:Rsidef} and denoted by $R_{s,k,i}(w)$, are assumed to be locally Lipschitz continuous in a small neighborhood around $w^o$:
\bs
\begin{align}
\| \nabla_{w}^2\,J_k(w^o + \delta w) - \nabla_{w}^2\,J_k(w^o) \| & \;\le\; \tau_{k,d} \| \delta w \| \\
\| R_{s,k,i}(w^o + \delta w) - R_{s,k,i}(w^o) \| & \;\le\; \tau_{k,s} \| \delta w \|^\kappa
\end{align}
for small perturbation $\| \delta w \| \le r_d$ and for some $\tau_{k,d}, \tau_{k,s} \ge 0$ and $1 \le \kappa \le 2$.
\es
}

\hfill\qd
\end{assumption}
}

{%\small
\begin{lemma}[{\bf Asynchronous network MSD performance}]
\label{lemma:asyncMSDnetwork}
Consider an asynchronous network of $N$ interacting agents running the asynchronous diffusion strategy (\ref{jahdk.1211a})--(\ref{jahdk.1211a.b}). Assume the conditions under Assumptions \ref{assumption:asyncnetwork}, \ref{assumption:costdistributed}, \ref{assumption:gradientnoisedistributed}, and \ref{assumption:primitive} hold. Assume further that the step-size parameter $\mu_{x,\max}$ is sufficiently small to ensure mean-square stability, as already ascertained by Lemma \ref{lemma:stabilitydistributed}. Then,
\bs
\be
\label{eqn:MSDasyncATC}
\mbox{\rm MSD}_{{\rm diff},k}^{\rm asyn} \approx \mbox{\rm MSD}_{{\rm diff, av}}^{\rm asyn} \approx \frac{1}{2} \Tr \left[ \left( \sum_{k=1}^N \bar{\mu}_k\bar{p}_kH_k \right)^{-1} \left( \sum_{k=1}^N (\bar{\mu}_k^2 + \sigma^2_{\mu,k})p_{c,kk}  R_{s,k} \right) \right]
\ee
where  $H_k = \nabla_w^2\, J_k(w^o)$. Moreover, for large enough $i$, the convergence rate towards the above steady-state value is well approximated by the scalar:
\be
\label{eqn:alphaasyncATC}
\alpha_{\rm dist}^{\rm asyn} \weq 1 - 2 \lambdamin\left( \sum_{k=1}^N \bar{\mu}_k\bar{p}_k H_k \right) + O\left(\mu_{x,\max}^{1+1/N^2}\right)
\ee
\es
%\hfill \qd
\end{lemma}
}

\bp See App. XII in \cite{xiao2013AA}.

\ep

{\small
\begin{example}[{\bf Gaussian regression data}] {\rm MSE performance expressions of the form (\ref{eqn:MSDasyncATC}) are accurate to first-order in the step-size parameters, i.e., they are on the order of $O(\bar{\mu}_k)$. The same is true of expression (112) from \cite{ProcIEEE2014} in the synchronous case. There are situations, however, where exact expressions for the MSE performance can be derived for multi-agent networks. We illustrate this possibility here for the case of MSE networks of the type described earlier in Example~\ref{example-9A}. We consider first synchronous networks and comment later on how the results should be adjusted to handle asynchronous behavior.

We refer to the same setting of
Example~9 from \cite{ProcIEEE2014} where we have $N$ agents observing streaming data $\{\d_k(i),\u_{k,i}\}$ that satisfy the regression model:
\bs\be
\d_k(i)\;=\;\u_{k,i}w^o+\v_k(i)\label{jad781kaldfka.new}
\ee
We assume the regression vectors are zero-mean Gaussian-distributed with diagonal covariance matrices denoted by $\Lambda_k$, say,
\be R_{u,k}=\Ex\u_{k,i}\tran\u_{k,i}\define\Lambda_k >0\ee
We further assume that the regression data is temporally white and independent over space so that
\be \Ex\u_{k,i}\tran\u_{\ell,j}=R_{u,k}\delta_{k,\ell}\delta_{i,j}\ee
We also assume that the measurement noise process $\v_k(i)$ is temporally white and independent over space so that
$\Ex\v_k(i)\v_{\ell}(j)=\sigma_{v,k}^2\delta_{k,\ell}\delta_{i,j}$ in terms of the Kronecker delta sequence $\delta_{m,n}$. Likewise, the measurement noise $\v_k(i)$ and the regression data $\u_{\ell,j}$ are assumed to be independent of each other for all $k,\ell,i,j$.

The Gaussian assumption on the regression data is useful because, for this case, fourth-order moments of the regression vectors can be evaluated in closed-form (and such fourth-order moment calculations will arise in the process of determining closed-form expressions for the MSE). Specifically, for any matrix $Q\geq 0$ of size $M\times M$, it holds for independent {\em real-valued} Gaussian regressors that \cite{Sayed08}  (there is a typo in reproducing this expression in the second equation of App. II in \cite{Cattivelli10}):
%the expression there is correct for complex data, and also for real data when $k\neq \ell$)
\be
\Ex\u_{k,i}\tran\u_{k,i} Q\u_{\ell,i}\tran\u_{\ell,i}=\label{kadl89173.12l}R_{u,k}Q R_{u,\ell}+\delta_{k,\ell}\left\{
R_{u,k}\mbox{\rm Tr}\left(QR_{u,k}\right)+R_{u,k}Q R_{u,k}
\right\}
\ee
\es
We assume the network is running a diffusion strategy of the form:
\bs\be
\left\{\begin{array}{l}\bm{\phi}_{k,i-1}=\displaystyle \sum_{\ell\in{\cal N}_k} a_{1,\ell k}\w_{\ell,i-1}\\
\bm{\psi}_{k,i} =   \bm{\phi}_{k,i-1}
+ 2\mu_k \u_{k,i}\tran \left[\d_{k}(i)-\u_{k,i}\bm{\phi}_{k,i-1}\right]\\
\w_{k,i}=\displaystyle \sum_{\ell\in{\cal N}_k} a_{2,\ell k}\bm{\psi}_{\ell,i}
\end{array}\right.
\label{dynahga,da.abc.new}\ee
which includes both the ATC and CTA LMS diffusion forms as special cases. We know from (77b) in \cite{ProcIEEE2014} that the error network vector evolves according to the dynamics
\be
\widetilde{\w}_i\;=\;{\cal A}_2\tran\left(I-{\cal M}\bm
{\cal R}_i\right){\cal A}_1\tran\widetilde{\w}_{i-1}\;-\;{\cal
A}_2\tran{\cal M}\z_i,\;\;i\geq 0
\label{like.lakd.ada.abc.new}\ee
which is defined in terms of the following $N\times 1$ column vectors whose entries are of size $M\times 1$ each:
\be
\z_i\define \ba{c}2\u_{1,i}\tran \v_1(i)\\
2\u_{2,i}\tran \v_2(i)\\\vdots\\
2\u_{N,i}\tran \v_N(i)\ea,\;\;\;\;
\s_i\define \ba{c}\s_{1,i}(\bm{\phi}_{1,i-1})\\\s_{2,i}(\bm{\phi}_{2,i-1})\\\vdots\\
\s_{N,i}(\bm{\phi}_{N,i-1})
\ea\label{defiadk13.adlk.new}\ee
\es
Moreover, \bs\be {\cal A}_1\define
A_1\otimes I_{M},\quad{\cal A}_2\define A_2\otimes I_{M}
\label{defioas.a.abc}\ee while the quantities
\bq \hspace{-0.2cm}{\cal M}\hspace{-0.2cm}&\define&\hspace{-0.2cm} \mbox{\rm diag}\{\;\mu_1
I_{M},\;\mu_2 I_{M},\;\ldots,\;\mu_N I_{M}\;\}\label{m.ladlk.abc}\\
\hspace{-0.2cm} \bm{\cal R}_i\hspace{-0.2cm}&\define&\hspace{-0.2cm}\mbox{\rm diag}\left\{2\u_{1,i}\tran \u_{1,i},
 2\u_{2,i}\tran \u_{2,i},\ldots,2\u_{N,i}\tran\u_{N,i}\;\right\}\\
\hspace{-0.2cm} {\cal R}\hspace{-0.2cm}&\define&\hspace{-0.2cm}\Ex\bm{\cal R}_i\;=\;\mbox{\rm diag}\left\{2R_{u,1},
 2R_{u,2},\ldots,2R_{u,N}\;\right\}\\
\hspace{-0.2cm} {\cal S}\hspace{-0.2cm}&\define&\hspace{-0.2cm}\Ex\z_i\z_i\tran\;=\;\mbox{\rm diag}\left\{4\sigma_{v,1}^2R_{u,1},
 \ldots,4\sigma_{v,N}^2 R_{u,N}\right\}\eq\es
\noindent are $N\times N$ {\em block}
diagonal matrices, whose individual entries are of size $M\times M$
each. Let $\Sigma$ be any $N\times N$ non-negative definite block matrix that we are free to choose, with blocks of size $M\times M$. Computing the $\Sigma-$weigthed squared norm of the error vector in (\ref{like.lakd.ada.abc.new}) under expectation gives (see the derivation leading to (269) in \cite{bookchapter} or (38)--(39) in \cite{Cattivelli10}):
\bs\be
\Ex\|\widetilde{\w}_i\|^2_{\Sigma}\;=\;\Ex\|\widetilde{\w}_{i-1}\|^2_{\Sigma'}\;+\;\Ex \left(\z_i\tran {\cal M} {\cal A}_2\Sigma {\cal A}_2\tran {\cal M}\z_i\right) \label{eq.1klajld.ne}
\ee
where the deterministic weighting matrix $\Sigma'$ is given by:
\be
\Sigma'={\cal A}_1\left\{{\cal A}_2\Sigma{\cal A}_2\tran-{\cal A}_2\Sigma{\cal A}_2\tran {\cal M}{\cal R}-{\cal R}{\cal M}{\cal A}_2\Sigma{\cal A}_2\tran+\right.\left.\Ex\left(\bm{\cal R}_i{\cal M}{\cal A}_2\Sigma{\cal A}_2\tran {\cal M}\bm{\cal R}_i\right)\right\}{\cal A}_1\tran\label{eq.1klajld.ne.2}
\ee
\es

We can evaluate the last expectations in (\ref{eq.1klajld.ne})--(\ref{eq.1klajld.ne.2}) in closed-form. But first we need to introduce a convenient block-vector notation, denoted by $\mbox{\rm bvec}(\cdot)$. Thus, given an $N\times N$ block matrix, with blocks of size $M\times M$ each, say for $N=3$,
\bs \be
X=\ba{ccc}X_{11}&X_{12}&X_{13}\\X_{21}&X_{22}&X_{23}\\X_{31}&X_{32}&X_{33}\ea
\ee
its block vectorization is obtained as follows. We first vectorize each of the block entries and define the column vector
$x_{k\ell}=\mbox{\rm vec}(X_{k\ell})$; this operation stacks the columns of $X_{k\ell}$ on top of each other. Subsequently, the quantity $\mbox{\rm bvec}(X)$ is obtained by stacking the vectors $\{x_{k\ell}\}$ on top of each other:
\be
\mbox{\rm bvec}(X)\define\mbox{\rm col}\{x_{11},x_{21},x_{31}, x_{12}, x_{22}, x_{32}, x_{13}, x_{23}, x_{33}\}
\ee
The following two useful properties can be easily verified for any block matrices $\{A,B,\Sigma\}$ of compatible dimensions \cite{Tracy72,koning}:
\bq
\mbox{\rm bvec}(A\Sigma B)&=&(B\tran\otimes_b A)\mbox{\rm bvec}(\Sigma)\label{prop.kron.1}\\
\mbox{\rm Tr}(A\tran B)&=&\left(\mbox{\rm bvec}(A)\right)\tran \mbox{\rm bvec}(B) \label{prop.kron.2}
\eq
where the notation $A\otimes_b B$ denotes the block Kronecker product of two block matrices $A$ and $B$ (assumed here to be both of size $N\times N$ with $M\times M$ blocks); the $k\ell-$th block of $A\otimes_b B$ has size $N M^2\times N M^2$ and is given by \cite{koning}:
\be
\left[A\otimes_b B\right]_{k\ell}=\ba{cccc}A_{k\ell}\otimes B_{11}&
A_{k\ell}\otimes B_{12}&\ldots&A_{k\ell}\otimes B_{1N}\\
A_{k\ell}\otimes B_{21}&
A_{k\ell}\otimes B_{22}&\ldots&A_{k\ell}\otimes B_{2N}\\
\vdots&\vdots&\ddots&\vdots\\
A_{k\ell}\otimes B_{N1}&
A_{k\ell}\otimes B_{N2}&\ldots&A_{k\ell}\otimes B_{NN}
\ea
\ee
\es
in terms of the traditional Kronecker product operation.  Using the block Kronecker properties (\ref{prop.kron.1})--(\ref{prop.kron.2}) we now find that the last expectation in (\ref{eq.1klajld.ne}) is given by:
\bs\bq
\Ex \left(\z_i\tran {\cal M} {\cal A}_2\Sigma {\cal A}_2\tran {\cal M}\z_i\right)&=&b\tran \sigma
\label{eq.subs.2}\\
\sigma&\define &\mbox{\rm bvec}(\Sigma)\\
b&\define& ({\cal A}_2\tran\otimes_b{\cal A}_2\tran)({\cal M}\otimes_b{\cal M})\mbox{\rm bvec}({\cal S})\label{eqlqlkqe.bb}
\eq
\es
This is the same expression (56) from \cite{Lopes08} for the case of CTA diffusion, and the same
expression (42) from \cite{Cattivelli10}) for CTA and ATC diffusion (except that this latter reference
used the traditional vec(.) and Kronecker notation $\otimes$ instead of bvec(.) and the block Kronecker notation $\otimes_b$).

Let us now evaluate the last expectation in (\ref{eq.1klajld.ne}). Let ${\cal Q}={\cal M}{\cal A}_2\Sigma{\cal A}_2\tran {\cal M}$ so that we can rewrite more compactly:
\bs\be
{\cal K}\define \Ex\left(\bm{\cal R}_i{\cal M}{\cal A}_2\Sigma{\cal A}_2\tran {\cal M}\bm{\cal R}_i\right)\;=\;\Ex\bm{\cal R}_i {\cal Q}\bm{\cal R}_i
\ee
The $k\ell-$th block entry in this matrix is given by
\be
{\cal K}_{k\ell}\;=\;4\;\Ex\u_{k,i}\tran\u_{k,i} {\cal Q}_{k\ell}\u_{\ell,i}\tran\u_{\ell,i}
\ee
in terms of the $k\ell-$th block of ${\cal Q}$. Using property (\ref{kadl89173.12l}) for Gaussian regressors, we get
\be
{\cal K}_{k\ell}=4R_{u,k}{\cal Q}_{k\ell} R_{u,\ell}+4\delta_{k\ell}\left\{
R_{u,k}\mbox{\rm Tr}\left[{\cal Q}_{kk} R_{u,k}\right]+R_{u,k}{\cal Q}_{kk} R_{u,k}
\right\}
\ee
It is clear from the above expression that the matrix ${\cal K}$ has the following general form involving two block diagonal matrices:
\bq
{1\over 4}{\cal K}&=&{\cal R}{\cal Q}{\cal R}+Z_1+Z_2\\
Z_1&\define& \mbox{\rm blkdiag}\left\{R_{u,k} {\cal Q}_{kk} R_{u,k}\right\}\\
Z_2&\define&\mbox{\rm blkdiag}\left\{R_{u,k} \mbox{\rm Tr}\left({\cal Q}_{kk} R_{u,k}\right)\right\}\eq
\es
Introduce the block diagonal matrices (written for $N=3$):
\bs\bq
{\cal L}_1&\define& \ba{ccccccccc}
I_{M^2}\\&0\\&&0\\\hline &&&0\\&&&&I_{M^2}\\&&&&&0\\\hline&&&&&&0\\&&&&&&&0\\&&&&&&&&I_{M^2}
\ea\\
{\cal L}_2&\define&
\ba{ccccccccc}
\lambda_1\lambda_1\tran\\&0\\&&0\\\hline &&&0\\&&&&\lambda_2\lambda_2\tran\\&&&&&0\\\hline&&&&&&0\\&&&&&&&0\\&&&&&&&&\lambda_3\lambda_3\tran
\ea\eq
\es where $\lambda_k=\mbox{\rm vec}(\Lambda_k)$. Then, it can be verified that
\bs\bq
\mbox{\rm bvec}(Z_1)&=&{\cal L}_1({\cal R}\otimes_b{\cal R})\mbox{\rm bvec}({\cal Q})\\
\mbox{\rm bvec}(Z_2)&=&{\cal L}_2\mbox{\rm bvec}({\cal Q})
\eq
\es
Noting that
\be
\mbox{\rm bvec}({\cal Q})=({\cal M}\otimes_b {\cal M})({\cal A}_2\otimes_b{\cal A}_2)\sigma
\ee
we conclude that
\bs\bq
\mbox{\rm bvec}({\cal K})&=&
 4{\cal X}({\cal M}\otimes_b {\cal M})({\cal A}_2\otimes_b{\cal A}_2)\sigma\label{eq.subs.222x}\\
 {\cal X}&\define&
 (I+{\cal L}_1)({\cal R}\otimes_b{\cal R})+{\cal L}_2
\eq
Note that the matrix ${\cal X}$ has the following block-diagonal structure:
\bq
{\cal X}&\define&\mbox{\rm diag}\left\{{\cal X}_1,\,{\cal X}_2,\,\ldots,{\cal X}_N\right\}\\
{\cal X}_k&\define& \mbox{\rm diag}\left\{{\cal X}_k^{(1)},\,{\cal X}_k^{(2)},\,\ldots,{\cal X}_k^{(N)}\right\}\\
{\cal X}_k^{(\ell)}&=&\left\{\begin{array}{ll}\Lambda_{k}\otimes \Lambda_{\ell},&\mbox{\rm when $k\neq \ell$}\\
\lambda_k\lambda_k\tran\;+\;2 \Lambda_k\otimes \Lambda_k,&\mbox{\rm when $k=\ell$}
\end{array}\right.
\eq
\es
which is the same structure derived through equations (57)--(67) in \cite{Lopes08}.  Substituting (\ref{eq.subs.222x}) into (\ref{eq.1klajld.ne.2}) and using again (\ref{prop.kron.1})--(\ref{prop.kron.2}) we find that the block vectorized forms of the weighting matrices $\{\Sigma,\Sigma'\}$ are related via (the expression for ${\cal F}$ below
fixes the typo in Eq.~(44) from \cite{Cattivelli10}; in the derivation in this example we considered the special case in which ${\cal S}_m=I$ in (44) from \cite{Cattivelli10}):
\bs
\bq
\sigma'&=&{\cal F}\sigma\\
{\cal F}&\define&({\cal A}_1\otimes_b{\cal A}_1)\left\{I-({\cal R}{\cal M}\otimes _b I)-(I\otimes_b {\cal R}{\cal M})+4({\cal M}\otimes_b {\cal M}){\cal X}\right\}({\cal A}_2\otimes_b{\cal A}_2)\label{eqs.a.f}
\eq
\es
This is the same expression for ${\cal F}$ in Eq. (69) from \cite{Lopes08} for the case of CTA diffusion (where $A_2=I_N$). This is also the same expression for ${\cal F}$ in Eq. (41) from \cite{Cattivelli10} for CTA and ATC diffusion (except that this latter reference wrote $\otimes$ instead $\otimes_b$). Substituting (\ref{eq.subs.2}) and the above expression for $\sigma'$ into
(\ref{eq.1klajld.ne}), and using the compact notation $\|x\|^2_{\sigma}$ for $\|x\|^2_{\Sigma}$, we rewrite (\ref{eq.1klajld.ne}) in the form
\be
\Ex\|\widetilde{\w}_i\|^2_{\sigma}\;=\;\Ex\|\widetilde{\w}_{i-1}\|^2_{{\cal F}\sigma}\;+\;b\tran \sigma\label{eq.asx1klajld.ne}
\ee
In seady-state, as $i\rightarrow\infty$, the mean-square-error approaches
\be
\lim_{i\rightarrow\infty}\;\Ex\|\widetilde{\w}_i\|^2_{(I-{\cal F})\sigma}\;=\;b\tran \sigma
\ee
As explained in Sec.~6.6 of \cite{bookchapter}, the network MSD can be assessed by selecting $\sigma$ to satisfy
\be
(I-{\cal F})\sigma\;=\;\frac{1}{N}\mbox{\rm bvec}(I_{NM})
\ee
which leads to the desired expression
\be
\mbox{\rm MSD}_{{\rm diff, av}}^{\rm sync}\;=\;\frac{1}{N}b\tran (I-{\cal F})^{-1}\mbox{\rm bvec}(I_{NM})
\ee
which is expression (105a) from \cite{Lopes08}.

To illustrate how these results are adjusted for asynchronous behavior, we consider the case in which the step-size parameters $\{\mu_k\}$ in (\ref{dynahga,da.abc.new}) are replaced by random values $\{\bm{\mu}_k(i)\}$. We denote the mean and variances of these random variables as follows:
\bs
\bq
\bar{\mu}_k&\define&\Ex\mu_k(i)\\
\bm{\cal M}_i&\define&\mbox{\rm diag}\{\bm{\mu}_1(i) I_M,\bm{\mu}_2(i)I_M,\ldots,\bm{\mu}_N(i)I_M\}\\
\bar{\cal M}&\define&\Ex\bm{\cal M}_i\\
{\cal C}_{\mu}&\define&\Ex(\bm{\cal M}_i-\bar{\cal M})\otimes_b(\bm{\cal M}_i-\bar{\cal M})
\eq
\es
If we now repeat the same analysis, expressions (\ref{eq.1klajld.ne})--(\ref{eq.1klajld.ne.2}) are replaced by
\bs\bq
\Ex\|\widetilde{\w}_i\|^2_{\Sigma}\hspace{-0.1cm}&=&\hspace{-0.1cm}\Ex\|\widetilde{\w}_{i-1}\|^2_{\Sigma'}\;+\;\Ex \left(\z_i\tran \bm{\cal M}_i {\cal A}_2\Sigma {\cal A}_2\tran \bm{\cal M}_i\z_i\right)\label{eq.1klajld.ne.222}\\
\Sigma'\hspace{-0.1cm}&=&\hspace{-0.1cm}{\cal A}_1\left\{{\cal A}_2\Sigma{\cal A}_2\tran-{\cal A}_2\Sigma{\cal A}_2\tran \bar{\cal M}{\cal R}-{\cal R}\bar{\cal M}{\cal A}_2\Sigma{\cal A}_2\tran+\Ex\left(\bm{\cal R}_i\bm{\cal M}_i{\cal A}_2\Sigma{\cal A}_2\tran \bm{\cal M}_i\bm{\cal R}_i\right)\right\}{\cal A}_1\tran\nn\\\label{eq.1klajld.ne.2.222}
\eq
\es
We can evaluate the last expectations in (\ref{eq.1klajld.ne.222})--(\ref{eq.1klajld.ne.2.222}) as follows. First we have:
\bs\bq
\Ex \left(\z_i\tran \bm{\cal M}_i {\cal A}_2\Sigma {\cal A}_2\tran \bm{\cal M}_i\z_i\right)&=&
b\tran \sigma
\label{eq.subs.2.222}\eq
where
\bq
b&\define& \Ex\mbox{\rm bvec}\left({\cal A}_2\tran \bm{\cal M}_i{\cal S}\bm{\cal M}_i {\cal A}_2\right)\nn\\
&=&({\cal A}_2\tran\otimes_b{\cal A}_2\tran)\Ex(\bm{\cal M}_i\otimes_b \bm{\cal M}_i)\mbox{\rm bvec}({\cal S})\nn\\
&=&({\cal A}_2\tran\otimes_b{\cal A}_2\tran)({\cal C}_{\mu}+\bar{\cal M}\otimes_b\bar{\cal M})\mbox{\rm bvec}({\cal S})
\eq\es
Second, we have:
\bq
\mbox{\rm bvec}\left\{\Ex\left(\bm{\cal R}_i\bm{\cal M}_i{\cal A}_2\Sigma{\cal A}_2\tran \bm{\cal M}_i\bm{\cal R}_i\right)\right\}&=&
\Ex(\bm{\cal M}_i\otimes_b \bm{\cal M}_i)\Ex\left(\bm{\cal R}_i\otimes_b \bm{\cal R}_i\right)({\cal A}_2\otimes_b {\cal A}_2)\sigma\nn\\
&=&
(C_{\mu}+\bar{\cal M}\otimes_b \bar{\cal M})\Ex\left(\bm{\cal R}_i\otimes_b \bm{\cal R}_i\right)({\cal A}_2\otimes_b {\cal A}_2)\sigma
\eq
so that we replace the quantities ${\cal M}$ and $({\cal M}\otimes_b {\cal M})$ in the expressions (\ref{eqlqlkqe.bb}) and  (\ref{eqs.a.f}) for $b$ and ${\cal F}$ by $\bar{\cal M}$ and $(C_{\mu}+\bar{\cal M}\otimes_b \bar{\cal M})$, respectively.

}
\qd\label{example-new-Gaussian}
\end{example}
}

\section{Concluding Remarks}
This chapter provides an overview of asynchronous strategies for adaptation, learning, and optimization over networks including non-cooperative, centralized, consensus, and diffusion strategies. Particular attention is given to the constant step-size case in order to examine solutions that are able to adapt and learn continuously from streaming data. The presentation complements the results from \cite{ProcIEEE2014,NOW2014}. We introduced a fairly general model for asynchronous behavior that allows for random step-sizes, link failures, random topology variations, and random combination coefficients. We examined the mean-square-error performance and stability properties under asynchronous events and recovered results for  synchronous operation as a special case. The results indicate that asynchronous networks are robust, resilient to failure, and remain mean-square stable for sufficiently small step-sizes.

There are of course several other aspects of distributed strategies that are not covered in this work. Comments on these aspects can be found in \cite{ProcIEEE2014,NOW2014,bookchapter}, including issues related to (a) the noisy exchange of information over links (e.g., \cite{bookchapter,Kar2009} and \cite{TS11b}--\cite{Mat2009}); (b) the use of gossip strategies (e.g., \cite{boyd2006,moura10,Aysal09,Lopes2008aax,dimakis10,shah2009,ror2010}); (c) the exploitation of sparsity constraints (e.g., \cite{Lor13}--\cite{Liu12}); (d) the solution of constrained optimization problems
 (e.g., \cite{Sriv2011,leenedic,Yan2013,theo11,tow2013YYY}); (e) the use of distributed solutions of the recursive least-squares type (e.g., \cite{bookchapter,Xiao06,Cattivelli08a}); (f) the development of distributed state-space solutions (e.g., \cite{Xiao06,Cattivelli08c} and \cite{OlfatiCDC}--\cite{rantzer06}); (g) the study of incremental-based strategies (e.g., \cite{Bertsekas97}--\cite{Predd2009}); (h) the study of distributed solutions under multi-task environments \cite{Roula2016a}--\cite{Chen2015R};  (i) the case of non-smooth risk functions in the context of subgradient learning \cite{Ying2017a,Ying2017b,NedicOzdal,Chou12}; and (j) the incorporation of proximal operators into the distributed setting \cite{Ozdaglar2012}--\cite{Roula2016b}.


\begin{thebibliography}{99}

\bibitem{ProcIEEE2014} A. H. Sayed, ``Adaptive networks,''
{\em Proceedings of the IEEE}, vol. 102, no. 4, pp. 460--497, April 2014.


\bibitem{xiao2013} X. Zhao and A. H. Sayed, ``Asynchronous adaptation and learning over networks --- Part I: Modeling and stability analysis,'' {\em IEEE Trans. Signal Processing}, vol.~63, no.~4, pp.~811--826, Feb. 2015.

\bibitem{xiao2013AA} X. Zhao and A. H. Sayed, ``Asynchronous adaptation and learning over networks -- Part II: Performance analysis,'' {\em IEEE Trans. Signal Processing}, vol.~63, no.~4, pp.~827--842, Feb. 2015.

\bibitem{xiao2013BB} X. Zhao and A. H. Sayed, ``Asynchronous adaptation and learning over networks -- Part III: Comparison analysis,'' {\em IEEE Trans. Signal Processing}, vol.~63, no.~4, pp.~843--858, Feb. 2015.

\bibitem{NOW2014} A. H. Sayed, {\em Adaptation, Learning, and Optimization over Networks},
Foundations and Trends in Machine Learning, vol. 7, issue 4--5, pp. 311--801, NOW Publishers, Boston-Delft, July 2014.

\bibitem{sayedSPM} A. H. Sayed, S-Y. Tu, J. Chen, X. Zhao, and Z. Towfic, ``Diffusion strategies for adaptation and learning over networks,'' {\em IEEE Signal Processing Magazine}, vol. 30, no. 3, pp.~155--171, May 2013.

\bibitem{bookchapter} A. H. Sayed, ``Diffusion adaptation over networks,'' in {\em  E-Reference Signal Processing}, vol. 3, R. Chellapa and S. Theodoridis, {\em Eds.}, pp. 323--454, Academic Press, 2014. Also available as arXiv:1205.4220v1 [cs.MA], May 2012.

\bibitem{barba13} P. Di Lorenzo, S. Barbarossa, and A. H. Sayed, ``Bio-inspired decentralized radio access based on swarming mechanisms over adaptive networks,'' {\em IEEE Trans. Signal Processing}, vol. 61, no. 12, pp. 3183--3197, June 2013.

\bibitem{richard2013} P. Chainais and C. Richard, ``Learning a common dictionary over a sensor network,'' {\em Proc. IEEE CAMSAP}, pp. 1--5, Saint Martin, Dec. 2013.

\bibitem{onlineicassp} J. Chen, A. H. Sayed, and Z. Towfic, ``Dictionary learning over distributed models,'' {\em IEEE Trans. Signal Process.} vol. 63, issue 4, pp.~1001--1016, February 2015.


\bibitem{TsiBB86} J. Tsitsiklis, D. Bertsekas, and M. Athans, ``Distributed asynchronous deterministic and stochastic gradient optimization algorithms,'' {\em IEEE Trans. Autom. Control}, vol. 31, no. 9, pp. 803--812, Sep. 1986.

\bibitem{boyd2006} S. Boyd, A. Ghosh, B. Prabhakar, and D. Shah, ``Randomized gossip algorithms,'' {\em IEEE Trans. on Information Theory,} vol. 52, no. 6, pp. 2508--2530, Jun. 2006.

\bibitem{Sriv2011} K. Srivastava and A. Nedic, ``Distributed asynchronous constrained
stochastic optimization,'' {\em IEEE J. Sel. Topics. Signal Process.}, vol. 5,
no. 4, pp. 772--790, Aug. 2011.

\bibitem{Kar2009} S. Kar and J. M. F. Moura, ``Distributed consensus algorithms in
sensor networks: Link failures and channel noise,'' {\em IEEE Trans. Signal
Process.}, vol. 57, no. 1, pp. 355--369, Jan. 2009.

\bibitem{moura10} S. Kar and J. M. F. Moura, ``Convergence rate analysis of distributed gossip (linear parameter) estimation: Fundamental limits and tradeoffs,'' {\em IEEE Journal on Selected Topics in Signal Processing}, vol. 5, no.~4, pp.~674--690, Aug. 2011.

\bibitem{karmour2008}  S. Kar and J. M. F. Moura, ``Sensor networks with random links: Topology design for distributed consensus,'' {\em IEEE Trans. Signal Process.}, vol.
56, no. 7, pp. 3315--3326, July 2008.

\bibitem{noisemour2010} D. Jakovetic, J. Xavier, and J. M. F. Moura, ``Weight optimization for
consensus algorithms with correlated switching topology,'' {\em IEEE Trans.
Signal Process.}, vol. 58, no. 7, pp. 3788--3801, July 2010.

\bibitem{noisemour2011} D. Jakovetic, J. Xavier, and J. M. F. Moura, ``Cooperative convex
optimization in networked systems: Augmented Lagranian algorithms
with directed Gossip communication,'' {\em IEEE Trans. Signal Process.,}
vol. 59, no. 8, pp. 3889--3902, Aug. 2011.

\bibitem{mourkar2010} S. Kar and J. M. F. Moura, ``Distributed consensus algorithms in sensor
netowrks: Quantized data and random link failures,'' {\em IEEE Trans. Signal
Process.}, vol. 58, no. 3, pp. 1383--1400, Mar. 2010.

\bibitem{Aysal09} T. C. Aysal, M. E. Yildiz, A. D. Sarwate, and A. Scaglione, ``Broadcast gossip algorithms for consensus,'' {\em IEEE Trans. on Signal Processing}, vol. 57, no. 7, pp. 2748--2761, July 2009.

\bibitem{Aysal09Allerton} T.~C. Aysal, A.~D. Sarwate, and A.~G. Dimakis, ``Reaching consensus in wireless networks with probabilistic  broadcast,'' in {\em Proc. Allerton Conf. Commun., Control, Comput.}, Allerton
  House, IL, Sept. and Oct. 2009, pp. 732--739.

\bibitem{Lopes2008aax} C. Lopes and A. H. Sayed, ``Diffusion adaptive networks with changing topologies,'' {\em Proc. IEEE ICASSP}, pp.~3285--3288, Las Vegas, April 2008.

\bibitem{tak2010} N. Takahashi and I. Yamada, ``Link probability control for probabilistic
diffusion least-mean squares over resource-constrained networks,'' in
{\em Proc. IEEE Int. Conf. Acoust., Speech, Signal Process.} (ICASSP), pp. 3518--3521, Dallas, TX, Mar. 2010.

\bibitem{Ying2017a} B. Ying and A. H. Sayed, ``Performance limits of stochastic sub-gradient learning, Part I: Single agent case,'' {\em submitted for publication}, Also available as arXiv1511.07902, April 2017.

\bibitem{Ying2017b} B. Ying and A. H. Sayed, ``Performance limits of stochastic sub-gradient learning, Part II: Multi-agent case,'' {\em submitted for publication}, Also available as arXiv1704.06025, April 2017.

\bibitem{Hay02} S. Haykin, {\em Adaptive Filter Theory}, Prentice Hall, NJ, 2002.

\bibitem{Widrow85} B.~Widrow and S.~D. Stearns, {\em Adaptive Signal Processing}, Prentice Hall, NJ, 1985.

\bibitem{Sayed08} A. H. Sayed, {\em Adaptive Filters}, Wiley, NJ, 2008.

\bibitem{Sayed03} A. H. Sayed, {\em Fundamentals of Adaptive Filtering}, Wiley, NJ, 2003.

\bibitem{Kai00} T. Kailath, A. H. Sayed, and B. Hassibi, {\em Linear Estimation}, Prentice Hall, NJ, 2000.

\bibitem{Bish2007} C. M. Bishop, {\em Pattern Recognition and Machine Learning}, Springer, 2007.

\bibitem{theo2008} S. Theodoridis and K. Koutroumbas, {\em Pattern Recognition}, 4th edition,
 Academic Press, 2008.

\bibitem{hos2000} D. W. Hosmer and S. Lemeshow, {\em Applied Logistic Regression}, 2nd edition, Wiley, NJ, 2000.

\bibitem{Pol73} B. T. Poljak and Y. Z. Tsypkin, ``Pseudogradient adaptation and
training algorithms,'' {\em Autom. Remote Control,} vol. 12, pp. 83--94,
1973.

\bibitem {Poljak87} B. Poljak, {\em Introduction to Optimization}, Optimization Software,
NY, 1987.

\bibitem{Ber97} D. P. Bertsekas and J. N. Tsitsiklis, {\em Parallel and
Distributed Computation: Numerical Methods}, 1st edition, Athena
Scientific, Singapore, 1997.

\bibitem{tsypkin1971} Y.~Z.~Tsypkin, {\em Adaptation and Learning
in Automatic Systems}, Academic Press, NY, 1971.

\bibitem{boyd2004convex} S. Boyd and L. Vandenberghe, {\em Convex
Optimization}, Cambridge University Press, 2004.

\bibitem{dimitri2003} D. Bertsekas, {\em Convex Analysis and Optimization}, Athena Scientific, 2003.

\bibitem{nesterov2004} Y. Nesterov, {\em Introductory Lectures on Convex Optimization: A Basic Course}, Kluwer Academic Publishers, 2004.

\bibitem{Chen10b} J. Chen and A. H. Sayed, ``Diffusion adaptation strategies for
distributed optimization and learning over networks,'' {\em IEEE Trans. Signal Processing},
vol.~60, no.~8, pp.~4289--4305, Aug. 2012.

\bibitem{chen2013} J. Chen and A. H. Sayed, ``On the learning behavior of adaptive networks --- Part I: Transient analysis,'' {\em  IEEE Trans. Information Theory}, vol. 61, no. 6, pp. 3487--3517, June 2015.


\bibitem{widrowhoff1960} B.~Widrow and M.~E.~Hoff, Jr., ``Adaptive
switching circuits,'' {\em IRE WESCON Conv.~Rec.}, Pt.~4,
pp.~96--104, 1960.

\bibitem {berst00} D. P. Bertsekas and J. N. Tsitsiklis, ``Gradient convergence in gradient methods with errors,'' {\em SIAM J. Optim.}, vol. 10, no. 3, pp. 627--642, 2000.

\bibitem{Feller1971} W. Feller, {\em An Introduction to Probability Theory and Its Applications}, vol. 2, Wiley, NY, 1971.

\bibitem{Hahn1994} G. J. Hahn and S. Shapiro, {\em Statistical Models in Engineering}, Wiley, NJ, 1994.

\bibitem{Abra1964} M. Abramowitz and I. Stegun, Eds., {\em Handbook of Mathematical Functions with Formulas, Graphs, and Mathematical Tables}, Dover, NY, 1972.

\bibitem{And1999} G. E. Andrews, R. Askey and R. Roy, {\em Special Functions}, Cambridge University Press, Cambridge, 1999.

\bibitem{Rudin76} W. Rudin, {\em Principles of Mathematical Analysis}, McGraw-Hill, 1976.

\bibitem{yousef2001} N. R. Yousef and A. H. Sayed, ``A unified approach to the steady-state and tracking analysis of adaptive filters,'' {\em IEEE Trans. Signal Processing}, vol. 49, no. 2, pp. 314--324, February 2001.

\bibitem{naffuori} T. Y. Al-Naffouri and A. H. Sayed, ``Transient analysis of data-normalized adaptive filters,''{\em IEEE Trans. Signal Processing}, vol. 51, no. 3, pp.~639--652, Mar. 2003.

\bibitem{chen2013AA} J. Chen and A. H. Sayed, ``On the learning behavior of adaptive networks --- Part II: Performance analysis,'' {\em IEEE Trans. Information Theory}, vol. 61, no. 6, pp. 3518--3548, June 2015.

\bibitem{Papoulis} A. Papoulis and S. U. Pilla, {\em Probability, Random Variables, and Stochastic Processes},  McGraw-Hill, NY, 2002.

\bibitem{Durret} R. Durret, {\em Probability Theory and Examples}, 2nd edition, Duxbury Press, 1996.

\bibitem{Dudley03} R. M. Dudley, {\em Real Analysis and Probability},
2nd edition, Cambridge University Press, 2003.

\bibitem{Widrow76Proc}
B.~Widrow, J.~M. McCool, M.~G. Larimore, and C.~R. Johnson~Jr., ``Stationary
  and nonstationary learning characterisitcs of the {LMS} adaptive filter,''
  \emph{Proc. {IEEE}}, vol.~64, no.~8, pp. 1151--1162, Aug. 1976.

\bibitem{Horowitz81TSP}
L.~Horowitz and K.~Senne, ``Performance advantage of complex {LMS} for
  controlling narrow-band adaptive arrays,'' \emph{{IEEE} Trans. Acoust.,
  Speech, Signal Process.}, vol.~29, no.~3, pp. 722--736, Jun. 1981.

\bibitem{Jones82TIT}
S.~Jones, R.~C. III, and W.~Reed, ``Analysis of error-gradient adaptive linear
  estimators for a class of stationary dependent processes,'' \emph{{IEEE}
  Trans. Inf. Theory}, vol.~28, no.~2, pp. 318--329, Mar. 1982.

\bibitem{Gardner84SP}
W.~A. Gardner, ``Learning characterisitcs of stochastic-gradient-descent
  algorithms: {A} general study, analysis, and critique,'' \emph{Signal
  Process.}, vol.~6, no.~2, pp. 113--133, Apr. 1984.

\bibitem{Feuer85TSP}
A.~Feuer and E.~Weinstein, ``Convergence analysis of {LMS} filters with
  uncorrelated {G}aussian data,'' \emph{{IEEE} Trans. Acoust., Speech, Signal
  Process.}, vol.~33, no.~1, pp. 222--230, Feb. 1985.

\bibitem{Foley88TSP}
J.~B. Foley and F.~M. Boland, ``A note on the convergence analysis of {LMS}
  adaptive filters with {G}aussian data,'' \emph{{IEEE} Trans. Acoust., Speech,
  Signal Process.}, vol.~36, no.~7, pp. 1087--1089, Jul. 1988.


\bibitem{vap2000} V. N. Vapnik, {\em The Nature of Statistical Learning Theory},
Springer, NY, 2000.

\bibitem{tow2012b} Z. Towfic, J. Chen, and A. H. Sayed, ``On the generalization ability of distributed online learners,'' {\em Proc. IEEE Workshop on Machine Learning for Signal Processing} (MLSP), Santander, Spain, pp.~1--6, Sep. 2012.

\bibitem{zhaoxiao2} X. Zhao and A. H. Sayed, ``Performance limits for distributed estimation over LMS adaptive networks,'' {\em IEEE Trans. Signal Processing}, vol. 60, no. 10, pp. 5107--5124, Oct. 2012.

\bibitem{Chen2013} J. Chen and A. H. Sayed, ``Distributed Pareto optimization via diffusion
strategies,'' {\em IEEE J. Selected Topics in Signal Processing}, vol.~7, no.~2, pp.~205--220, April 2013.

\bibitem{Xiao04} L. Xiao and S. Boyd, ``Fast linear iterations for distributed averaging,''
{\em Syst. Control Lett.}, vol. 53, no. 1, pp. 65--78, Sep. 2004.

\bibitem{Tsi84} J. Tsitsiklis and M. Athans, ``Convergence and asymptotic agreement
in distributed decision problems,'' {\em IEEE Trans. Autom. Control}, vol.
29, no. 1, pp. 42--50, Jan. 1984.

\bibitem{Nedic2010} A. Nedic and A. Ozdaglar, ``Cooperative distributed multi-agent
optimization,'' in {\em Convex Optimization in Signal Processing and Communications},
 Y. Eldar and D. Palomar (Eds.), Cambridge University Press, pp. 340-386, 2010.

\bibitem{NedicOzdal} A. Nedic and A. Ozdaglar, ``Distributed subgradient methods for multi-agent optimization,'' {\em IEEE Trans. Autom. Control}, vol.~54, no.~1, pp.~48--61, Jan. 2009.

\bibitem{johansson} B. Johansson, T. Keviczky, M. Johansson, and K. Johansson, ``Subgradient methods and consensus algorithms for solving convex optimization problems,'' {\em Proc. IEEE CDC}, pp. 4185--4190, Cancun, Mexico, December 2008.

\bibitem{dimakis10} A. G. Dimakis, S. Kar, J. M. F. Moura, M. G. Rabbat, and A. Scaglione, ``Gossip algorithms for distributed signal processing,'' {\em Proc. IEEE}, vol. 98, no. 11, pp.~1847--1864, Nov. 2010.

\bibitem{Xiao06} L. Xiao, S. Boyd and S. Lall, ``A space-time diffusion scheme peer-to-peer least-squares-estimation,'' {\em Proc. Information Processing in Sensor Networks (IPSN)}, pp.~168--176, Nashville, TN, April 2006.

\bibitem{Ren05}
W. Ren and R. W. Beard, ``Consensus seeking in multi-agent systems
under dynamically changing interaction topologies,'' \emph{IEEE
Trans. on Automatic Control}, vol. 50, pp. 655--661, May 2005.

\bibitem{Shamma05} R. Olfati-Saber and J. Shamma, ``Consensus filters for sensor networks and distributed sensor fusion,'' {\em Proc. 44th IEEE Conference on Decision and Control (CDC)}, pp.~6698--6703, Seville, Spain, Dec. 2005.

\bibitem{barba07} S. Barbarossa, and G. Scutari, ``Bio-inspired sensor network design,''
{\em IEEE Signal Processing Magazine}, vol. 24, no. 3, pp. 26--35, May 2007.

\bibitem{barba10} S. Sardellitti, M. Giona, and S. Barbarossa, ``Fast distributed average consensus algorithms based on advection-diffusion processes,'' {\em IEEE Trans. Signal Processing}, vol. 58, no. 2, pp. 826--842, Feb. 2010.

\bibitem{braca08} P. Braca, S. Marano, and V. Matta, ``Running consensus in wireless sensor networks,'' {\em Proc. 11th International Conference on Information Fusion}, pp.~1--6, Cologne, Germany, June 2008.

\bibitem{Lopes08} C. G. Lopes and A. H. Sayed,
``Diffusion least-mean squares over adaptive networks: Formulation and performance analysis,'' {\em IEEE Trans. Signal Process.}, vol. 56, no. 7, pp.~3122--3136, July 2008.

\bibitem{Cattivelli10} F. S. Cattivelli and A. H. Sayed, ``Diffusion LMS strategies for distributed estimation,'' {\em IEEE Trans. Signal Process.}, vol. 58, no. 3, pp.~1035--1048, Mar. 2010.

\bibitem{yusayed12bb}  S-Y. Tu and A. H. Sayed, ``Diffusion strategies outperform consensus strategies for distributed estimation over adaptive networks,'' {\em IEEE Trans. on Signal Processing}, vol. 60, no. 12, pp. 6217--6234, Dec. 2012.

\bibitem{horn} R. A. Horn and C. R. Johnson, {\em Matrix Analysis}, Cambridge University Press, 2003.

\bibitem{Tu11} S.-Y. Tu and A. H. Sayed, ``Mobile adaptive networks,'' {\em IEEE J.
Sel. Topics. Signal Process.}, vol. 5, no. 4, pp.~649--664, Aug.
2011.

\bibitem{catmay11} F. Cattivelli and A. H. Sayed, ``Modeling bird flight formations
using diffusion adaptation,'' {\em IEEE Transactions on Signal
Processing}, vol. 59, no. 5, pp.~2038--2051, May 2011.

\bibitem{Dekel8} O. Dekel, R. Gilad-Bachrach, O. Shamir, and L. Xiao, ``Optimal
distributed online prediction,'' {\em Proc. International
Conference on Machine Learning (ICML)}, Bellevue, WA, pp. 713--720, Jun.~2011.

\bibitem{Agarwal8} A. Agarwal and J. Duchi, ``Distributed delayed stochastic optimization,''
{\em Proc. Neural Information Processing Systems (NIPS)}, Granada, Spain, pp. 873--881, Dec. 2011.

\bibitem{Predd06} J. B. Predd, S. B. Kulkarni, and H. V. Poor, ``Distributed learning in wireless sensor networks,'' {\em IEEE Signal Processing
Magazine}, vol. 23, no. 4, pp. 56--69, Jul. 2006.

\bibitem{Tow11} Z. J. Towfic, J. Chen, and A. H. Sayed, ``Collaborative learning of
mixture models using diffusion adaptation,'' {\em Proc. IEEE Workshop
Mach. Learn. Signal Process. (MLSP)}, Beijing, China, pp.~1--6, Sep. 2011.

\bibitem{golub} G. H. Golub and C. F. Van Loan, {\em Matrix
Computations}, 3rd edition, The John Hopkins University Press,
Baltimore, 1996.

\bibitem{BermanPF}
A.~Berman and R.~J. Plemmons, \emph{Nonnegative Matrices in the Mathematical
  Sciences}, SIAM, PA, 1994.

\bibitem{Pillai05} S. U. Pillai, T. Suel, and S. Cha, ``The Perron--Frobenius theorem: Some
of its applications,'' {\em IEEE Signal Process. Mag.}, vol. 22, no. 2, pp. 62--75, Mar.
2005.

\bibitem{Lopes06} C. G. Lopes and A. H. Sayed, ``Distributed processing over
adaptive networks,'' in {\em Proc. Adaptive Sensor Array Processing
Workshop}, MIT Lincoln Laboratory, MA, pp.1--5, June 2006.

\bibitem{Sayed07} A. H. Sayed and C. G. Lopes, ``Adaptive processing over
distributed networks,'' {\em IEICE Trans. Fund. of Electron.,
Commun. and Comput. Sci.}, vol. E90-A, no. 8, pp.~1504--1510, 2007.

\bibitem{Lopes07a} C. G. Lopes and A. H. Sayed, ``Diffusion least-mean-squares over
adaptive networks,'' {\em Proc. IEEE ICASSP}, Honolulu, Hawaii, vol. 3,
pp. 917-920, April 2007.

\bibitem{Lopes07b} C. G. Lopes and A. H. Sayed, ``Steady-state performance of
adaptive diffusion least-mean squares,'' {\em Proc. IEEE Workshop on
Statistical Signal Processing} (SSP), pp. 136-140, Madison, WI,
Aug. 2007.

\bibitem{Cattivelli08}
F. S. Cattivelli and A. H. Sayed, ``Diffusion LMS algorithms with
information exchange,'' {\em Proc. Asilomar Conf. Signals, Syst.
Comput.}, Pacific Grove, CA,  pp.~251--255, Nov. 2008.

\bibitem{Cattivelli07}
 F. S. Cattivelli, C. G. Lopes, and A. H. Sayed, ``A diffusion RLS
scheme for distributed estimation over adaptive networks,''  {\em
Proc. IEEE Workshop on Signal Process. Advances Wireless Comm.
(SPAWC)}, Helsinki, Finland, pp.~1--5, June 2007.


\bibitem{Cattivelli08c}
F. S. Cattivelli and A. H. Sayed, ``Diffusion mechanisms for
fixed-point distributed Kalman smoothing,'' {\em Proc. EUSIPCO},
Lausanne, Switzerland, pp.~1--4, Aug. 2008.

\bibitem{Cattivelli08a}
 F. S. Cattivelli, C. G. Lopes, and A. H. Sayed, ``Diffusion
recursive least-squares for distributed estimation over adaptive
networks,'' {\em IEEE Trans. Signal Process.}, vol. 56, no. 5, pp.~1865--1877, May 2008.

\bibitem{ramdistributed} S. S. Ram, A. Nedic, and V. V. Veeravalli,
``Distributed stochastic subgradient projection algorithms for convex
optimization,'' {\em J. Optim. Theory Appl.}, vol. 147, no. 3, pp.~516--545, 2010.


\bibitem{leenedic} S. Lee and A. Nedic, ``Distributed random projection algorithm for
convex optimization,'' {\em IEEE J. Selected Topics in Signal
Processing}, vol. 7, no. 2, pp. 221--229, Apr. 2013.

\bibitem{bianchi} P. Bianchi, G. Fort, and W. Hachem, ``Performance of a distributed
stochastic approximation algorithm,'' {\em IEEE Trans. Information
Theory}, vol. 59, no. 11, pp. 7405--7418, Nov. 2013.

\bibitem{Stankovic11} S. S. Stankovic, M. S. Stankovic, and D. S. Stipanovic,
``Decentralized parameter estimation by consensus based stochastic
approximation,'' {\em IEEE Trans. on Autom. Control}, vol.~56, no.~3, pp.~531--543, Mar. 2011.

\bibitem{ChenSayeAller} J. Chen and A. H. Sayed, ``On the limiting behavior of distributed optimization strategies,'' {\em Proc. 50th Annual Allerton Conference on Communication, Control, and Computing}, pp.~1535--1542, Monticello, IL, Oct. 2012.


\bibitem{Bondy08}
J.~A. Bondy and U.~S.~R. Murty, {\em Graph Theory}, Springer, 2008.

\bibitem{zhaoxiao} X. Zhao, S-Y. Tu, and A. H. Sayed, ``Diffusion adaptation over networks under imperfect information exchange and non-stationary data,'' {\em IEEE Trans. Signal Processing}, vol. 60, no. 7, pp. 3460--3475, July 2012.

\bibitem{Takahashi10}
N.~Takahashi, I.~Yamada, and A.~H. Sayed, ``Diffusion least-mean-squares with adaptive combiners: Formulation and   performance analysis,'' {\em IEEE Trans. Signal Process.}, vol. 58, no. 9, pp. 4795--4810, Sep. 2010.

\bibitem{Leskovec10JMLR}
J.~Leskovec, D.~Chakrabarti, J.~Kleinberg, C.~Faloutsos, and Z.~Ghahramani, ``Kronecker graphs: {A}n approach to modeling networks,'' vol. 11, pp. 985--1042, Sept. 2010.

\bibitem{TS11b} S-Y. Tu and A. H. Sayed, ``Adaptive networks with noisy links,''
{\em Proc. IEEE Globecom}, pp.~1--5, Houston, TX, December 2011.


\bibitem{abdolee} R. Abdolee and B. Champagne, ``Diffusion LMS algorithms for sensor
networks over non-ideal inter-sensor wireless channels,'' {\em Proc. IEEE Int. Conf. Dist. Comput. Sensor Systems} (DCOSS), pp.~1--6, Barcelona, Spain, June 2011.

\bibitem{khalili} A. Khalili, M. A. Tinati, A. Rastegarnia, and J. A. Chambers,
``Steady state analysis of diffusion LMS adaptive networks with noisy links,'' {\em IEEE Trans. Signal Processing}, vol.~60, no.~2, pp.~974--979, Feb. 2012.

\bibitem{zhaoxiao2abc} X. Zhao and A. H. Sayed, ``Combination weights for diffusion strategies with imperfect information exchange,'' {\em Proc. IEEE ICC}, pp.~648--652, Ottawa, Canada, June 2012.

\bibitem{Mat2009} G. Mateos, I. D. Schizas, and G. B. Giannakis, ``Performance analysis of the consensus-based distributed LMS algorithm,'' {\em EURASIP J. Adv. Signal Process.}, pp. 1--19, 2009, 10.1155/2009/981030, Article ID 981030.

\bibitem{shah2009} D. Shah, ``Gossip algorithms,'' {\em Found. Trends
Netw.,} vol.~3, pp.~1--125, 2009.

\bibitem{ror2010} O. L. Rortveit, J. H. Husoy, and A. H. Sayed, ``Diffusion LMS with communications constraints,'' {\em Proc. 44th Asilomar Conference on Signals, Systems and Computers,} Pacific Grove, CA, pp. 1645--1649, Nov. 2010.

\bibitem{Lor13} P. Di Lorenzo and A. H. Sayed, ``Sparse distributed learning based on diffusion adaptation,''  {\em IEEE Trans. Signal Processing},
    vol. 61, no. 6, pp.~1419--1433, March 2013.

\bibitem{Chou12} S. Chouvardas, K. Slavakis, Y. Kopsinis, S. Theodoridis,
``A sparsity-promoting adaptive algorithm for distributed learning,'' {\em IEEE
Transactions on Signal Processing}, vol. 60, no. 10, pp. 5412--5425, Oct. 2012.

\bibitem{Chou13} S. Chouvardas, G. Mileounis, N. Kalouptsidis, and S. Theodoridis,  ``A greedy sparsity-promoting LMS for distributed adaptive learning in diffusion networks,'' {\em Proc. ICASSP}, pp. 5415--5419, Vancouver, BC, Canada, 2013.

\bibitem{Liu12} Y. Liu, C. Li and Z. Zhang, ``Diffusion sparse least-mean squares over
networks,'' {\em IEEE Transactions on Signal Processing}, vol. 60, no.
8, pp. 4480--4485, Aug. 2012.

\bibitem{Yan2013} F. Yan, S. Sundaram, S. V. N. Vishwanathan, and Y. Qi, ``Distributed
autonomous online learning: Regrets and intrinsic privacy-preserving
properties,'' {\em IEEE Trans. Knowledge and Data Engineering},
vol. 25, no. 11, pp. 2483--2493, Nov. 2013.

\bibitem{theo11} S. Theodoridis, K. Slavakis, and I. Yamada, ``Adaptive learning in a world of projections: A unifying framework for linear and nonlinear classification and regression tasks,'' {\em IEEE Signal Processing Magazine}, vol.~28, no.~1, pp.~97--123, Jan. 2011.

\bibitem{tow2013YYY} Z. Towfic and A. H. Sayed, ``Adaptive penalty-based distributed
stochastic convex optimization,'' {\em IEEE Trans. Signal Process.}, vol. 62, no. 15, pp. 3924--3938, August 2014.

\bibitem{OlfatiCDC} R. Olfati-Saber, ``Kalman-consensus filter: Optimality, stability, and performance,'' {\em Proc. IEEE CDC}, pp.~7036--7042, Shangai, China, 2009.

\bibitem{Saber07} R. Olfati-Saber, ``Distributed Kalman filtering for sensor networks,'' {\em
Proc. 46th IEEE Conf. Decision Control}, pp. 5492--5498, New Orleans, LA, Dec. 2007.

\bibitem{Cattivelli08b}
 F. S. Cattivelli, C. G. Lopes, and A. H. Sayed, ``Diffusion
strategies for distributed Kalman filtering: Formulation and
performance analysis,'' {\em  Proc. IAPR Workshop on Cognitive Inf.
Process.}(CIP), Santorini, Greece,  pp.~36--41, June 2008.

\bibitem{Cattivelli10b} F. Cattivelli and A. H. Sayed, ``Diffusion
strategies for distributed Kalman filtering and smoothing,'' {\em
IEEE Trans. Automatic Control}, vol. 55, no. 9, pp.~2069--2084, Sep. 2010.

\bibitem{asilomar09} F. Cattivelli and A. H. Sayed,
``Diffusion distributed Kalman filtering with adaptive weights,'' {\em
Proc. Asilomar Conference on Signals, Systems and Computers}, pp. 908--912, Pacific Grove, CA, Nov. 2009.

\bibitem{Khan08} U. A. Khan and J. M. F. Moura, ``Distributing the Kalman filter for
large-scale systems,'' {\em IEEE Trans. Signal Processing}, vol. 56, no. 10,
 pp. 4919--4935, Oct. 2008.

\bibitem{rantzer06} P. Alriksson and A. Rantzer, ``Distributed Kalman filtering using weighted
averaging,'' {\em Proc. 17th Int. Symp. Math. Thy Net. Sys} (MTNS), pp.~1--6, Kyoto, Japan, 2006.

\bibitem{Bertsekas97} D. P. Bertsekas, ``A new class of incremental gradient methods
for least squares problems,'' {\em SIAM J. Optim.}, vol. 7, no. 4,
pp.~913--926, 1997.

\bibitem{Bertsekas99} D. P. Bertsekas, {\em Nonlinear Programming}, 2nd edition, Athena Scientific, Belmont, MA, 1999.

\bibitem{Nedic01} A. Nedic and D. P. Bertsekas, ``Incremental
 subgradient methods for nondifferentiable optimization,'' {\em SIAM J. Optim.,}
vol. 12, no. 1, pp.~109--138, 2001.

\bibitem{Rabbat05}  M. G. Rabbat and R. D. Nowak, ``Quantized incremental
algorithms for distributed optimization,'' {\em IEEE J. Sel. Areas
Commun.}, vol. 23, no. 4, pp.~798--808, 2005.

\bibitem{Helou09} E. S. Helou and A. R. De Pierro, ``Incremental subgradients for constrained convex optimization: A unified framework and new methods,'' {\em SIAM J. on Optimization}, vol. 20, pp.~1547--1572, 2009.

\bibitem{Johansson09} B. Johansson, M. Rabi, and M. Johansson, ``A randomized incremental subgradient method for distributed optimization in networked systems,'' {\em SIAM J. on Optimization}, vol. 20, pp.~1157--1170, 2009.

\bibitem{Blatt2008} D. Blatt, A. O. Hero, and H. Gauchman, ``A convergent incremental gradient method with a constant step size,'' {\em SIAM J. Optimization}, vol.~18, pp.~29--51, 2008.

\bibitem{Sayed2006L} A. H. Sayed and C. Lopes, ``Distributed recursive least-squares strategies over adaptive networks,'' {\em Proc. 40th Asilomar Conference on Signals, Systems and Computers}, Pacific Grove, CA, pp. 233--237, Oct.-Nov., 2006.

\bibitem{Lopes07} C. G. Lopes and A. H. Sayed, ``Incremental adaptive
strategies over distributed networks,'' {\em IEEE Trans. Signal
Process.,} vol. 55, no. 8, pp.~4064--4077, Aug. 2007.

\bibitem{Sayed09} A. H. Sayed and F. Cattivelli, ``Distributed adaptive learning
mechanisms,'' {\em Handbook on Array Processing and Sensor
Networks,} S. Haykin and K. J. Ray Liu, {\em Eds.}, pp.~695--722, Wiley,
NJ, 2009.

\bibitem{Li2010} L. Li, J. Chambers, C. G. Lopes, and A. H. Sayed, ``Distributed estimation over an adaptive incremental network based on the affine projection algorithm,'' {\em IEEE Trans. Signal Processing}, vol. 58, no. 1, pp. 151--164, Jan. 2010.

\bibitem{Catt2011} F. Cattivelli and A. H. Sayed, ``Analysis of spatial and incremental LMS processing for distributed estimation,'' {\em IEEE Trans. Signal Processing}, vol. 59, no. 4, pp. 1465-1480, April 2011.

\bibitem{Predd2009} J. B. Predd, S. R. Kulkarni, and H. V. Poor, ``A collaborative training algorithm for distributed learning,'' {\em IEEE Trans. Information Theory}, vol.~55, no. 4, pp.~1856--1871, April 2009.


\bibitem{Roula2016a} R. Nassif, C. Richard, A. Ferrari, and A. H. Sayed, ``Multitask diffusion adaptation over asynchronous networks,'' {\em IEEE Trans. Signal Processing}, vol. 64, no. 11, pp. 2835-2850, June 2016.


\bibitem{Bertrand2010} A. Bertrand and M. Moonen, ``Distributed adaptive node-specific signal estimation in fully connected sensor networks -- Part I: Sequential node updating,'' {\em IEEE Trans. Signal Process.}, vol. 58, no. 10, pp. 5277--5291, Oct. 2010.

\bibitem{Berberidis2014} N. Bogdanovic, J. Plata-Chaves, and K. Berberidis, ``Distributed diffusion-based LMS for node-specific parameter estimation over adaptive networks,'' {\em Proc. IEEE ICASSP}, pp. 7223--7227, Florence, Italy, May 2014.

\bibitem{Chen2014R} J. Chen, C. Richard, and A. H. Sayed, ``Multitask diffusion adaptation over networks,'' {\em IEEE Trans. Signal Processing}, vol. 62, no. 16, pp. 4129--4144, August 2014.

\bibitem{Chen2015R} J. Chen, C. Richard, and A. H. Sayed, ``Diffusion LMS over multitask networks,'' {\em IEEE Trans. Signal Process.}, vol. 63, no. 11, pp. 2733--2748, June 2015.

\bibitem{Ozdaglar2012}  A. I. Chen and A. Ozdaglar, ``A fast distributed proximal gradient
method,'' {\em Proc. Annual Allerton Conference on
Communication, Control, and Computing}, pp. 601--608, Allerton, USA, Oct.
2012.

\bibitem{Wee2013} W. Wee and I. Yamada, ``A proximal splitting approach to regularized
distributed adaptive estimation in diffusion network,'' {\em  in Proc. IEEE
ICASSP},  pp. 5420--5424, Vancouver, Canada, May 2013,

\bibitem{Vlaski2016} S. Vlaski, L. Vandenberghe, and A. H. Sayed, ''Diffusion stochastic optimization with non-smooth regularizers,'' {\em Proc. IEEE ICASSP}, pp. 4149--4153, Shanghai, China, March 2016.

\bibitem{Vlaski2015} S. Vlaski and A. H. Sayed, ``Proximal diffusion for stochastic costs with non-differentiable regularizers,'' {\em Proc. IEEE ICASSP}, pp. 3352--3356, Brisbane, Australia, April 2015.

\bibitem{Roula2016b} R. Nassif, A. Ferrari, C. Richard, and A. H. Sayed, ``Proximal multitask learning over networks with sparsity-inducing coregularization,'' {\em IEEE Trans. Signal Processing}, vol. 64, no. 23, pp. 6329--6344, Dec. 2016.

\bibitem{Tracy72} D. S. Tracy and R. P. Singh, ``A new matrix product and its applications
in partitioned matrix differentiation,'' {\em Statistica Neerlandica}, vol. 26,
no. 4, pp. 143--157, 1972.

\bibitem{koning}  R. H. Koning and H. Neudecker and T. Wansbeek, ``Block Kronecker products and the vecb operator,
{\em Linear Algebra Appl.}, vol. 149, pp. 165--184, Apr. 1991.

\end{thebibliography}
\end{document}